\documentclass[11pt]{amsart}
\usepackage{amssymb,amsfonts,amsmath}
\usepackage[dvips]{graphicx}
\newtheorem{theo}{Theorem}[section]
\newtheorem{prop}[theo]{Proposition}
\newtheorem{lemma}[theo]{Lemma}
\newtheorem{coro}[theo]{Corollary}

\newcommand{\sm}{\setminus}

\newcommand{\eps}{\varepsilon}

\def\proof{\noindent{\bf Proof.}\ }

\def\COMMENT#1{}
\def\noproof{{\unskip\nobreak\hfill\penalty50\hskip2em\hbox{}\nobreak\hfill%
       $\square$\parfillskip=0pt\finalhyphendemerits=0\par}\goodbreak}
\def\endproof{\noproof\bigskip}
\newdimen\margin   
\def\textno#1&#2\par{%
   \margin=\hsize
   \advance\margin by -4\parindent
          \setbox1=\hbox{\sl#1}%
   \ifdim\wd1 < \margin
      $$\box1\eqno#2$$%
   \else
      \bigbreak
      \hbox to \hsize{\indent$\vcenter{\advance\hsize by -3\parindent
      \sl\noindent#1}\hfil#2$}%
      \bigbreak
   \fi}
\def\proof{\removelastskip\penalty55\medskip\noindent{\bf Proof. }}
\def\enddiscard{}
\long\def\discard#1\enddiscard{}

\def\textno#1&#2\par{%
   \margin=\hsize
   \advance\margin by -4\parindent
          \setbox1=\hbox{\sl#1}%
   \ifdim\wd1 < \margin
      $$\box1\eqno#2$$%
   \else
      \bigbreak
      \hbox to \hsize{\indent$\vcenter{\advance\hsize by -3\parindent
      \sl\noindent#1}\hfil#2$}%
      \bigbreak
   \fi}

\title{A proof of Sumner's universal tournament conjecture for large 
tournaments} 
\author{Daniela K\"uhn \and Richard Mycroft \and Deryk Osthus}
\date{}
\thanks{D.~K\"uhn and D.~Osthus were partially supported by the EPSRC, grant
no.~EP/F008406/1.}

\begin{document}

\vspace*{-0.8cm}
\begin{abstract} Sumner's universal tournament conjecture states that any 
tournament on $2n-2$ vertices contains any directed tree on $n$ vertices. In  
this paper we prove that this conjecture holds for all sufficiently large~$n$. 
The proof makes extensive use of results and ideas 
from a recent paper by the same authors, in which an approximate version of
the conjecture was proved.
\end{abstract}

\maketitle

\vspace*{-0.6cm}

\section{Introduction}

\subsection{Introduction}

A tournament is an orientation of a complete graph.
Obviously one cannot guarantee any substructures which contain a cycle within an arbitrary tournament.
On the other hand, Sumner's universal tournament conjecture states that one can find any directed tree $T$ within an arbitrary tournament $G$, 
even if the order of $T$ is rather large compared to that of $G$.
More precisely, the conjecture states that any 
tournament on $2n-2$ vertices contains any directed tree on $n$ vertices.
Many partial results towards this conjecture (made in 1971) have been proved
-- some of them are described below.
Here we prove this conjecture for all large~$n$.
\begin{theo} \label{main}
There exists $n_0$ such that the following holds. Let $T$ be a directed tree 
on $n \geq n_0$ vertices, and $G$ a tournament on $2n-2$ vertices. Then $G$ 
contains a copy of $T$. 
\end{theo}
To see that the bound is best possible, let $T$ be a star with all edges directed 
inwards, and let $G$ be a regular tournament on $2n-3$ vertices. Then 
every vertex of $G$ has $n-2$ inneighbours and $n-2$ outneighbours, and so 
$G$ does not contain a copy of $T$, whose central vertex has $n-1$ 
inneighbours. There are also  `near-extremal' examples which have a different structure
to the one given above:
let $T$ be obtained from a directed path on $\ell\ge 1$ vertices by adding $y:=(n- \ell)/2$
outneighbours to the terminal vertex of the path and $y$ inneighbours to the initial vertex of the path.
Let $G$ consist of regular tournaments $Y$ and $Z$, each on $2y-1$ vertices, together with an
arbitrary tournament $X$ on $\ell-1$ vertices
so that all edges are oriented from $Z$ to $X$, from $X$ to $Y$ and from $Z$ to $Y$.
Then $|G|=2n-\ell-3$ as well as $|T|=n$, and it is easy to see that $G$ does not contain $T$. These examples will play a significant role in the proof (see Section~\ref{proofsketch}).

In~\cite{KMO}, we used a randomised embedding 
algorithm to prove an approximate version of Sumner's universal tournament 
conjecture, and also a stronger result for directed trees of bounded degree. 
Both of these results will be important tools in this paper.

\begin{theo}[\cite{KMO}, Theorem 1.4] \label{approxversion}
Let $\alpha >0$. Then the following properties hold.
\begin{enumerate}
\item[(i)] There exists $n_0$ such that for any $n \geq n_0$, any
tournament $G$ on $2(1+\alpha)n$ vertices contains any directed tree $T$ on $n$ vertices.
\item[(ii)] Let $\Delta$ be any positive integer. Then there exists $n_0$ such that
for any $n \geq n_0$, any tournament $G$ on $(1+\alpha)n$
vertices contains any directed tree $T$ on $n$ vertices with $\Delta(T) \leq \Delta$.
\end{enumerate}
\end{theo}
Let $f(n)$ denote the smallest integer such that any tournament on $f(n)$ 
vertices contains any directed tree on $n$ vertices. So Sumner's 
conjecture states that $f(n) = 2n-2$. Chung (see 
\cite{W}) observed that $f(n) \leq n^{1+o(1)}$, and Wormald~\cite{W} improved 
this to $f(n) \leq O(n \log n)$. The first linear bound on $f(n)$ 
was established by H\"aggkvist and Thomason~\cite{HT}. 
Havet~\cite{H} then showed 
that $f(n) \leq 38n/5$, and later Havet and Thomass\'e~\cite{HTh} used their 
notion of median orders to improve this to $f(n) \leq 7n/2$. Finally El 
Sahili used the same notion to prove the best known bound for general $n$, namely that $f(n) = 3n-3$. 
We shall make extensive 
use of this result in this paper (actually, any linear bound would suffice for 
our purposes; the factor of 3 is not essential.) 

\begin{theo}[El Sahili~\cite{ES}] \label{bestsofar}
Let $T$ be a directed tree on $n$ vertices, and let $G$ be a tournament on 
$3n-3$ vertices. Then $G$ contains a copy of $T$. 
\end{theo}

Sumner's conjecture is also known to hold for special classes of trees (see e.g.~\cite{RW}). 
In particular, Havet and Thomass\'e~\cite{HTh} proved it for `outbranchings',
again using median orders. 
Here an $\emph{outbranching}$ is a directed tree $T$ in which we may choose a root 
vertex $t \in T$ so that for any vertex $t' \in T$, the path between $t$ and 
$t'$ in $T$ is directed from $t$ to $t'$. 
(Outbranchings are also known as arborescences.)

\begin{theo}[Havet and Thomass\'e~\cite{HTh}] \label{outbranchers}
Let $T$ be an outbranching on $n$ vertices, and let $G$ be a tournament on 
$2n-2$ vertices. Then $G$ contains a copy of $T$. 
\end{theo}

For many types of trees, Sumner's conjecture holds with room to spare.
A classical result of this type is Redei's theorem.

\begin{theo}[Redei~\cite{redei}] \label{redeithm}
Any tournament contains a spanning directed path.
\end{theo}
This was generalised considerably by Thomason~\cite{T} who showed that whenever $n$ 
is sufficiently large, every tournament on $n$ vertices contains every orientation 
of the path on $n$ vertices (this was a conjecture of Rosenfeld). Havet and 
Thomass\'e~\cite{HTrosenfeld} proved that this even holds for all $n \neq 3,5,7$. 
They also proposed
the following generalisation of Sumner's conjecture (see~\cite{Havetconj}):
Let $T$ be a directed tree on $n$ vertices with $k$ leaves. Then every tournament on $n+k-1$
vertices contains a copy of $T$.
Some special cases are known (see e.g.~\cite{Havet3leaves}).
It would be interesting to know whether our methods can be used to prove this conjecture.

As illustrated in the next section, 
our proof relies on all of the above theorems (i.e.~Theorems~\ref{approxversion}--\ref{redeithm}),
as well as a directed version of Szemer\'edi's regularity lemma and several structural results proved in~\cite{KMO}.

\subsection{Outline of the proof} \label{proofsketch}

In Section~\ref{sec:defs}, we shall introduce some notation, 
before introducing some key ideas and lemmas. 
In particular we shall define the core tree~$T_\Delta$ of a tree $T$. This is a subtree of $T$ 
consisting of all the `central' vertices of~$T$, which has the important property that every 
component of $T-T_\Delta$ is small. This is useful for the problem of embedding $T$ in a tournament $G$, 
as we may first embed $T_\Delta$ and then proceed to embed the components of $T - T_\Delta$ one 
by one, using the fact that each such component is small. We also introduce the notion of an 
`almost-regular' tournament $G$, which is a tournament in which every vertex has in- and 
outdegree approximately equal to $|G|/2$. 
Section~\ref{sec:defs} also contains three auxiliary lemmas for embedding a directed tree $T$ 
in a tournament $G$ which are derived from Theorems~1.2 and~1.3 and which we shall use extensively in later sections:
\begin{itemize}
\item Lemma~\ref{roundtheback} is designed to embed a directed tree $T$ which is similar to an outstar, 
in the sense that $T$ contains a vertex $t$ with no inneighbours such that every component of 
$T-t$ is small. 
\item In Lemma~\ref{onebyone}, we consider a subtree~$T_c$ of $T$ with the property 
that every component of $T - T_c$ is small, showing that a suitable embedding of~$T_c$ in $G$ can be extended to an embedding of $T$ in $G$. 
\item In Lemma~\ref{component_by_component} we consider the case where the vertices of $G$ 
can be partitioned into disjoint sets~$Y$ and $Z$ such that almost all edges between $Y$ 
and $Z$ are directed the same way. Here we show that if the vertices of $T$ are 
partitioned appropriately between forests $F^-$ and~$F^+$, then to be able to embed 
$T$ in $G$ it is sufficient to embed the largest component of~$F^+$ within~$Y$.
\end{itemize}
We begin the proof of Theorem~\ref{main} in Section~\ref{sec:coretreesize1}, by proving the case where $|T_\Delta|=1$
(Lemma~\ref{core_tree_size_1}). 
Note that the extremal case when $T$ is a star is covered by this case.
To do this, we first embed the single vertex of $T_\Delta$ to a vertex of $G$ with appropriate in- and outdegree. 
We then use Lemma~\ref{roundtheback}, Lemma~\ref{onebyone} and Theorem~\ref{outbranchers} to embed the components of $T - T_\Delta$ 
appropriately among the remaining vertices of $G$ to obtain a copy of~$T$ in $G$. 

Then in Section~\ref{sec:reglem} we introduce the digraph regularity lemma, which yields
a partition of the vertex set of $G$ into clusters so that the edges between pairs of 
clusters of $G$ form quasi-random bipartite subgraphs. We use the regularity lemma to prove
\begin{itemize}
\item Lemma~\ref{ar_tourns_small_core_tree}, which states that Theorem~\ref{main} 
holds in the case where $G$ is almost-regular and~$T_\Delta$ is small enough to be 
embedded within a single cluster of $G$. 
\end{itemize}
To prove this, we first select an appropriate cluster or pair of clusters of $G$ 
in which to embed~$T_\Delta$, and then use Lemma~\ref{onebyone} to extend this embedding 
of $T_\Delta$ to an embedding of $T$ in~$G$. 
We also prove that if we additionally assume that $|T_\Delta| \geq 2$ then the result holds with room to spare,
i.e.~we can allow $G$ to be of order $(2-\alpha) n$, where $\alpha$ is small.

Next, in Section~\ref{sec:robustexpander} we consider the case when the tournament $G$ is a `robust outexpander'.
The latter implies that every set $S$ of reasonable size has a large outneighbourhood.
A key lemma in~\cite{KMO} showed that if $G$ is a robust outexpander tournament on at least $(2+\alpha)n$ 
vertices with large minimum semidegree, then~$G$ contains any directed tree $T$ on $n$ vertices. 
However, the $\alpha n$ error term was only required in the case where $T_\Delta$ is small. 
In Section~\ref{sec:robustexpander} we modify the argument from~\cite{KMO} to prove
\begin{itemize}
\item Lemma~\ref{ar_tourns_large_core_tree}, 
which states that if $T_\Delta$ is large, then any robust outexpander tournament on at least 
$(2-\alpha) n$ vertices with large minimum semidegree contains a copy of $T$. 
\end{itemize}
(The proof relies on further results from~\cite{KMO}.)
It is easy to see that any almost-regular tournament is a robust outexpander tournament.
So we can combine Lemmas~\ref{ar_tourns_small_core_tree} and~\ref{ar_tourns_large_core_tree}
to deduce
\begin{itemize}
\item Lemma~\ref{sumner_for_ar_tourns}, which states that Theorem~\ref{main} holds with a little room to spare 
if $G$ is a large almost-regular tournament and $|T_\Delta|\ge 2$. 
\end{itemize}
We also prepare the ground for the proof of Theorem~\ref{main} by modifying an algorithm from~\cite{KMO} 
to prove Lemma~\ref{tournament_split}. This states that any tournament $G$ may be split into disjoint 
subtournaments, each of which is either small or a robust outexpander with large minimum semidegree. 
This will allow us to apply our results on robust outexpander tournaments to (subtournaments of) 
general tournaments~$G$.

In Section~\ref{sec:smallcore} we prove Lemma~\ref{main_small_core}, which states that Theorem~\ref{main} 
holds for all directed trees $T$ for which $T_\Delta$ is small.
In particular, the `near extremal' construction described in the introduction is dealt with in this part of the proof.
Lemma~\ref{main_small_core} is proved in four steps. 
Firstly, in Lemma~\ref{smallcorepart1} we show that we may assume the tournament $G$ contains two almost-regular 
subtournaments on vertex sets $Y$ and $Z$ which between them contain almost all of the vertices of $G$.
Using this structural information,
we show in Lemmas~\ref{smallcorepart1.5} and~\ref{smallcorepart2} 
that we may assume that $T_\Delta$ is a short directed path and that most of the remainder of $T$ is attached
to the endvertices of this path. 
(Lemma~\ref{sumner_for_ar_tourns} is used as a tool here: we can apply it to embed a suitable subforest of $T$
into $Y$ or $Z$, and afterwards use Lemma~\ref{component_by_component} to embed the remainder of $T$.)
We then consider the case $|T_\Delta| = 2$ separately, proving that Theorem~\ref{main} holds for such $T$. 
This allows us to assume for the proof of Lemma~\ref{main_small_core} that $|T_\Delta| \geq 3$.
Since $T_\Delta$ is a directed path, we can use Redei's theorem to embed $T_\Delta$ within a set $W$
of $|T_\Delta|$  vertices which have high in- and outdegree, and then apply Lemmas~\ref{roundtheback} and~\ref{onebyone}
to complete the embedding again.

Finally, in Section~\ref{sec:mainproof} we complete the proof of Theorem~\ref{main}. 
By Lemma~\ref{main_small_core} we may assume for this that $T_\Delta$ is large. 
None of the extremal or near-extremal cases satisfy this condition, so we will always have a little  room to spare in our calculations
in this part of the proof.
We proceed by using Lemma~\ref{tournament_split} to split the tournament $G$ into disjoint robust 
outexpander subtournaments of large minimum semidegree. If there is just one such subtournament 
then this subtournament contains a copy of $T$ by Lemma~\ref{ar_tourns_large_core_tree}. 
By using Lemma~\ref{component_by_component} we prove Lemma~\ref{largecorepart2}, 
which shows that if there are two such subtournaments then these must also together contain a copy of $T$. 
We may therefore assume in the proof of Theorem~\ref{main} that there are at least three such subtournaments of $G$.
In this case we use Lemma~\ref{ar_tourns_large_core_tree}, Theorem~\ref{approxversion}
and Theorem~\ref{bestsofar} to embed $T$ into these subtournaments.

\section{Definitions and basic tools} \label{sec:defs}

\subsection{Notation}
For a graph $G$, we write $V(G)$ and $E(G)$ to denote the vertex set and edge 
set of $G$ respectively. Then $|G| := |V(G)|$ denotes the number of vertices 
of $G$, and $e(G) := |E(G)|$ is the number of edges of $G$. We shall 
sometimes  write $v \in G$ to mean $v \in V(G)$. A \emph{tree} is a connected 
graph which does not contain any cycles, and we say that a vertex of a tree 
is a \emph{leaf} if it has degree one.  

A \emph{directed graph} $G$, or digraph, consists of a vertex set $V(G)$ and 
an edge set $E(G)$, where each edge $e \in E$ is an ordered pair $(u, v)$ of 
vertices of $G$. For vertices $u, v \in V(G)$ we write $u \rightarrow v$ or 
$v \leftarrow u$ to denote that $(u, v) \in E(G)$. If $u \rightarrow v$ then 
we say that $v$ is an \emph{outneighbour} of $u$, that $u$ is an 
\emph{inneighbour} of $v$, and that the edge $(u, v)$ is \emph{directed from $u$ to $v$}. Sometimes we shall use the term \emph{neighbour} 
of $v$ to mean a vertex which is either an inneighbour or an outneighbour of 
$v$. For any vertex $v \in G$, we denote the set of all outneighbours of $v$ 
by $N^+_G(v)$, or simply $N^+(v)$ when $G$ is clear from the context. 
Similarly we write $N^-_G(v)$ or $N^-(v)$ to denote the set of all 
inneighbours of $v$. Then the \emph{outdegree} of $v$, denoted $d^+_G(v)$, is 
defined by $d^+_G(v) := |N^+_G(v)|$. Similarly the \emph{indegree} of $v$, 
denoted $d^-_G(v)$, is defined by $d^-_G(v) := |N^-_G(v)|$. Again we may write 
$d^+(v)$ or $d^-(v)$ when $G$ is clear from the context. We define the 
minimum outdegree of $G$, denoted $\delta^+(G)$, to be the minimum of 
$d^+(v)$ taken over all vertices $v \in G$, and the minimum indegree, denoted 
$\delta^-(G)$, to be the minimum of $d^-(v)$ taken over all vertices $v \in 
G$. Then the minimum \emph{semidegree} of $G$, denoted $\delta^0(G)$, is the 
minimum of $\delta^-(G)$ and $\delta^+(G)$. We write $G[U \rightarrow V]$ to denote the bipartite subgraph of $G$ formed by edges directed from $U$ to $V$.

We say that a directed graph $G$ is an \emph{oriented graph} if for any $u, v 
\in G$ at most one of $u \rightarrow v$ and $u \leftarrow v$ holds. So an 
oriented graph may be obtained by assigning a direction to each edge of an 
undirected graph. We call this undirected graph the \emph{underlying graph}, 
and denote it by $G_{under}$. An oriented graph is a \emph{tournament} if for 
any distinct $u, v \in V(G)$ precisely one of $u \rightarrow v$ and $u 
\leftarrow v$ holds. Equivalently, the underlying graph of a tournament is a 
complete graph. A \emph{directed tree} is an oriented graph $T$ for which the 
underlying graph $T_{under}$ is a tree. The maximum degree of $T$, denoted 
$\Delta(T)$, is defined to be equal to $\Delta(T_{under})$. A tree or directed 
tree $T$ may be \emph{rooted} by identifying a specific vertex $r$ as the 
\emph{root} of $T$.  

Let $T$ be a directed tree, and let $x$ be a vertex of $T$. Then for any edge 
$e \in E(T)$ incident to $x$, the \emph{weight of $e$ at $x$}, denoted 
$w_e(x)$, is the number of vertices $y$ of $T$ for which $e$ (ignoring the 
orientation) is the first edge of the path in $T_{under}$ from $x$ to $y$. We 
say that a component of $T - x$ is an \emph{incomponent of $x$} if the unique 
edge between $x$ and this component is directed towards $x$, and an \emph{outcomponent 
of $x$} if this edge is directed away from $x$. The \emph{inweight of $x$}, 
denoted $w^-(x)$, is then the number of vertices in incomponents of $x$, and 
the \emph{outweight of $x$}, denoted $w^+(x)$, is the number of vertices in 
outcomponents of $x$. Equivalently, the inweight of $x$ is the sum of 
$w_e(x)$ taken over all edges $e$ incident to $x$ which are directed towards 
$x$, and the outweight can be defined similarly.  

In the same way we define incomponents and outcomponents for a subtree $T_c$ 
of $T$. Indeed, for any component $T'$ of $T - T_c$ there is precisely one 
edge between $T'$ and $T_c$. If this edge is directed towards a vertex of 
$T'$ then we say that $T'$ is an \emph{outcomponent of $T_c$}, whereas if 
this edge is directed towards $T_c$ we say that $T'$ is an \emph{incomponent 
of $T_c$}. As when $T_c$ is a single vertex we define the \emph{inweight of 
$T_c$}, denoted $w^-(T_c)$, to be the number of vertices in incomponents of 
$T_c$, and the \emph{outweight of $T_c$}, denoted $w^+(T_c)$, to be the 
number of vertices in outcomponents of $T_c$. Again these inweights and 
outweights can equivalently be defined as the sum of the weights of the 
appropriate edges of $T$. 

Throughout this paper we shall write $x \ll y$ to indicate that for any $y > 
0$ there exists $x_0 > 0$ such that for any $0 < x \leq x_0$ the subsequent 
statements hold. Such statements with more variables are defined similarly.

\subsection{The core tree}

Let $T$ be a tree on $n$ vertices, and let $\Delta \geq 2$ be fixed. Then we 
say that a vertex $x$ of $T$ is \emph{$\Delta$-core} if every edge $e$ 
incident to $x$ has $w_e(x) \leq (1-1/\Delta)n$. We call the subgraph of $T$ 
induced by $\Delta$-core vertices of $T$ the \emph{core tree of $T$ with 
parameter $\Delta$}, and denote it by $T_\Delta$. With this definition, for 
any tree $T$, the core tree $T_\Delta$ is the same as the $\Delta$-heart of 
$T$ considered by H\"aggkvist and Thomason in~\cite{HT}. The following 
proposition from~\cite{KMO} gives some important properties of the core tree (these properties are also stated in~\cite{HT}). 

\begin{prop}[\cite{KMO}, Proposition 4.2] \label{coretreeprops}
Let $T$ be a tree on $n$ vertices and let $\Delta \geq 2$. Then:
\begin{itemize}
\item[(i)] $T_\Delta$ is a tree containing at least one vertex.
\item[(ii)] $w_e(x) \geq n/\Delta$ if $e = xy$ is an edge of $T_\Delta$.
\item[(iii)] $\Delta(T_\Delta) \leq \Delta$.
\item[(iv)] Every component subtree $T'$ of $T-T_\Delta$ has $|T'| \leq 
n/\Delta$. 
\end{itemize}
\end{prop}

Note that $T_\Delta$ is an undirected tree obtained from an undirected tree 
$T$. However we will frequently refer to the core tree of a directed tree 
$T$; this means the directed tree formed by taking the core tree $T_\Delta$ 
of the underlying graph $T_{under}$ (an undirected tree) of $T$ and directing 
each edge of~$T_\Delta$ as it is directed in $T$. 

The following proposition is needed in the proof of Lemma~\ref{twocoretrees}. 
Essentially the latter states that if trees $T^1$ and $T^2$ almost partition 
a tree $T$, then the core tree $T_\Delta$ is not much larger than $T^1_\Delta 
\cup T^2_\Delta$. 

\begin{prop} \label{deleteleaf}
Let $T$ be a tree on $n$ vertices, let $x$ be a leaf of $T$, and let $\Delta 
\geq 2$. Then $|(T-x)_\Delta| \geq |T_\Delta| - 1$. 
\end{prop}

\proof Let $y$ be a vertex of $T_\Delta - (T-x)_\Delta$, and let $z$ be an 
arbitrary vertex of $(T-x)_\Delta$. Then for some edge $e$ incident to $y$ we 
have $w_e(y) > (1-1/\Delta)(n-1)$ in $T-x$. Since by 
Proposition~\ref{coretreeprops}(iv) the component of $(T-x)-(T-x)_\Delta$ 
containing $y$ contains at most $(n-1)/\Delta$ vertices, this edge must in 
fact be the first edge of the path in $T$ from $y$ to $z$. If $e$ is also the 
first edge of the path in~$T$ from $y$ to $x$ then we have $w_e(y) > 
(1-1/\Delta)(n-1)+1 \geq (1-1/\Delta)n$ in $T$, and so $y \notin T_\Delta$, 
giving a contradiction. So $y$ must lie on the path in $T$ from $x$ to $z$. 
Since $y \in T_\Delta$ we must have $w_e(y) \leq (1-1/\Delta)n$ in $T$, and 
so in $T$ we have  
$$(1-\frac{1}{\Delta})n-1 \leq (1-\frac{1}{\Delta})(n-1) < w_e(y) \leq 
(1-\frac{1}{\Delta})n.$$ 
Clearly this can hold for at most one vertex $y$ on the path from $x$ to $z$. 
So $|T_\Delta - (T-x)_\Delta| \leq 1$, as desired. 
\endproof

\begin{lemma} \label{twocoretrees}
Let $T$ be a tree on $n$ vertices, let $\Delta \geq 2$ and let $\gamma, 
\alpha > 0$. Also let $T^1$ and $T^2$ be subtrees of $T$ such that $|T^1 \cup 
T^2| \geq (1- \gamma) n$. Suppose also that $|T^1_\Delta|, |T^2_\Delta| \leq 
\alpha n$. Then $|T_\Delta| \leq \gamma n + 2 \alpha n + 2n/\Delta$. 
\end{lemma}

\proof
Arbitrarily choose vertices $x_1 \in T^1_\Delta$ and $x_2 \in T^2_\Delta$, 
and let $P$ be the path from $x_1$ to $x_2$ (so $P$ is also a subtree of 
$T$). Then let $T^* := T^1 \cup P \cup T^2$, so $|T^*| \geq (1-\gamma)n$. 
Furthermore, $T^*$ can be formed from $T$ by repeated leaf-deletions. So by 
Proposition~\ref{deleteleaf} we must have $|T| - |T^*| \geq |T_\Delta| - 
|T^*_\Delta|$, and so  
\begin{equation}\label{eq:sizeofcoretree}
|T_\Delta| \leq |T| - |T^*| + |T^*_\Delta| \leq \gamma n - |P - (T^1 \cup 
T^2)| + |T^*_\Delta|. 
\end{equation}

Let $T^*_c := T^1_\Delta \cup P \cup T^2_\Delta$. We claim that $T^*_\Delta 
\subseteq T^*_c$. Indeed, suppose for a contradiction that there exists a 
vertex $y \in T^*_\Delta - T^*_c$. Since $T^*_c$ is a subtree of $T$, every 
vertex of $T^*_c$ lies in the same component $C$ of $T^*-y$. Note that $T^* - 
C$ is a tree. Now, $T^1_\Delta$ and $T^2_\Delta$ are subtrees of $C$, so by 
Proposition~\ref{coretreeprops}(iv) $T^* - C$ contains at most $|T^1|/\Delta$ 
vertices of $T^1$ and at most $|T^2|/\Delta$ vertices of $T^2$. Let $e$ be 
the edge of $T^*$ between $y$ and $C$. Then since $y \in T^*_\Delta$, $w_e(y) 
\leq (1-1/\Delta)|T^*|$ in $T^*$.  So at least $|T^*|/\Delta$ vertices of 
$T^*$ lie in components of $T^*-y$ other than $C$. As every vertex of $P$ 
lies in $C$, either at least $|T^1|/\Delta$ vertices of $T^1$ lie in 
components of $T^*-y$ other than $C$, or at least $|T^2|/\Delta$ vertices of 
$T^2$ lie in components of $T^*-y$ other than $C$. In the former case this 
implies that $T^* - C$ contains more than $|T^1|/\Delta$ vertices of $T^1$, 
and in the latter case this implies that $T^* - C$ contains more than 
$|T^2|/\Delta$ vertices of $T^2$. In either case this yields a contradiction. 

Now, $|T^*_c| \leq 2\alpha n + |P - (T^1_\Delta \cup T^2_\Delta)|$. Since $(P 
\cap T^1) - T^1_\Delta$ is contained within a single component of $T^1 - 
T^1_\Delta$, $|(P \cap T^1) - T^1_\Delta| \leq |T^1|/\Delta$, by 
Proposition~\ref{coretreeprops}(iv). Similarly $|(P \cap T^2) - T^2_\Delta| 
\leq |T^2|/\Delta$. So  
$$|T_\Delta^*| \leq |T_c^*| \leq 2\alpha n + (|T^1|+|T^2|)/\Delta + |P - (T^1 
\cup T^2)|.$$  
So by (\ref{eq:sizeofcoretree})
$$|T_\Delta| \leq \gamma n + (|T^1|+|T^2|)/\Delta + 2\alpha n \leq \gamma n + 
2n/\Delta + 2\alpha n.$$ 
\endproof

\subsection{Almost-regular tournaments.}
In a regular directed graph $G$, every vertex $v$ has $d^+(v) = d^-(v) = 
e(G)/|G|$. We say that a directed graph $G$ is \emph{$\gamma$-almost-regular} 
if every vertex $v \in G$ has $d^+(v), d^-(v) \geq (1 - \gamma) e(G)/|G|$. In 
particular, if $G$ is a tournament then $G$ is $\gamma$-almost-regular if and 
only if every vertex $v \in G$ has $d^+(v), d^-(v) \geq (1 - 
\gamma)(|G|-1)/2$. The next proposition shows that for a large tournament $G$ 
only one of these two bounds is needed to ensure that $G$ contains an 
almost-spanning almost-regular tournament.  

\begin{prop} \label{findartourn}
Suppose that $1/n \ll \alpha \ll \gamma \ll 1$. Let $G$ be a tournament on 
$n$ vertices in which at least one of the following holds: 
\begin{enumerate}
\item[(i)] $d^+(v) \geq (1 - \alpha)(n-1)/2$ for every $v \in G$,
\item[(ii)] $d^-(v) \geq (1 - \alpha)(n-1)/2$ for every $v \in G$, 
\item[(iii)] $d^+(v) \leq (1 + \alpha)(n-1)/2$ for every $v \in G$,
\item[(iv)] $d^-(v) \leq (1 + \alpha)(n-1)/2$ for every $v \in G$.
\end{enumerate}
Then $G$ contains a $\gamma$-almost-regular subtournament $G'$ on at least 
$(1-\gamma)n$ vertices. 
\end{prop}

\proof
We shall prove (i); then (ii), (iii) and (iv) follow immediately.
Suppose that $G$ has at least $\sqrt{\alpha}n$ vertices with $d^+(v) > 
(1+\sqrt{\alpha})(n-1)/2$. Then 
$$\binom{n}{2} = e(G) = \sum_{v \in G} d^+(v) > (1-\alpha)\binom{n}{2} + 
\sqrt{\alpha}n \cdot \sqrt{\alpha}(n-1)/2 = \binom{n}{2},$$ 
giving a contradiction. So there are at most $\sqrt{\alpha}n$ vertices of $G$ 
with $d^+(v) > (1+\sqrt{\alpha})(n-1)/2$. Delete all of these vertices of 
$G$, and let $G'$ be the obtained subtournament. Then $n-\sqrt{\alpha}n \leq 
|G'| \leq n$. Also, every vertex of $G'$ has  
$$d^+_{G'}(v) \geq \frac{(1-\alpha)(n-1)}{2} - \sqrt{\alpha}n \geq 
\frac{(1-\gamma)(|G'|-1)}{2}$$  
and 
$$d^-_{G'}(v) \geq n-1-\sqrt{\alpha}n - \frac{(1+\sqrt{\alpha})(n-1)}{2} \geq 
\frac{(1-\gamma)(|G'|-1)}{2}.$$ 
So $G'$ is a $\gamma$-almost-regular tournament on at least $(1-\gamma)n$ 
vertices, as desired. 
\endproof

\subsection{Some embedding results}

The following three lemmas will be the main tools we shall use to embed directed trees in tournaments. We use Theorem~\ref{bestsofar} in the proofs of all three lemmas, although the factor of 3 in Theorem~\ref{bestsofar} is not critical to our proof; any linear bound would suffice. For the proof of Lemma~\ref{component_by_component} we also require the use of Theorem~\ref{approxversion}. 

\begin{lemma} \label{roundtheback}
Let $T$ be a directed tree on $n$ vertices, rooted at $t$, such that $t$ has 
no inneighbours in~$T$, and every component of $T - t$ contains at most $d$ 
vertices. Let $G$ be a tournament whose vertex set is partitioned into three 
sets, $\{v\}, N$ and $X$, where $|N| \geq n-1$, every vertex of $N$ is an 
outneighbour of $v$, and at least $3d$ vertices of $N$ each have at least 
$6d$ inneighbours in $X$ and at least $6d$ outneighbours in $X$. Then $T$ can 
be embedded in $G$ in such a way that $t$ is embedded to $v$ and at most $4d$ 
vertices of $X$ are occupied by this embedding. 
\end{lemma}

\proof
Let $N' \subseteq N$ consist of all vertices of $N$ with at least $6d$ 
inneighbours in $X$ and at least $6d$ outneighbours in $X$. Then $|N'| \geq 
3d$. We begin by embedding $t$ to the vertex $v$. Now let $T_1, \dots, T_r$ 
be the components of $T - t$, in order of decreasing order. For each $i$, let $t_i$ be the single vertex of $T_i$ which is 
an outneighbour of $t$. Then we shall embed $T_1, \dots, T_r$ in turn in $N 
\cup X$, with each $t_i$ embedded in $N$ and each $T_i$ embedded in the 
vertices not occupied by the embeddings of $T_1, \dots, T_{i-1}$. This will 
give an embedding of $T$ in $G$. So suppose that we have embedded $T_1, 
\dots, T_{i-1}$ in this manner, and we now wish to embed $T_i$. Then at most 
$n-1$ vertices of $T$ have been embedded. At least one of these vertices 
(namely $t$) was not embedded in $N$, so at least one vertex of $N$ must be 
unoccupied. 

Suppose that $N'$ contains at least one unoccupied vertex $v_i$, and also 
that fewer than $3d$ vertices of $X$ have been occupied. Then $v_i$ has at 
least $3d$ unoccupied inneighbours in $X$ and at least $3d$ unoccupied 
outneighbours in $X$. Embed $t_i$ to $v_i$. We then proceed through the 
outcomponents of $t_i$ in $T_i$ in turn. Suppose that when we come to embed 
an outcomponent of $t_i$ we have previously embedded $m$ vertices of $T_i$. 
Then the current outcomponent has order at most $d-m$. Also, $v_i$ has at 
least $3d-m \geq 3(d-m)$ outneighbours in $X$ which have not yet been 
occupied, so by Theorem~\ref{bestsofar} we may embed this outcomponent 
amongst the outneighbours of $v_i$ in $X$. Similarly we may embed the 
incomponents of $t_i$ in turn amongst the inneighbours of $v_i$ in $X$, and 
so we obtain an embedding of $T_i$ in the unoccupied vertices of $G$. Note 
that all vertices of $T_i$ apart from $t_i$ are embedded in $X$. 

Now suppose instead that every vertex of $N'$ has been occupied, but still 
that fewer than $3d$ vertices of $X$ have been occupied. Then at least one of 
the $T_j$ with $j<i$ must have had $|T_j| = 1$, and so $T_i$ consists of one 
single vertex, namely $t_i$. We may therefore embed $t_i$ to any unoccupied 
vertex of $N$ (recall that there is at least one such vertex). 

Finally, suppose that at least $3d$ vertices of $X$ have been occupied. Then 
at least $3d+1$ vertices of $T$ have been embedded outside $N$, and so $N$ 
contains at least $n-1-(n-(3d+1)) = 3d$ unoccupied vertices. Since $|T_i| 
\leq d$, by Theorem~\ref{bestsofar} we may embed $T_i$ among these unoccupied 
vertices. 

By embedding each $T_i$ in this fashion we obtain an embedding of $T$ in $G$ 
with $t$ embedded to $v$. Furthermore, the only vertices embedded in $X$ are 
those in some $T_i$ such that when we came to embed $T_i$, $N'$ contained at 
least one unoccupied vertex $v_i$, and fewer than $3d$ vertices of $X$ had 
been occupied. The embedding of $T_i$ occupied at most another $d$ vertices 
of $X$, and so at most $4d$ vertices of $X$ can have been occupied in total. 
\endproof

\begin{lemma} \label{onebyone}
\begin{itemize}
\item[(a)] Let $T$ be a directed tree, and let $T_c$ be a subtree of $T$ such 
that every component of $T - T_c$ contains at most $d$ vertices. Let $G$ be a 
tournament whose vertices are partitioned into two sets $S$ and $N$ such that 
for every vertex $v \in S$ we have 
\begin{itemize}
\item[(i)] $|N^+(v) \cap N| \geq |T - T_c| + 2d$, and
\item[(ii)] $|N^-(v) \cap N| \geq |T - T_c| + 2d$.
\end{itemize}
Then any embedding of $T_c$ in $G[S]$ can be extended to an embedding of $T$ 
in $G$.  
\item[(b)] Suppose that in addition to the above assumptions we choose a set 
$N' \subseteq N$ and an integer $r \leq |T - T_c|$, so that every vertex $v 
\in S$ satisfies 
\begin{itemize}
\item[(iii)] $|N^+(v) \cap N'| \geq r + 2d$, and
\item[(iv)] $|N^-(v) \cap N'| \geq r + 2d$.
\end{itemize}
Then any embedding of $T_c$ in $G[S]$ can be extended to an embedding of $T$ 
in $G$ such that at least $r$ vertices of $T$ are embedded in $N'$. 
\item[(c)] Suppose that no edges of $T$ are directed from $T_c$ to $T - T_c$. 
Then conditions (i) and (iii) may be dropped without affecting the validity 
of the above result. Likewise if no edges of $T$ are directed from $T - T_c$ 
to $T_c$, then the above results hold even without conditions (ii) and (iv). 
\end{itemize}
\end{lemma}

\proof
Let $n := |T|$. We shall prove (b) and (c); for (a), apply (b) with $r := |T 
- T_c|$ and $N' := N$. Let $T_1, \dots, T_q$ be the components of $T - T_c$, 
so $|T_i| \leq d$ for each $i$. Suppose now that we have successfully 
extended the embedding of $T_c$ in $G[S]$ to an embedding of $T_c \cup T_1 
\cup \dots \cup T_{s-1}$ in $G$. We shall demonstrate how to extend this 
embedding to an embedding of $T_c \cup T_1 \cup \dots \cup T_s$ in $G$. 
Indeed, there is precisely one edge between $T_c$ and $T_s$. Let $t \in T_c$ 
and $t_s \in T_s$ be the endvertices of this edge, and let $v$ be the vertex 
in $S$ to which $t$ is embedded.  

Suppose that $t_s$ is an outneighbour of $t$. By (i), $v$ has at least $|T - 
T_c| + 2d$ outneighbours in $N$. At most $|T_1| + \dots + |T_{s-1}|$ of these 
outneighbours are occupied by the embedding of $T_c \cup T_1 \cup \dots \cup 
T_{s-1}$, and so $v$ has at least $|T_s|+2d \geq 3|T_s|$ outneighbours in $N$ 
which are not occupied by this embedding. Now, by (iii), $v$ has at least $r 
+ 2d$ outneighbours in $N'$. If at most $r-|T_s|$ of these outneighbours are 
occupied by the embedding of $T_c \cup T_1 \cup \dots \cup T_{s-1}$, then by 
Theorem~\ref{bestsofar} we may embed $T_s$ amongst the at least $2d+|T_s| 
\geq 3|T_s|$ unoccupied outneighbours of $v$ in $N'$. If instead $r-k$ of 
these outneighbours are occupied, for some $1 \leq k \leq |T_s|-1$, then by 
Theorem~\ref{bestsofar} we may embed $T_s$ amongst the $2|T_s|+k$ unoccupied 
outneighbours in $N'$ and some arbitrary $|T_s|-k$ outneighbours of $v$ in $N 
\sm N'$. Then at least $k$ vertices of $N'$ will be occupied by this 
embedding of $T_s$. Finally, if at least $r$ outneighbours of $v$ in $N'$ 
have been occupied by this embedding, then we may embed $T_s$ within the at 
least $3|T_s|$ unoccupied outneighbours of $v$ in $N$.  

If instead $t_s$ is an inneighbour of $t$, then we may extend the embedding 
similarly, using (ii) and (iv) rather than (i) and (iii). So we may extend 
the embedding of $T_c$ in $G[S]$ to an embedding of $T$ in $G$ by proceeding 
through each $T_i$ in this manner. Also conditions (i) and (iii) will only be 
required if at least one edge of $T$ is directed from~$T_c$ to $T - T_c$, and 
conditions (ii) and (iv) will only be required if at least one edge of $T$ is 
directed from $T - T_c$ to~$T_c$. Finally, note that after each $T_s$ is 
embedded, either every vertex of $T_1 \cup \dots \cup T_s$ will have been 
embedded in $N'$, or at least $r$ vertices of $T_1 \cup \dots \cup T_s$ will 
have been embedded in $N'$. Since $|T_1 \cup T_2 \cup \dots \cup T_q| = |T - 
T_c| \geq r$, we can be sure that at least $r$ vertices of $N'$ will be 
occupied by the embedding of $T$, as desired.  
\endproof

\begin{lemma} \label{component_by_component}
Suppose that $1/n \ll \gamma \ll \alpha \ll 1$. Let $T$ be a directed tree on $n$ 
vertices, and let forests $F^-$ and $F^+$ be induced subgraphs of $T$ such 
that $V(F^-)$ and $V(F^+)$ partition $V(T)$ and every edge between $F^-$ and 
$F^+$ is directed from $F^-$ to $F^+$. Let $T^+_1$ and $T^+_2$ be the largest 
and second largest components of $F^+$ respectively. Also, let $Y$ and $Z$ be 
disjoint sets such that  
$$|Y| \geq |F^+| + |T^+_2| + \alpha n \textrm{ and } |Z| \geq 2|F^-| + \alpha n.$$ 
Let $G$ be a tournament on vertex set $Y \cup Z$ such that every vertex of 
$Y$ has at most $\gamma n$ outneighbours in $Z$, and every vertex of $Z$ has 
at most $\gamma n$ inneighbours in $Y$. Then any embedding of $T^+_1$ in 
$G[Y]$ can be extended to an embedding of $T$ in $G$. 
\end{lemma}

\proof
Let $T_1, \dots, T_r$ be the components of $F^-$ and $F^+$, ordered so that 
$T_1 = T^+_1$ and so that for each $2 \leq i \leq r$ there is exactly one 
edge of $T$ between $T_i$ and $T_1 \cup \dots \cup T_{i-1}$. Then we have an 
embedding of $T_1$ in $G[Y]$. We shall proceed through the trees $T_i$ in 
turn, embedding each $T_i$ in $G[Y]$ if $T_i$ is a component of~$F^+$, or in 
$G[Z]$ if $T_i$ is a component of $F^-$. Each $T_i$ will be embedded so that 
the embeddings of $T_1, \dots, T_i$ form an embedding of the subtree of $T$ 
induced by the vertices of $T_1, \dots, T_i$. Suppose that we have 
successfully embedded $T_1, \dots, T_{i-1}$ in this manner, and we wish to 
extend this embedding to include $T_i$. Note that there is precisely one edge 
$e$ between $T_i$ and $T_1 \cup \dots \cup T_{i-1}$. Let $t$ be the endvertex 
of $e$ in $T_1 \cup \dots \cup T_{i-1}$, and let $v$ be the vertex to which 
$t$ was embedded.  

If $T_i$ is a component of $F^+$, then $t \in F^-$, so $v \in Z$. In this 
case we will embed $T_i$ within the unoccupied outneighbours of $v$ in $Y$. 
Since $v \in Z$, $|N^+(v) \cap Y| \geq |Y| - \gamma n \geq |F^+| + |T^+_2| + 
\alpha n/2$. At most $|F^+|-|T_i|$ of these vertices are occupied by the 
embeddings of $T_1, \dots, T_{i-1}$. Since $i \geq 2$, $T_i$ is not the 
largest component of $F^+$, and so has order $|T_i| \leq |T_2^+|$. So at 
least $2|T_i| + \alpha n/2$ outneighbours of $v$ in $Y$ remain unoccupied. So 
if $|T_i| \geq \alpha n/2$ then by Theorem~\ref{approxversion}(i) we may 
embed $T_i$ in these unoccupied vertices of $N^+(v) \cap Y$. On the other 
hand, if $|T_i| < \alpha n/2$ then by Theorem~\ref{bestsofar} we may embed 
$T_i$ in these unoccupied vertices of $N^+(v) \cap Y$.  

Now suppose instead that $T_i$ is a component of $F^-$. Then $t \in F^+$, so 
$v \in Y$. Here we will embed $T_i$ within the unoccupied inneighbours of $v$ 
in $Z$. Since $v \in Y$, $|N^-(v) \cap Z| \geq |Z| - \gamma n \geq 2|F^-| + 
\alpha n/2$, and at most $|F^-|-|T_i|$ of these vertices are occupied by the 
embeddings of $T_1, \dots, T_{i-1}$. So at least $2|T_i| +\alpha n/2$ such 
vertices remain unoccupied. So as before, if $|T_i| \geq \alpha n/2$ then by 
Theorem~\ref{approxversion}(i) we may embed $T_i$ in these unoccupied 
vertices of $N^-(v) \cap Z$, whereas if $|T_i| < \alpha n/2$ then by 
Theorem~\ref{bestsofar} we may embed $T_i$ in these unoccupied vertices of 
$N^-(v) \cap Z$. By proceeding through all of the trees $T_i$ in this manner 
we will obtain an embedding of $T$ in $G$. 
\endproof

Observe that if in the statement of Lemma~\ref{component_by_component} we let 
$T^-_1$ and $T^-_2$ be the largest and second-largest components of $F^-$ 
respectively, and replaced the conditions on the sizes of $Z$ and $Y$ by the 
conditions that $|Y| \geq 2|F^+| + \alpha n$ and $|Z| \geq |F^-| + |T_2^-| + 
\alpha n$, then we could conclude that any embedding of $T^-_1$ in $G[Z]$ can 
be extended to an embedding of $T$ in $G$. To see this, either note that the 
proof will still be valid with appropriate changes (switching inneighbours 
and outneighbours and so forth) or observe that this is the effect of 
reversing the direction of every edge of $T$ and every edge of $G$, in which 
case the embedding problem is the same. Sometimes when referring to 
Lemma~\ref{component_by_component} we will implicitly mean this `dual' of 
Lemma~\ref{component_by_component} instead.

\section{Embedding trees whose core tree is a single vertex} \label{sec:coretreesize1}

In this section we shall verify that Sumner's universal tournament conjecture 
holds for large directed trees $T$ whose core tree $T_\Delta$ contains only one 
vertex, that is, trees which are `star-shaped'.  Such trees can be embedded 
by selecting an appropriate vertex to which to embed the single vertex of 
$T_\Delta$, and then embedding the components of $T - T_\Delta$ one by one.  

\begin{lemma}\label{core_tree_size_1}
Suppose that $1/n \ll 1/\Delta \ll 1$. Let $T$ be a directed tree on $n$ vertices with 
$|T_\Delta| = 1$, and let $G$ be a tournament on $2n-2$ vertices. Then $G$ 
contains a copy of $T$. 
\end{lemma}

\proof
Introduce constants $\alpha$ and $\gamma$ with $1/\Delta \ll \alpha \ll 
\gamma \ll 1$. Let $t$ be the single vertex of~$T_\Delta$, let $y$ be the 
outweight of $T_\Delta$, and let $z$ be the inweight of $T_\Delta$. Also, let 
$T_1$ be the subtree of $T$ formed by $t$ and all of its outcomponents, and 
let $T_2$ be the subtree of $T$ formed by $t$ and all of its incomponents. 
Then $y+z = n-1$, $|T_1| = y+1$ and $|T_2| = z+1$. Now, suppose that~$G$ 
contains a vertex $v$ such that  
\begin{itemize}
\item[(i)] either $d^+(v) \geq y + 2n/\Delta$ or $y = 0$, and
\item[(ii)] either $d^-(v) \geq z + 2n/\Delta$ or $z = 0$.
\end{itemize}
Then embed $t$ to $v$. By Proposition~\ref{coretreeprops} each component of 
$T - t$ contains at most $n/\Delta$ vertices. So by Lemma~\ref{onebyone} we 
may extend the embedding of~$t$ in $\{v\}$ to an embedding of $T_1$ in $\{v\} 
\cup N^+(v)$ (since if $y=0$ then $T_1$ consists of the single vertex $t$). 
Also by Lemma~\ref{onebyone}, we may extend the embedding of $t$ in $\{v\}$ 
to an embedding of $T_2$ in $\{v\} \cup N^-(v)$ (since if $z=0$ then $t$ is 
the only vertex of $T_2$). These two embeddings only overlap in the vertex 
$v$, and so combining these two embeddings gives an embedding of $T$ in $G$. 

So we may assume that every vertex $v \in G$ has either $d^+(v) < y 
+2n/\Delta$ or $d^-(v) < z + 2n/\Delta$. Let $Y := \{v \in G: d^+(v) < y 
+2n/\Delta\}$ and let $Z := \{v \in G: d^-(v) < z + 2n/\Delta\}$. Then every 
vertex of $G$ lies in precisely one of $Y$ and $Z$, so $|Y|+|Z|=2n-2$. Thus 
we must have either $|Y| \geq 2y$ or $|Z| \geq 2z$.  Furthermore, if $y = 0$ 
and $|Y| \geq 1$ then each $v \in Y$ has $d^+(v) < 2n/\Delta$ and therefore 
$d^-(v) \geq z + 2n/\Delta$, and so satisfies (ii). We may therefore assume 
that if $y=0$ then $|Y| = 0$ and similarly that if $z = 0$ then $|Z| = 0$. 
So without loss of generality we may assume that $|Y| \geq 2y$ and $y > 0$ 
(otherwise reverse the direction of every edge of $T$ and every edge of $G$; 
then we would have $|Y| \geq 2y$ and $y>0$ at this stage, and the embedding 
problem is the same). Observe that by definition of $Y$ we must also have 
$|Y| \leq 2y + 4n/\Delta + 1$. 

Now suppose that $y \geq \alpha n$. Since $y \in \mathbb{N}$ and $|Y| \geq 
2y$, $Y$ must contain a vertex $v$ which satisfies $|N^+(v) \cap Y| \geq y$. 
Choose a subset $N' \subseteq N^+(v) \cap Y$ of size $y$. For any vertex $u 
\in Y$,  
$$d^+_{G[Y]}(u) = |N^+(u) \cap Y| \leq d_G^+(u) < y+2n/\Delta \leq (1 + 
\alpha) (|Y|-1)/2.$$  
So by Proposition~\ref{findartourn} $G[Y]$ contains a $\gamma$-almost-regular 
tournament on at least $2(1-\gamma)y$ vertices. So at most $|Y|  - 
2(1-\gamma)y \leq 3\gamma y$ vertices of~$Y$ have fewer than $(1-2\gamma)y$ 
inneighbours in~$Y$ or fewer than $(1-2\gamma)y$ outneighbours in $Y$. Since 
$|N'| = y$, at most $6 \gamma y+1$ vertices of~$N'$ have more than 
$(1-3\gamma)y$ inneighbours in $N'$, and at most $6 \gamma y+1$ vertices of 
$N'$ have more than $(1-3\gamma)y$ inneighbours in $N'$. So at least 
$(1-16\gamma)y$ vertices of $N'$ have at least~$\gamma y$ inneighbours in $Y 
\sm N'$ and at least $\gamma y$ outneighbours in $Y \sm N'$. Certainly 
therefore at least~$3n/\Delta$ vertices of $N'$ have at least $6n/\Delta$ 
inneighbours in $Y \sm (\{v\} \cup N')$ and at least $6n/\Delta$ 
outneighbours in $Y \sm (\{v\} \cup N')$. So by Lemma~\ref{roundtheback} we 
may embed $T_1$ in $Y$, with $t$ embedded to $v$, and at most~$4n/\Delta$ 
vertices embedded outside $N' \cup \{v\}$. Let $V'$ be the set of vertices 
of~$G$ not occupied by this embedding of $T_1$. Since $v$ has at least 
$|G|-1-(y+2n/\Delta) \geq z+6n/\Delta$ inneighbours in~$G$, all outside $N' 
\cup \{v\}$, $v$ must have at least $z+2n/\Delta$ unoccupied inneighbours in 
$V'$. So by Lemma~\ref{onebyone} we may extend the embedding of $t$ in 
$\{v\}$ to an embedding of $T_2$ in $\{v\} \cup V'$. These two embeddings 
only overlap in the vertex $v$, and so combine to give an embedding of~$T$ in $G$. 

So we may assume that $1 \leq y < \alpha n$. Then every vertex $v \in Y$ has 
\begin{equation} \label{eq:coretreesize1} 
d^-(v) \geq |G| - 1 - y - 2 n/\Delta \geq n + 2n/\Delta.
\end{equation}
Let $T_3$ be the subtree of $T$ formed by every vertex $t' \in T$ for which 
$T$ contains a directed path from from $t$ to $t'$. Then $t \in T_3$, and 
(taking $t$ as the root vertex) $T_3$ is an outbranching. Also $T_3 \subseteq 
T_1$, so $|T_3| \leq y+1$, and so by Theorem~\ref{outbranchers}, we may embed 
$T_3$ in $G[Y]$. Since $T_\Delta \subseteq T_3$, by 
Proposition~\ref{coretreeprops}(iv) each component of $T - T_3$ contains at 
most $n/\Delta$ vertices. So as every edge of $T$ between $T - T_3$ and $T_3$ 
is directed from $T - T_3$ to $T_3$, and also since 
by~(\ref{eq:coretreesize1}) every vertex of~$Y$ has at least $|T - T_3| + 
2n/\Delta$ inneighbours which were not occupied by the embedding of $T_3$, we 
may extend the embedding of $T_3$ in $G[Y]$ to an embedding of $T$ in~$G$ by 
Lemma~\ref{onebyone}. 
\endproof

\section{The regularity lemma and its applications to embedding trees} \label{sec:reglem}

In this section we shall present a degree form of the regularity lemma for 
directed graphs, and show how this may be used to embed trees. In particular, 
the regularity lemma is useful for embedding directed trees $T$ for which 
$T_{\Delta}$ is substantially smaller than the size of a cluster obtained by 
applying the regularity lemma to a tournament $G$; our approach here is 
essentially to select an appropriate cluster in $G$ in which to embed 
$T_{\Delta}$ so that we may then embed the components of $T - T_{\Delta}$ 
in the remaining clusters of $G$. By using this method we shall prove~Lemma~\ref{ar_tourns_small_core_tree}, which states that Theorem~\ref{main} holds in the case where $G$ is a 
large and almost-regular tournament, and $T$ is a directed tree such that 
$T_{\Delta}$ is small. 

Let $U$ and $V$ be disjoint sets, and let $G$ be a directed graph on vertex 
set $U \cup V$. Recall that $G[U \rightarrow V]$ denotes the bipartite 
subgraph of $G$ formed by edges directed from $U$ to $V$. The \emph{density 
from $U$ to $V$}, denoted $d(G[U \rightarrow V])$, is then defined by $$d(G[U 
\rightarrow V]) := \frac{e(G[U \rightarrow V])}{|U||V|}.$$ We say that $G[U 
\rightarrow V]$ is \emph{$\eps$-regular} if for any $U' \subseteq U$ and $V' 
\subseteq V$ with $|U'| > \eps |U|$ and $|V'| > \eps |V|$ we have $d(G[U' 
\rightarrow V']) = d(G[U \rightarrow V]) \pm \eps$.  

The next lemma is the degree form of the regularity lemma which we shall use. 
A regularity lemma for digraphs was proven by Alon and Shapira~\cite{AS}. The 
degree form follows from this in the same way as in the undirected case 
(see~\cite{BCCsurvey} for a sketch of the latter). 

\begin{lemma}[Regularity Lemma for directed graphs] \label{regularity_lemma}
Suppose that $1/n \ll 1/M \ll 1/M' \ll \eps$. Let~$G$ be a directed graph on 
$n$ vertices. Then there exists a partition of $V(G)$ into $V_0, \dots, V_k$ 
and a spanning subgraph $G'$ of $G$ such that 
\begin{itemize}
\item[(1)] $M' \leq k \leq M$,
\item[(2)] $|V_0| \leq \eps n$,
\item[(3)] $|V_1| = \dots = |V_k|$,
\item[(4)] $d^+_{G'}(x) > d^+_{G}(x) - \eps n$ for all vertices $x \in V(G)$, 
\item[(5)] $d^-_{G'}(x) > d^-_{G}(x) - \eps n$ for all vertices $x \in V(G)$, 
\item[(6)] for all $i \in [k]$ the directed graph $G'[V_i]$ is empty,
\item[(7)] for all $i, j \in [k]$ with $i \neq j$ the directed graph $G'[V_i
\rightarrow V_j]$ is $\eps$-regular.
\end{itemize}
\end{lemma}

We say that an oriented graph $G$ on clusters $V_1, \dots, V_k$ of equal size 
is an \emph{$\eps$-regular cluster tournament} if for any $i, j \in [k]$ with 
$i \neq j$ the subdigraph $G[V_i \rightarrow V_j]$ is $\eps$-regular and for 
any $i \in [k]$ the subdigraph $G[V_i]$ is a tournament. If $G$ is a cluster 
tournament on clusters $V_1, \dots, V_k$ then we shall denote the density of 
$G[V_i \rightarrow V_j]$ by $d_{ij}$ for any $i, j \in [k]$ (the tournament 
$G$ will be clear from the context). The following corollary of the 
regularity lemma shows that any sufficiently large tournament $G$ contains an 
almost-spanning $\eps$-regular cluster tournament $G^*$ such that vertices 
have similar in- and outdegrees in both $G$ and $G^*$. 

\begin{coro} \label{clustertourns}
Suppose that $1/n \ll 1/M \ll 1/M' \ll \eps$. Let $G$ be a tournament on $n$ 
vertices. Then there exist disjoint subsets $V_1, \dots, V_k \subseteq V(G)$ 
of equal size and a subgraph $G^* \subseteq G$ on vertex set $V_1 \cup \dots 
\cup V_k$ such that: 
\begin{itemize}
\item[(i)] $M' \leq k \leq M$,
\item[(ii)] $G^*$ is an $\eps$-regular cluster tournament,
\item[(iii)] $\bigcup_{i \in [k]} V_i \geq (1-\eps)n$,
\item[(iv)] $d^+_{G^*}(x) > d^+_{G}(x) - 2\eps n$ for all vertices $x \in 
V(G)$, and  
\item[(v)] $d^-_{G^*}(x) > d^-_{G}(x) - 2\eps n$ for all vertices $x \in 
V(G)$. 
\end{itemize}
\end{coro}

\proof
Apply Lemma~\ref{regularity_lemma} to obtain a partition $V_0, \dots, V_k$ of 
$V(G)$ and a subgraph $G' \subseteq G$ which satisfy the conditions of 
Lemma~\ref{regularity_lemma}. In particular (i) and (iii) are satisfied. Now 
form $G^*$ from $G'[V_1 \cup \dots \cup V_k]$ by adding every edge of $G$ for 
which both endvertices lie in the same cluster $V_i$. So $G^* \subseteq G$, 
and by (7) of Lemma~\ref{regularity_lemma} and the fact that $G^*[V_i]$ is a 
tournament for each $i \in [k]$ we have (ii). Finally note that using (4) of 
Lemma~\ref{regularity_lemma} we have  
$$d^+_{G^*}(x) \geq d^+_{G'}(x) - |V_0| \geq d^+_{G}(x) - 2\eps n.$$
Similarly $d^-_{G^*}(x) \geq d^-_{G}(x) - 2\eps n$ using (5) of 
Lemma~\ref{regularity_lemma}. 
\endproof

It follows immediately from the definition of regularity that if $U$ and $V$ 
are sets of size $m$, and $G[U \rightarrow V]$ is $\eps$-regular with density 
$d$, then all but at most $2\eps m$ vertices of $U$ have $(d \pm \eps)m$ 
outneighbours in $V$. The next lemma is a generalisation of this fact, 
considering the number of outneighbours of vertices in one cluster within a 
cluster tournament. 

\begin{lemma} \label{regularnbrhoods}
Suppose that $1/m \ll 1/k \ll \eps \ll \eps' \ll 1$.
Let $G$ be an $\eps$-regular cluster tournament on clusters $V_1, \dots, 
V_k$, each of size $m$. Let $V'_j \subseteq V_j$ for each $j \in [k]$ be 
fixed. Then for any~$i$, all but at most $\eps' m$ vertices of $V_i$ have 
$\sum_{j \in [k] \sm \{i\}} d_{ij} |V_j'| \pm \eps' km$ outneighbours in 
$\bigcup_{j \in [k] \sm \{i\}} V_j'$ and $\sum_{j \in [k] \sm \{i\}} d_{ji} 
|V_j'| \pm \eps' km$ inneighbours in $\bigcup_{j \in [k] \sm \{i\}} V_j'$. 
\end{lemma}

\proof Fix some $i \in [k]$. Then let $L$ be the set of all $j \in [k] \sm 
\{i\}$ such that $|V_j'| \geq \eps m$ and $d_{ij} \geq \sqrt{\eps}$. For each 
$j \in L$, let $A_j$ denote the set of vertices of $V_i$ which have fewer 
than $(1-\sqrt{\eps}) d_{ij} |V_j'|$ outneighbours in $V_j'$. 
Then for each $j \in L$, the subdigraph of $G[V_i \rightarrow V_j]$ induced 
by $A_j$ and $V_j'$ has density less than $(1-\sqrt{\eps}) d_{ij} \leq 
d_{ij}-\eps$. Since $G[V_i \rightarrow V_j]$ is $\eps$-regular with density 
$d_{ij}$, and $|V_j'| \geq \eps m$, we must have $|A_j| < \eps m$. 

Now, fix a vertex $v \in V_i$. Suppose that $v$ appears in at most 
$\sqrt{\eps}|L|$ of the sets $A_j$ with $j \in L$. Then 
\begin{align*}
|N^+(v) \cap \bigcup_{j \in L} V_j'| 
&\geq \sum_{j \in L: v \notin A_j}(1-\sqrt{\eps}) d_{ij}|V_j'| 
\\& \geq \sum_{j \in [k] \sm \{i\}}(1-\sqrt{\eps}) d_{ij}|V_j'| - \sum_{j \in 
[k] \sm (L \cup \{i\})} d_{ij}|V_j'|- \sum_{j \in L: v \in A_j} d_{ij}|V_j'| 
\\ &\geq \sum_{j \in [k] \sm \{i\}} d_{ij}|V_j'| - \sqrt{\eps} km - 
\sqrt{\eps} km - \sqrt{\eps}|L|m  
\\ & \geq \sum_{j \in [k] \sm \{i\}} d_{ij}|V_j'|- 3\sqrt{\eps} km.
\end{align*}

Since at most $\sqrt{\eps}m$ vertices $v \in V_i$ appear in more than 
$\sqrt{\eps}|L|$ of the sets $A_j$ with $j \in L$, we may conclude that there 
are at most $\sqrt{\eps}m$ vertices $v \in V_i$ with fewer than $\sum_{j \in 
[k] \sm \{i\}} d_{ij}|V_j'|- 3\sqrt{\eps} km$ outneighbours in $\bigcup_{j 
\in [k] \sm \{i\}} V_j'$. A similar argument shows that there are at most 
$\sqrt{\eps}m$ vertices $v \in V_i$ with more than $\sum_{j \in [k] \sm 
\{i\}} d_{ij}|V_j'| + 3\sqrt{\eps} km$ outneighbours in $\bigcup_{j \in [k] 
\sm \{i\}} V_j'$. 

Now, let $L'$ be the set of all $j \in [k]$ such that $|V_j'| \geq \eps m$ 
and $d_{ji} \geq \sqrt{\eps}$. Then the same argument applied to inneighbours 
rather than outneighbours shows that there are at most $\sqrt{\eps}m$ 
vertices $v \in V_i$ with fewer than $\sum_{j \in [k] \sm \{i\}} 
d_{ji}|V_j'|- 3\sqrt{\eps} km$ inneighbours in $\bigcup_{j \in [k] \sm \{i\}} 
V_j'$ and at most $\sqrt{\eps}m$ vertices $v \in V_i$ with more than $\sum_{j 
\in [k] \sm \{i\}} d_{ji}|V_j'| + 3\sqrt{\eps} km$ inneighbours in 
$\bigcup_{j \in [k] \sm \{i\}} V_j'$. Since $\eps \ll \eps'$, this completes 
the proof. 
\endproof

The next two lemmas will be used in the proof of 
Lemma~\ref{ar_tourns_small_core_tree}; we state them separately as we shall 
also refer to them in Section~\ref{sec:smallcore}. Both of these consider an 
$\eps$-regular cluster tournament~$G$ on~$k$ clusters with the property that 
for some cluster~$V_i$ the density of edges leaving~$V_i$ and the density of 
edges entering~$V_i$ are each roughly $1/2$. Lemma~\ref{ar_tourns_part_one} considers the case where for many clusters~$V_j$ the density of edges between $V_i$ and $V_j$ is large in both directions, showing that in this case $G$ contains a copy of a directed tree $T$ of the type considered. Lemma~\ref{ar_tourns_part_two} considers the alternative, namely that for almost all clusters~$V_j$ the density of edges between $V_i$ and $V_j$ is small in one direction, showing that in this case $G$ contains a copy of $T$ provided that $T_\Delta$ has large inweight and large outweight.

\begin{lemma} \label{ar_tourns_part_one}
Suppose that $1/n \ll 1/{\Delta'}, \beta \ll 1/k \ll \eps \ll \gamma \ll 
\alpha \ll 1/\Delta \ll 1$. Let $T$ be a directed tree on $n$ vertices with 
$|T_{\Delta'}| \leq \beta n$ and $|T_{\Delta}| \geq 2$, and let $G$ be an 
$\eps$-regular cluster tournament on clusters $V_1, \dots, V_k$, each of size 
$m \geq 2(1-\gamma)n/k$. Suppose also that for some~$i \in [k]$ we have  
$$\sum_{j \in [k] \sm \{i\}} d_{ij} \geq \frac{(1-3\gamma)k}{2} 
\textrm{\hspace{5mm}and } \sum_{j \in [k] \sm \{i\}} d_{ji} \geq 
\frac{(1-3\gamma)k}{2},$$ 
and also that there are at least $\alpha k$ values of $j \in [k] \sm \{i\}$ 
such that $d_{ij} \geq \alpha$ and $d_{ji} \geq \alpha$. 
Then $G$ contains a copy of $T$.
\end{lemma}

\proof 
Fix such a value of $i$, and introduce a new constant $\eps'$ with $\eps \ll 
\eps' \ll \gamma$. Since $\Delta \leq {\Delta'}$, we must have $T_{\Delta} 
\subseteq T_{\Delta'}$. Also, since $|T_{\Delta}| \geq 2$, we may choose an 
edge $t^- \rightarrow t^+$ of~$T_{\Delta}$, which therefore is also an edge 
of $T_{\Delta'}$. Let $T^+$ and $T^-$ be the two components formed when this 
edge is deleted from $T$, labelled so that $t^+ \in T^+$ and $t^- \in T^-$. 
Similarly, let~$T_{{\Delta'}}^+$ and $T_{{\Delta'}}^-$ be the two components 
formed by the deletion of the edge $t^- \rightarrow t^+$ from 
$T_{{\Delta'}}$, labelled with $t^+ \in T^+_{\Delta'}$ and $t^- \in 
T^-_{\Delta'}$. Then $T^+$ and $T^-$ partition the vertices of $T$, and there 
is precisely one edge of $T$ between $T^+$ and $T^-$, which is directed 
towards $T^+$. Furthermore, since $t^- \rightarrow t^+$ was an edge of 
$T_{\Delta}$, by Proposition~\ref{coretreeprops}(ii) we have $|T^+|, |T^-| 
\geq n/\Delta$.  

Let $J \subseteq [k] \sm \{i\}$ satisfy $|J| \geq \alpha k$ and also that for 
any $j \in J$ we have $d_{ij} \geq \alpha$ and $d_{ji} \geq \alpha$. 
Then $\sum_{j \in J} d_{ij} \geq \alpha^2 k$ and $\sum_{j \in J} d_{ji} \geq 
\alpha^2 k$. By Lemma~\ref{regularnbrhoods} (applied with $V_j' = \emptyset$ 
for each $j \notin J$) at most $\eps' m$ vertices of~$V_i$ have fewer than  
\begin{equation} \label{eq:densitiesinJ}
\sum_{j \in J} d_{ij}m - \eps' km \geq \alpha^2 km - \eps' km \geq 
\frac{\alpha^2km}{2} 
\end{equation}
outneighbours in $\bigcup_{j \in J} V_j$ or fewer than $\sum_{j \in J} d_{ji} 
m - \eps' km \geq \alpha^2km/2$ inneighbours in $\bigcup_{j \in J} V_j$. Also 
by Lemma~\ref{regularnbrhoods} at most $\eps' m$ vertices of $V_i$ have fewer 
than  
\begin{equation} \label{eq:densitiesoverall} 
\sum_{j \in [k] \sm \{i\}} d_{ij}m - \eps' km \geq \frac{(1-3\gamma-2\eps') 
km}{2} \geq (1-5\gamma) n 
\end{equation}
outneighbours in $\bigcup_{j \in [k] \sm \{i\}} V_j$ or fewer than $\sum_{j 
\in [k] \sm \{i\}} d_{ji} m - \eps' km \geq (1-5\gamma)n$ inneighbours in 
$\bigcup_{j \in [k] \sm \{i\}} V_j$. Finally, at most $m/2+1$ vertices of 
$V_i$ have fewer than $m/4$ inneighbours in $V_i$. So we may choose a set 
$S^+$ of $m/10$ vertices of $V_i$ which do not fall into any of these 
categories. Since $|T_{\Delta'}^+| \leq |T_{\Delta'}| \leq \beta n \leq 
m/30$, by Theorem~\ref{bestsofar} we may embed $T_{\Delta'}^+$ in $S^+$. Let 
$S_{\Delta'}^+$ be the set of vertices of $S^+$ occupied by this embedding of 
$T_{\Delta'}^+$, and let $v^+$ be the vertex to which $t^+$ was embedded.  
Recall that $|T^-| \geq n/\Delta$, so 
$$|T^+| = n-|T^-| \leq (1-1/\Delta)n.$$
Furthermore, every component of $T^+ - T^+_{\Delta'}$ is a component of $T - 
T_{\Delta'}$ and thus has order at most $n/{\Delta'}$ by 
Proposition~\ref{coretreeprops}. So by (\ref{eq:densitiesinJ}) and 
(\ref{eq:densitiesoverall}), and since $\gamma \ll 1/\Delta$, we may apply 
Lemma~\ref{onebyone}(b) to extend the embedding of $T_{\Delta'}^+$ in 
$S_{\Delta'}^+$ to an embedding of $T^+$ in $S_{\Delta'}^+ \cup \bigcup_{j 
\in [k] \sm \{i\}} V_j$ so that at least $\alpha^2n/3$ vertices of 
$\bigcup_{j \in J} V_j$ are occupied by this embedding of $T^+$. 

Now, at least $m/4-m/10 = 3m/20$ vertices of $V_i \sm S^+_{\Delta'}$ are 
inneighbours of $v^+$. For each $j \in [k] \sm \{i\}$, let $o_j$ denote the 
number of vertices of $V_j$ which are occupied by our embedding of $T^+$, and 
let $V_j' \subseteq V_j$ consist of those vertices of $V_j$ which are not 
occupied by this embedding. So $|V_j'| = m-o_j$ for each $j$. Note that since 
$d_{ij} + d_{ji} \leq 1$ we have $d_{ij} \leq 1-\alpha$ for each $j \in J$. 
Then by Lemma~\ref{regularnbrhoods}, at most $\eps' m$ vertices of $V_i$ have 
fewer than  
\begin{align} \label{eq:num_of_outnbrs}
\nonumber\sum_{j \in [k] \sm \{i\}} d_{ij} (m-o_j) - \eps' km &\geq
\sum_{j \in [k] \sm \{i\}} d_{ij} m - \eps' km - \sum_{j \in J} d_{ij} o_j - 
\sum_{j \in [k] \sm (\{i\} \cup J)} d_{ij} o_j  
\\\nonumber &\stackrel{(\ref{eq:densitiesoverall})}{\geq}  (1-5 \gamma)n - 
(1-\alpha)\sum_{j \in J} o_j - \sum_{j \in [k] \sm (\{i\} \cup J)} o_j 
\\\nonumber &\geq (1-5 \gamma)n - \sum_{j \in [k] \sm \{i\}} o_j + \alpha 
\sum_{j \in J} o_j  
\\ &\geq (1-5 \gamma) n - \sum_{j \in [k] \sm \{i\}} o_j + \alpha^3n/3 
\geq n- \sum_{j \in [k] \sm \{i\}} o_j + 2n/{\Delta'}
\end{align}
outneighbours in $\bigcup_{j \in [k] \sm \{i\}} V_j'$ or fewer than 
$$\sum_{j \in [k] \sm \{i\}} d_{ji} (m-o_j) - \eps' km \geq n- \sum_{j \in 
[k] \sm \{i\}} o_j + 2n/{\Delta'},$$  
inneighbours in $\bigcup_{j \in [k] \sm \{i\}} V_j'$. So we may choose a set 
$S^-$ of $m/10$ vertices of $V_i \sm S_{\Delta'}^+$, none of which fall into 
these two categories, and all of which are inneighbours of $v^+$. Since 
$|T_{\Delta'}^-| \leq |T_{\Delta'}| \leq \beta n \leq m/30$, by 
Theorem~\ref{bestsofar} we may embed $T_{\Delta'}^-$ in $S^-$. Let 
$S_{\Delta'}^-$ be the set of vertices of $S^-$ occupied by this embedding of 
$T_{\Delta'}^-$. Then since  
$$|T^-| = n-|T^+| \leq n - \sum_{j \in [k] \sm \{i\}} o_j,$$
the right hand side of (\ref{eq:num_of_outnbrs}) is at least $|T^-| + 
2n/{\Delta'}$. Also every component of $T^- - T^-_{\Delta'}$ is a component 
of $T - T_{\Delta'}$ (and so has order at most $n/{\Delta'}$ by 
Proposition~\ref{coretreeprops}(iv)). So by Lemma~\ref{onebyone} we may 
extend the embedding of $T_{\Delta'}^-$ in $S_{\Delta'}^-$ to an embedding of 
$T^-$ in $S_{\Delta'}^- \cup \bigcup_{j \in [k] \sm \{i\}} V_j'$. Then the 
embeddings of $T^+$ and $T^-$ do not overlap, and so together these 
embeddings form an embedding of $T$ in $G$. 
\endproof

Given an $\eps$-regular cluster tournament $G$ on clusters $V_1, \dots, V_k$, 
we define the \emph{reduced digraph of $G$ with parameter $d$}, denoted 
$R_{G}(d)$, to be the directed graph on vertex set $[k]$ in which $i 
\rightarrow j$ if and only if $d_{ij} \geq d$. Observe that since 
$d_{ij}+d_{ji} \leq 1$ for any $i$ and $j$, if $d > 1/2$ then $R_{G}(d)$ is 
an oriented graph.

\begin{lemma} \label{ar_tourns_part_two}
Suppose that $1/n \ll 1/{\Delta'}, \beta \ll 1/k \ll \eps \ll \gamma \ll 
\alpha \ll 1$. Let $T$ be a directed tree on $n$ vertices with $|T_{\Delta'}| 
\leq \beta n$, and let $y$ and $z$ be the outweight and inweight 
of~$T_{\Delta'}$ respectively. Let $G$ be an $\eps$-regular cluster 
tournament on clusters $V_1, \dots, V_k$, each of size $m \geq 
2(1-\gamma)n/k$.  Suppose that for some $i \in [k]$ we have  
$$\sum_{j \in [k] \sm \{i\}} d_{ij} \geq \frac{(1-3\gamma)k}{2} 
\textrm{\hspace{5mm}and } \sum_{j \in [k] \sm \{i\}} d_{ji} \geq 
\frac{(1-3\gamma)k}{2},$$ 
and also that there are at most $\alpha k$ values of $j \in [k] \sm \{i\}$ 
such that $d_{ij} \geq \alpha$ and $d_{ji} \geq \alpha$. Then: 
\begin{itemize}
\item[(i)] There are at most $2\alpha k$ values of $j \in [k] \sm \{i\}$ such 
that $d_{ij} < 1 - 2\alpha$ and $d_{ji} < 1 - 2\alpha$. 
\item[(ii)] Let $R := R_G(1-2\alpha)$. Then $|N^+_R(i)|, |N^-_R(i)| \geq 
(1-10\alpha)k/2$. 
\item[(iii)] If $y, z \geq 14 \alpha n$, then $G$ contains a copy of $T$.
\end{itemize}
\end{lemma}

\proof
Fix such an $i$, and introduce a new constant $\eps'$ with $\eps \ll \eps' 
\ll \gamma$. For (i), note that since $d_{ij} + d_{ji} \leq 1$ for any $j \in 
[k] \sm \{i\}$, and  
$$\sum_{j \in [k] \sm \{i\}} (d_{ij} + d_{ji}) \geq (1-3\gamma)k, $$
there are at most $3 \sqrt{\gamma}k \leq \alpha k$ values of $j \in [k] \sm 
\{i\}$ for which $d_{ij} + d_{ji} < 1-\sqrt{\gamma}$. So there are at most 
$2\alpha k$ values of $j \in [k] \sm \{i\}$ for which $d_{ij} < 1- \alpha - 
\sqrt{\gamma}$ and $d_{ji} < 1 - \alpha - \sqrt{\gamma}$, so (i) holds.  

For (ii), observe that by (i) we have
\begin{align*}
\frac{(1-3\gamma)k}{2} 
&\leq \sum_{j \in [k] \sm \{i\}} d_{ij}
\leq \sum_{\substack{j \in [k] \sm \{i\}\\d_{ij} \geq 1-2\alpha}} d_{ij} + 
\sum_{\substack{j \in [k] \sm \{i\}\\d_{ij}, d_{ji} <1-2\alpha}} d_{ij} + 
\sum_{\substack{j \in [k] \sm \{i\}\\d_{ij} \leq 2\alpha}} d_{ij}
\\&\leq |N^+_R(i)| +2\alpha k + 2\alpha k,
\end{align*}
so $|N^+_R(i)| \geq (1-10\alpha)k/2$. A similar calculation shows that 
$|N^-_R(i)| \geq (1-10\alpha)k/2$.

For (iii), let $N^+$ and $N^-$ denote $N_R^+(i)$ and $N_R^-(i)$ respectively, 
and let $ V^+ := \bigcup_{j \in N^+} V_j$ and $V^- := \bigcup_{j \in N^-} 
V_j$, so $V^+$ and $V^-$ are disjoint. By Lemma~\ref{regularnbrhoods}, $V_i$ 
contains at most $\eps' m$ vertices with fewer than  
\begin{align*} 
\sum_{j \in N^+} d_{ij} m - \eps' km &\geq |N^+_R(i)| (1-2\alpha)m - \eps' km 
\geq (1-10\alpha)(1-2\alpha)km/2 -\eps' km  
\\\nonumber&\geq (1-12\alpha - 2\eps')km/2 \geq (1-13\alpha) n
\end{align*}
outneighbours in $V^+$ and at most $\eps' m$ vertices with fewer than 
$\sum_{j \in N^-} d_{ji} m - \eps' km \geq (1-13\alpha) n$ inneighbours in 
$V^-$. Choose a set $S$ of $m/2$ vertices of $V_i$, not including any of 
these at most $2 \eps' m$ vertices. Since $|T_{\Delta'}| \leq \beta n \leq 
m/6$, by Theorem~\ref{bestsofar} we may embed $T_{\Delta'}$ in $S$. Let 
$S_{\Delta'}$ be the set of vertices of $S$ occupied by this embedding of 
$T_{\Delta'}$. Also let $T_1$ be the tree formed by $T_{\Delta'}$ and all of 
its outcomponents, and let $T_2$ be the tree formed by $T_{\Delta'}$ and all 
of its incomponents. Note that all of these out- and incomponents have order 
at most $n/{\Delta'} \ll \alpha n$ by Proposition~\ref{coretreeprops}(iv). In 
addition $|T_1| = n-z \leq (1-14\alpha) n$ and $|T_2| = n-y \leq 
(1-14\alpha)n$. So by Lemma~\ref{onebyone} we may extend the embedding of 
$T_{\Delta'}$ in $S_{\Delta'}$ to an embedding of $T_1$ in $S_{\Delta'} \cup 
V^+$. Similarly by Lemma~\ref{onebyone} we may extend the embedding of 
$T_{\Delta'}$ in $S_{\Delta'}$ to an embedding of $T_2$ in $S_{\Delta'} \cup 
V^-$. Then these embeddings do not overlap outside $T_{\Delta'}$, so we may 
combine them to form an embedding of $T$ in $G$. 
\endproof

To finish this section we shall show how Lemma~\ref{regularity_lemma} can be 
used to show that Sumner's universal tournament conjecture holds for any 
large and almost-regular tournament with a small core tree. Actually we shall 
prove a slightly stronger result in this case, considering a tournament on 
fewer than $2n-2$ vertices. Later on we shall make use of the fact that we 
have a little room to spare in the order of the tournament. Much of the work 
for this lemma is done by the two previous lemmas.

\begin{lemma} \label{ar_tourns_small_core_tree}
Suppose that $1/n \ll 1/{\Delta'}, \beta \ll \gamma \ll 1/\Delta \ll 1$. Let 
$T$ be a directed tree on $n$ vertices such that $|T_{\Delta'}| \leq \beta n$ and 
$|T_{\Delta}| \geq 2$. Let $G$ be a $\gamma$-almost-regular tournament on at 
least $(2-\gamma)n$ vertices. Then $G$ contains a copy of $T$. 
\end{lemma}

\proof
Introduce new constants $\eps, \eps', \alpha, M,$ and $M'$ with 
$$1/n \ll 1/{\Delta'}, \beta \ll 1/M \ll 1/M' \ll \eps \ll \eps' \ll \gamma 
\ll \alpha \ll 1/\Delta \ll 1.$$ 
If $|G| \geq (2+\gamma) n$, then $G$ contains a copy of $T$ by 
Theorem~\ref{approxversion}(i). So we may assume that $|G| = (2 \pm 
\gamma)n$. Observe that $d^+(v), d^-(v) \geq (1-\gamma)(|G|-1)/2 \geq 
(1-2\gamma)n$ for all $v \in G$. 

Since $\Delta \leq {\Delta'}$, we must have $T_{\Delta} \subseteq 
T_{\Delta'}$. Also, since $|T_{\Delta}| \geq 2$, we may choose an edge $t^- 
\rightarrow t^+$ of $T_{\Delta}$, which must also lie in $T_{\Delta'}$. Let 
$T^+$ and $T^-$ be the two components formed when this edge is deleted from 
$T$, labelled so that $t^+ \in T^+$ and $t^- \in T^-$. Similarly, let 
$T_{{\Delta'}}^+$ and $T_{{\Delta'}}^-$ be the two components formed by the 
deletion of the edge $t^- \rightarrow t^+$ from $T_{{\Delta'}}$, labelled 
with $t^+ \in T^+_{\Delta'}$ and $t^- \in T^-_{\Delta'}$. Then $T^+$ and 
$T^-$ partition the vertices of $T$, and there is precisely one edge of $T$ 
between $T^+$ and $T^-$, which is directed towards $T^+$. Furthermore, 
$|T^+|, |T^-| \geq n/\Delta$.  

Let disjoint subsets $V_1, \dots, V_k$ and a subgraph $G^* \subseteq G$ 
satisfy the conditions of Corollary~\ref{clustertourns}. So $M' \leq k \leq 
M$, and $G^*$ is an $\eps$-regular cluster tournament on clusters $V_1, 
\dots, V_k$ of equal size $m$, where 
\begin{equation} \label{eq:boundsonm}
\frac{2(1-\gamma)n}{k} \leq \frac{(2-\gamma)n-3\eps n}{k} \leq m \leq 
\frac{(2+\gamma)n}{k}. 
\end{equation}

Also, for each $v \in G^*$ we have $d^+_{G^*}(v) \geq d^+_G(v) - 2\eps |G| 
\geq  d^+_G(v) - 5\eps n$ and $d^-_{G^*}(v) \geq  d^-_G(v) - 5\eps n$. So for 
each $i \in [k]$ we have  
\begin{align} \label{eq:sumofdij}
\sum_{j \in [k] \sm \{i\}} d_{ij} 
&= \sum_{j \in [k] \sm \{i\}} \frac{e_{G^*}(V_i \rightarrow V_j)}{m^2} 
\geq \sum_{v \in V_i} \frac{d^+_{G^*}(v)-m}{m^2} 
\\ \nonumber&\geq \sum_{v \in V_i} \frac{d^+_{G}(v) - 5\eps n-m}{m^2} 
\geq \frac{(1 - 2\gamma)n - 5\eps n-m}{m} 
\stackrel{(\ref{eq:boundsonm})}{\geq} \frac{(1 - 3\gamma) k}{2}, 
\end{align}
and similarly $\sum_{j \in [k] \sm \{i\}} d_{ji} \geq (1 - 3\gamma) k/2$.

So if there exists some $i \in [k]$ for which there are at least $\alpha k$ 
values of $j \in [k] \sm \{i\}$ such that $d_{ij} \geq \alpha$ and $d_{ji} 
\geq \alpha$, then by Lemma~\ref{ar_tourns_part_one} we may embed $T$ in 
$G^*$, and therefore in $G$. So we may assume that for each $i \in [k]$ fewer 
than $\alpha k$ values of $j \in [k] \sm \{i\}$ satisfy $d_{ij} \geq \alpha$ 
and $d_{ji} \geq \alpha$. Then by Lemma~\ref{ar_tourns_part_two} we may 
assume that $R := R_G(1-2\alpha)$ has  
\begin{equation} \label{eq:semidegofR}
\delta^0(R) \geq (1-10\alpha)k/2.
\end{equation}
Let $y$ be the number of vertices in outcomponents of $T_{\Delta'}$, and let 
$z$ be the number of vertices in incomponents of $T_{\Delta'}$, so $y + z + 
|T_{\Delta'}| = n$. So if $y, z \geq 14 \alpha n$ then $G^*$ (and therefore 
$G$) contains a copy of $T$ by Lemma~\ref{ar_tourns_part_two}. We may 
therefore assume without loss of generality that $z < 14 \alpha n$.%
\COMMENT{Since if $y < 14 \alpha n$, reversing the direction of all edges 
in $T$ and in $G$ gives $z <14 \alpha n$ and the same embedding problem} 

Now, since $|R| = k$ we may choose a vertex $i \in R$ with $d^+_R(i) \leq 
k/2$. Then we may choose a vertex $j \in N^+_R(i)$ with at most $d^+_R(i)/2$ 
outneighbours in $N^+_R(i)$. So $i \rightarrow j$ and $|N^+_R(i) \cap 
N^+_R(j)| \leq k/4$. For this choice of $i$ and $j$, let  
\begin{align*}
A &:= N^+_R(i) \cap N^+_R(j), 
\\ B &:= N^+_R(i) \sm N^+_R(j), 
\\ C &:=  N^+_R(j) \sm N^+_R(i).
\end{align*} 
Then $A, B$ and $C$ are disjoint, and $|B|, |C| \geq k/2 - 5\alpha k -|A| 
\geq k/4 - 5\alpha k$ by (\ref{eq:semidegofR}). 
Now, choose a set $S^+$ of $m/2$ vertices of $V_j$ such that each vertex $v 
\in S^+$ has  
\begin{enumerate}
\item[(i)] at least $m/2$ inneighbours in $V_i$,
\item[(ii)] at least $\sum_{\ell \in A} d_{j\ell} m - \eps' km \geq 
(1-2\alpha)m|A| - \eps' km$ outneighbours in $\bigcup_{\ell \in A} V_\ell$, and 
\item[(iii)] at least $\sum_{\ell \in C} d_{j\ell} m - \eps' km \geq 
(1-2\alpha)m|C| - \eps' km$ outneighbours in $\bigcup_{\ell \in C} V_\ell$. 
\end{enumerate}
We can be sure that such a choice is possible, as by 
Lemma~\ref{regularnbrhoods} there are at most $2\eps' m$ vertices of $V_j$ 
which fail either of (ii) and (iii), and since $G^*[V_i \rightarrow V_j]$ is 
$\eps$-regular with density $d_{ij} \geq 1-2\alpha$ there are at most $\eps 
m$ vertices of $V_j$ which fail (i). Then since $|T_{\Delta'}^+| \leq \beta n 
\leq m/6$, by Theorem~\ref{bestsofar} we can embed $T_{\Delta'}^+$ in $S^+$. 
Let $v^+$ be the vertex to which $t^+$ is embedded. Then $v^+$ has at least 
$m/2$ inneighbours in $V_i$. Choose a set $S^-$ of $m/3$ of these 
inneighbours so that every vertex $v \in S^-$ has at least  
\begin{equation} \label{eq:outnbrsofS-}
\sum_{\ell \in A \cup B = N^+_R(i)} d_{i\ell} m - \eps' km \geq 
(1-2\alpha)m|N_R^+(i)| - \eps' km \stackrel{(\ref{eq:semidegofR})}{\geq} 
(1-13\alpha)n 
\end{equation} outneighbours in $\bigcup_{\ell \in A \cup B} V_\ell.$ Again 
we can be sure that such a choice is possible, since by 
Lemma~\ref{regularnbrhoods} at most $\eps' m$ vertices of $V_i$ fail this 
condition. Then since $|T_{\Delta'}^-| \leq \beta n \leq m/9$, by 
Theorem~\ref{bestsofar} we can embed $T_{\Delta'}^-$ in $S^-$. Let 
$S_{\Delta'}^+$ and $S_{\Delta'}^-$ be the sets of vertices of $G$ occupied 
by $T_{\Delta'}^+$ and $T_{\Delta'}^-$ respectively. 

Let $T_3$ be the tree formed by $T_{\Delta'}$ and all of its incomponents. 
Let $T_4$ be the tree formed by $T_{\Delta'}^+$ and all of its outcomponents, 
and let $T_5$ be the tree formed by $T_{\Delta'}^-$ and all of its 
outcomponents in $T^-$ (i.e. all of its outcomponents except $T^+$). Note 
that $T_3 \cup T_4 \cup T_5 = T$. Then $|T_3| = |T_{\Delta'}| + z < 15\alpha 
n$, $|T_4| \leq |T^+| \leq n - |T^-| \leq (1-1/\Delta)n$, and similarly 
$|T_5| \leq (1-1/\Delta)n$. Every vertex of $G$ has at least $(1-2\gamma)n$ 
inneighbours in $G$, so by Lemma~\ref{onebyone}(c) we may extend the 
embedding of $T_{\Delta'}$ in $S_{\Delta'}^+ \cup S_{\Delta'}^-$ to an 
embedding of $T_3$ in $G$. For each $\ell \in [k]\sm \{i\}$, let $V_\ell' 
\subseteq V_\ell$ consist of the vertices of $V_\ell$ which are not occupied 
by this embedding. 

By (ii) and (iii), every vertex of $S_{\Delta'}^+$ then has at least 
$(1-2\alpha)(|A|+|C|)m - 2\eps' km - |T_3| \geq (1-28\alpha)n$ outneighbours 
in $\bigcup_{\ell \in A \cup C} V_\ell'$ (here we also use the fact that 
$|A|+|C| = |N_R^+(j)| \geq (1-10\alpha)k/2$ by (\ref{eq:semidegofR})). Since 
also $1/{\Delta'} \ll \alpha \ll 1/\Delta$ and every component of $T_4 - 
T_{\Delta'}^+$ has order at most $n/{\Delta'}$, by Lemma~\ref{onebyone} we 
may extend the embedding of $T_{\Delta'}^+$ in $S_{\Delta'}^+$ to an 
embedding of $T_4$ in $S_{\Delta'}^+ \cup \bigcup_{\ell \in A \cup C} 
V_\ell'$. Furthermore, since every vertex of $S_{\Delta'}^+$ has at least 
$(1-2\alpha)|C|m - \eps' km - |T_3| \geq n/\Delta$ outneighbours in 
$\bigcup_{\ell \in C} V_\ell'$, and $|T_4 - T_{\Delta'}^+| = |T^+ - T_3| \geq 
n/2\Delta$, by Lemma~\ref{onebyone}(b) we can ensure that this embedding of 
$T_4$ occupies at least $n/2\Delta$ vertices of $\bigcup_{\ell \in C} 
V_\ell'$. So crucially at most $|T_4| - n/2\Delta$ vertices of $T_4$ are 
embedded in $\bigcup_{\ell \in A \cup B} V_\ell$. For each $\ell \in A \cup 
B$, let $V''_\ell \subseteq V_\ell$ consist of those vertices which are not 
occupied by the embedding of $T_3$ and $T_4$.  

Finally, by (\ref{eq:outnbrsofS-}), every vertex of $S_{\Delta'}^-$ has at 
least  
$$(1-13\alpha)n - (|T_4| - n/2\Delta) - |T_3| \geq n - |T_4| + n/3\Delta$$ 
outneighbours in $\bigcup_{\ell \in A \cup B} V''_\ell$. Since $|T_5 - 
T_{\Delta'}^-| \leq n-|T_4|$, by Lemma~\ref{onebyone}(c) we can extend the 
embedding of $T_{\Delta'}^-$ in $S_{\Delta'}^-$ to an embedding of $T_5$ in 
$S_{\Delta'}^- \cup \bigcup_{\ell \in A \cup B} V''_\ell$. Then the 
embeddings of $T_3$, $T_4$ and $T_5$ do not overlap outside $S_{\Delta'}^+ 
\cup S_{\Delta'}^-$, and so together form an embedding of $T$ in~$G$. 
\endproof

\section{Embedding trees in robust outexpander tournaments} \label{sec:robustexpander}

Let $G$ be a 
tournament on $n$ vertices, and let $\mu \leq \nu$ be positive constants. 
Then the \emph{robust outneighbourhood} $RN_\mu^+(S)$ of a set $S \subseteq 
V(G)$ is the set of vertices of $G$ with at least $\mu n$ inneighbours in 
$S$. We say that $G$ is a \emph{robust $(\mu, \nu)$-outexpander} if for any 
$S \subseteq V(G)$ with $\nu n \leq |S| \leq (1-\nu) n$ we have 
$|RN_\mu^+(S)| \geq |S| + \mu n$. 

If a tournament $G$ is not a robust outexpander, then the following lemma 
shows that $G$ contains two subtournaments which partition the vertices of 
$G$ and which have almost all edges between them directed the same way. 

\begin{lemma}[\cite{KMO}, Lemma 2.8] \label{not_robust_expander_split}
Suppose that $1/n \ll \mu \ll \nu$, that $G$ is a tournament on $n$ vertices
and that $G$ is not a robust $(\mu, \nu)$-outexpander. Then we can partition
$V(G)$ into sets $S$ and $S'$ such that $\nu n < |S|, |S'| < (1-\nu)n$
and $e(G[S \rightarrow S']) \leq 4 \mu n^2$.
\end{lemma}

By iterating this split,  we obtain a decomposition of $G$ into sets $S_i$ which 
either induce robust expanders or are small, and where for all $i < j$, almost all edges are
directed from $S_i$ to $S_j$. (So if all the $S_i$ are small, then $G$ is close to 
being a transitive tournament.) 
We will use this decomposition in Section~\ref{sec:mainproof} to prove Theorem~\ref{main}.

\begin{lemma} \label{tournament_split} Suppose that $1/n \ll \mu \ll \nu \ll 
\eta \ll \gamma \ll 1$. Let $G$ be a tournament on $n$ vertices. Then we may 
choose disjoint subsets $S_1, \dots, S_r$ of $V(G)$ such that: 
\begin{itemize}
\item[(i)] $|\bigcup_{i \in [r]} S_i| \geq (1-\gamma)n$,
\item[(ii)] for each $i \in [r]$, any vertex $v \in S_i$ has at most $\gamma 
n$ inneighbours in $\bigcup_{j > i} S_j$ and at most $\gamma n$ outneighbours 
in $\bigcup_{j < i} S_j $, and 
\item[(iii)] for each $i \in [r]$, either $G[S_i]$ is a robust $(\mu, 
\nu)$-outexpander with $\delta^0(G[S_i]) \geq \eta n$ or $|S_i| < \gamma n$. 
\end{itemize}
\end{lemma}

\proof
We shall use a modified version of an algorithm from~\cite{KMO}, which keeps 
track of an ordered family 
$\mathcal{S}^\tau$ of disjoint subsets of $V(G)$, and a set~$B^\tau$ of bad
edges of $G$, at each time $\tau$. The analysis of this algorithm is also 
similar to the analysis in~\cite{KMO}. Initially, let $\mathcal{S}^1 := 
(V(G))$, and 
let $B^1 := \emptyset$. Then at time $\tau \geq 1$, we have $\mathcal{S}^\tau 
= (S_1^\tau, 
\dots, S_\tau^\tau)$, and the algorithm proceeds as follows.
\begin{enumerate}
\item Let $S_\ell^\tau$ be the largest member of $\mathcal{S^\tau}$ which is 
not a robust $(\mu, \nu)$-outexpander with $\delta^0(G[S_\ell^\tau]) \geq 
\eta n$. If there is no such member of $\mathcal{S^\tau}$, or if 
$|S_\ell^\tau| < \gamma n$, then terminate. If there is more than one largest 
such member, then choose one of these arbitrarily. 
\item If some $v \in S_\ell^\tau$ has $d^+_{G[S_\ell^\tau]}(v) < \eta n$, 
then let 
$$\mathcal{S}^{\tau+1} := (S_1^\tau, \dots, S_{\ell-1}^\tau, S_\ell^\tau \sm
\{v\},
\{v\}, S_{\ell+1}^\tau, \dots, S_{\tau}^\tau), $$
let $B^{\tau+1} := B^\tau \cup E(\{v\} \rightarrow S_\ell^\tau \sm \{v\})$, and 
proceed to step (5).
\item Similarly, if some $v \in S_\ell^\tau$ has $d^-_{G[S_\ell^\tau]}(v) < \eta 
n$, then let
$$\mathcal{S}^{\tau+1} := (S_1^\tau, \dots, S_{\ell-1}^\tau, \{v\}, 
S_\ell^\tau \sm 
\{v\}, S_{\ell+1}^\tau, \dots, S_{\tau}^\tau),$$
let $B^{\tau+1} := B^\tau \cup E(S_\ell^\tau \sm \{v\} \rightarrow \{v\})$, and 
proceed to step (5).
\item If $G[S_\ell^\tau]$ is not a robust $(\mu, \nu)$-outexpander then apply 
Lemma~\ref{not_robust_expander_split} to partition the vertices of 
$S_\ell^\tau$ 
into sets $S'$ and $S''$ such
that $\nu|S_\ell^\tau| \leq |S'|, |S''| \leq (1-\nu)|S_\ell^\tau|$ and at most 
$4\mu |S_\ell^\tau|^2$ edges of $G[S_\ell^\tau]$ are directed from $S''$ to 
$S'$. 
Then let
$$\mathcal{S}^{\tau+1} := (S_1^\tau, \dots, S_{\ell-1}^\tau, S', S'',
S_{\ell+1}^\tau,
\dots, S_{\tau}^\tau)$$
and let $B^{\tau+1} := B^\tau \cup E(S'' \rightarrow S')$.
\item Finally, for each $i \in [\tau+1]$, delete from $S_i^{\tau+1}$ any 
vertex 
$v$ which lies in more than $\sqrt{\eta}n$ edges of $B^{\tau+1}$.
\end{enumerate}

At any time $\tau$, if the algorithm does not terminate at step (1) then 
$S_\ell^\tau$ will be split in precisely one of steps (2), (3) and (4). So at 
each time $\tau$, either the algorithm terminates or $|\mathcal{S}^\tau|$ 
increases from $\tau$ to $\tau+1$ (in forming $\mathcal{S}^{\tau+1}$) by 
reducing the size of the largest piece. Therefore the algorithm must 
terminate at some time $\tau_{end} \leq n$. Take $r := \tau_{end}$, and $S_i 
:= S_i^r$ for each $i$. Then since the algorithm terminated at step (1) of 
time $r$, (iii) must hold. 

To see (i), observe that the split in step (4) will occur for at most 
$1/\gamma \nu$ times $\tau < \tau_{end}$. This is because any set obtained by 
a split in step (4) must have size at least $\gamma \nu n$ (since 
$|S^\tau_\ell| \geq \gamma n$, and the sets $S', S''$ obtained have $|S'|, 
|S''| \geq \nu |S^\tau_\ell|$). Also, at each time $\tau \leq \tau_{end}$, 
the number of edges added to form $B^{\tau+1}$ from $B^\tau$ is at most $\eta 
n$ if the algorithm carried out the split in step (2) or (3), and at most 
$4\mu n^2$ if the algorithm carried out the split in step (4). Since 
$\tau_{end} \leq n$, and the split in step (4) is carried out in at most 
$1/\gamma \nu$ steps, we must have  
$$|B^{\tau_{end}}| \leq \eta n^2 + 4 \mu n^2/\nu \gamma \leq 2 \eta n^2.$$ 
Since $B^1 \subseteq \dots \subseteq B^{\tau_{end}}$, any vertex of $G$ which 
was ever deleted in step (5) must lie in at least $\sqrt{\eta} n$ edges of 
$B^{\tau_{end}}$, and so at most $4 \sqrt{\eta} n \leq \gamma n$ vertices of 
$G$ can have been deleted in step (5) over the entire course of the 
algorithm. But any vertex which was not deleted lies in some $S_i$, and so 
(i) holds. 

Finally, for (ii) fix any $i \in [r]$ and any $v \in S_i$. Observe that all 
edges directed from $v$ to $\bigcup_{j < i} S_j$ and all edges directed from 
$\bigcup_{j > i} S_j$ to $v$ are contained in $B^r$. This means that there 
are at most $\sqrt{\eta}n$ such edges, as otherwise $v$ would have been 
deleted in step (5) at some point. Since $i$ and $v$ were arbitrary, (ii) 
must hold. 
\endproof

We now consider the case when $G$ is a robust 
outexpander. Lemma~4.1 of~\cite{KMO} stated that if $T$ is a directed tree on $n$ 
vertices, and $G$ is a robust outexpander tournament on at least 
$(2+\alpha)n$ vertices with large minimum semidegree, then $G$ contains a 
copy of $T$. However, in the proof of this lemma, the $\alpha n$ error term 
was only needed in the case when $T_\Delta$ is small. Indeed, in this section 
we modify this proof to show that Sumner's universal tournament conjecture 
holds for such $G$ in the case when $T_\Delta$ is large. This is the 
following lemma. 

\begin{lemma} \label{ar_tourns_large_core_tree}
Suppose that $1/n \ll 1/\Delta \ll \mu \ll \nu \ll \eta \ll \gamma \ll \beta 
\ll 1$. Let $T$ be a directed tree on $n$ vertices such that $|T_\Delta| \geq \beta 
n$, and let $G$ be a robust $(\mu, \nu)$-outexpander tournament on at least 
$(2-\gamma)n$ vertices, with $\delta^0(G) \geq \eta |G|$. Then $G$ contains a 
copy of $T$. 
\end{lemma}

Before we can present the proof of this lemma, we must give some definitions 
from \cite{KMO}. Let $V_1, \dots, V_k$ be disjoint sets of equal size. A 
digraph $G$ on vertex set $V_1, \dots, V_k$ is a \emph{$\eps$-regular 
$d$-dense cycle of cluster tournaments} if for each $i$, $G[V_i]$ is a 
tournament and $G[V_i \rightarrow V_{i+1}]$ is $\eps$-regular with density at 
least $d$ (where addition on the index of $V_{i+1}$ is taken modulo $k$). The 
following lemma from~\cite{KMO} (an immediate consequence of two results 
from~\cite{KOT}) will help us to find such digraphs. 

\begin{lemma}[\cite{KMO}, Lemma 2.7] \label{KOT_combined}
Suppose that $1/n \ll 1/M \ll 1/M' \ll \eps \ll d \ll \mu \ll \nu \ll \eta 
\ll 1$. 
Let~$G$ be a tournament on $n$ vertices which is a robust $(\mu,
\nu)$-outexpander with $\delta^0(G) \geq \eta n$. Then~$G$ contains an
$\eps$-regular $d$-dense cycle of cluster tournaments on clusters $V_1, 
\dots, 
V_k$, where $|\bigcup_{i=1}^k V_i| > (1-\eps)n$, and $M' \leq k \leq M$.
\end{lemma}

Let $T$ be a directed tree. Then the \emph{distance} between vertices $u, v 
\in T$, denoted $d(u, v)$, is the length of the shortest path connecting~$u$ 
and $v$ in the underlying graph $T_{under}$. Similarly for a set $X$ of 
vertices of $T$, the distance $d(u, X)$ is the minimum of $d(u, x)$ taken 
over all vertices $x \in X$. If $T$ is a rooted tree with root $r$, then the 
\emph{children} of a vertex $u \in T$ are those neighbours $v$ of $u$ for 
which $d(r, u) = d(r, v)+1$. 

Let $T$ be a tree on $n$ vertices, rooted at $t_1$, and let $H \subseteq 
V(T)$. Also let $k$ be a positive integer. For any vertex $x \in T$, there is 
a unique path in $T$ from $x$ to $t_1$; let~$P_x$ denote the set of the first 
$k$ vertices of this path, starting from $x$. Let $H^1 := \bigcup_{x \in H} P_x$, and then for each $i \geq 1$ let 
$H^{i+1}$ be formed from $H^{i}$ by adding the vertices of $P_x$ for any $x 
\in H^{i}$ with at least two children in $H^{i}$. After at most $n$ steps we 
must have $H^{i} = H^{i+1}$, when we terminate the process. We refer to this 
final $H^{i}$ as \emph{$H$ with leading paths included}, 
denoted~$\mathcal{P}_k(H)$. So $H \subseteq \mathcal{P}_k(H) \subseteq V(T)$. 
Note that $\mathcal{P}_k(H)$ depends on both the value of $k$ and the root 
$t_1$ of $T$. 

We may now present the key lemma from~\cite{KMO} we shall use to prove 
Lemma~\ref{ar_tourns_large_core_tree}. This says that a directed tree of bounded degree 
can be embedded in a robust outexpander tournament of large minimum 
semidegree such that the vertices in a small set $H$ of vertices of $T$ are 
embedded within a chosen set $U \subseteq V(G)$.\COMMENT{I have this as lemma 
4.6 - my version doesn't have a lemma 4.7. I hope nothing has gone wrong...} 

\begin{lemma}[\cite{KMO}, Lemma~4.6] \label{embed_with_restrictions}
Suppose that $1/n \ll 1/\Delta, 1/k \ll \eps \ll d \ll \alpha, \lambda \leq 
1/2$, that $m := n/k$, that $\lambda \leq \alpha/4$ and that $\delta := 
d\lambda/8k$. Let $T$ be a directed tree on $n$ vertices rooted at $t_1$ and 
with $\Delta(T) \leq \Delta$. Let $H \subseteq V(T)$ be such that $|H| \leq 
\delta n/7k$ and $|\{x \in T \colon 1 \leq d(x, \mathcal{P}_k(H)) \leq k^3\}| 
\leq \delta n$. Let $G$ be an $\eps$-regular $d$-dense cycle of cluster 
tournaments on clusters $V_1, \dots, V_k$, each of size $(1+\alpha)m$, and 
let $U \subseteq V_1 \cup \dots \cup V_k$ have size $|U| \geq \lambda n$. 
Then $T$ can be embedded in~$G$ so that each vertex $t \in H$ is embedded to 
some $u \in U$. 
\end{lemma}

We will also use the following lemma, again from~\cite{KMO}. This shows that 
we can extend $T_\Delta$ to an `extended tree' $T_{ext}$, with desired 
properties. We will apply Lemma~\ref{embed_with_restrictions} to $T_{ext}$ 
and embed $H$ within a set $U$ of vertices of high in- and outdegree. 

\begin{prop}[\cite{KMO}, Lemma 4.5] \label{extendedtree}
Suppose that $1/n, 1/\Delta^* \ll 1/\Delta, 1/k, \omega \ll 1$. Let $T$ be a 
directed tree on $n$ vertices. Choose any vertex $t_1 \in T_\Delta$ as the root of 
$T$. Then there exists a subtree $T_{ext}$ of $T$ and a subset $H \subseteq 
V(T_{ext})$ which satisfy the following properties. 
\begin{itemize}
\item[(i)] $T_\Delta \subseteq T_{ext}$.  
\item[(ii)]  $\Delta(T_{ext}) \leq \Delta^*$.
\item[(iii)]  For any edge $e$ between $T - T_{ext}$ and $T_{ext}$, the 
endvertex of $e$ in $T_{ext}$ lies in $H$. 
\item[(iv)]  The number of
vertices $v \in T_{ext}$ which satisfy $1 \leq d(v, \mathcal{P}_k(H)) \leq
k^3$ is at most $\omega n$.
\item[(v)]  $|H| \leq n/\Delta^{k^{1/\omega}}$.
\end{itemize}
\end{prop}

The final lemma we shall need to prove Lemma~\ref{ar_tourns_large_core_tree} 
gives standard Chernoff-type bounds for the binomial and hypergeometric 
distributions. The binomial random variable $X$ with parameters $(n, p)$ is 
defined to be the number of successes in $n$ independent trials, each of 
which has probability $p$ of success. So $\mathbb{E}X = np$. The 
hypergeometric random variable $Y$ with parameters $(n,m,k)$ is defined as 
follows. Let $N$ be a set of size $n$, and fix a set $S \subseteq N$ of size 
$|S|=m$. Now choose a set $T \subseteq N$ of size $|T|=k$ uniformly at 
random. Then $Y=|T \cap S|$. Note that $\mathbb{E}Y = km/n$. 

\begin{prop} [\cite{JLR}, Corollary 2.3 and Theorem 2.10]\label{chernoff}
Suppose $X$ has binomial or hypergeometric distribution and $0<a<3/2$. Then
$\mathbb{P}(|X - \mathbb{E}X| \ge a\mathbb{E}X) \le 2
e^{-\frac{a^2}{3}\mathbb{E}X}$.
\end{prop}

\medskip
\noindent {\bf Proof of Lemma~\ref{ar_tourns_large_core_tree}.}
We begin by introducing new constants $\Delta^*, M, M', \eps, d$ and $\alpha$ 
which satisfy 
$$ 1/n \ll 1/\Delta^* \ll 1/M \ll 1/M' , 1/\Delta \ll \eps \ll d \ll \mu \ll 
\nu \ll \eta \ll \gamma \ll \alpha \ll \beta \ll 1.$$ 
Now, if $|G| \geq (2+\gamma)n$, then by Theorem~\ref{approxversion}(i), $G$ 
contains a copy of $T$. So we may assume that $|G| = (2 \pm \gamma)n$. Since 
$G$ is a robust $(\mu, \nu)$-outexpander with $\delta^0(G) \geq \eta |G|$, 
Lemma~\ref{KOT_combined} implies that $G$ contains an $\eps$-regular 
$d$-dense cycle of cluster tournaments on clusters $V_1, \dots, V_k$ each of 
equal size between $(1-\eps)|G|/k \geq (1-\eps)(2-\gamma)m \geq 2(1-\gamma)m$ 
and $|G|/k \leq (2+\gamma)m$, where $m := n/k$ and $M' \leq k \leq M$. So we 
may remove vertices from each $V_i$ to obtain a $2\eps$-regular $(d/2)$-dense 
cycle of cluster tournaments $G'$ on clusters $V_1', \dots, V_k'$ each of 
size $2(1-\gamma) m$. So $|G'| = 2(1-\gamma)n$. Let 
$$\delta := d\alpha\beta/160k.
$$ 

Choose any vertex $t_1 \in T_\Delta$ as the root of $T$. Then let $T_{ext}$ 
and $H$ satisfy the properties of Proposition~\ref{extendedtree}, with 
$\omega := \delta \beta$. Let $T_1$ denote the subtree of $T$ formed by 
$T_{ext}$ and all of its outcomponents, and let $T_2$ denote the subtree of 
$T$ formed by $T_{ext}$ and all of its incomponents. Since $T_\Delta 
\subseteq T_{ext}$ (this is (i) of Proposition~\ref{extendedtree}), all of 
these incomponents and outcomponents have order at most $n/\Delta$ by 
Proposition~\ref{coretreeprops}. Let $x := |T_{ext}|, y := |T_1 - T_{ext}|, z 
:= |T_2 - T_{ext}|$, so $x+y+z=n$. Since $T_\Delta \subseteq T_{ext}$, we 
have $x \geq \beta n$. Also, all but at most $2y+x-\alpha n/2$ vertices of 
$G$ have at least $y+x/2-\alpha n/4$ outneighbours, and all but at most 
$2z+x-\alpha n/2$ vertices of $G$ have at least $z + x/2 - \alpha n/4$ 
inneighbours. So at least $(2-\gamma)n - 2y - 2z - 2x + \alpha n \geq \alpha 
n/2$ vertices of $G$ satisfy both of these conditions. Let $U_0$ be the set 
of these vertices, so $|U_0| \geq \alpha n/2$, and each $v \in U_0$ has at 
least $y+ x/2 - \alpha n/4$ outneighbours and at least $z + x/2 - \alpha n/4$ 
inneighbours. 

From each cluster $V_i'$ of $G'$ choose a set $X_i$ of $(1+\alpha)x/k$ 
vertices uniformly at random, and let $X := X_1 \cup \dots \cup X_k$. Then 
$|X| = (1+\alpha)x$. For any single vertex $u \in G'$, the probability 
that~$u$ is included in $X$ is $(1+\alpha)x/|G'| \geq x/2n$, so by 
Proposition~\ref{chernoff}, with probability at least $2/3$ the set $U := X \cap 
U_0$ satisfies $|U| \geq \alpha x/5 \geq \alpha \beta n/5$. Also, for any 
vertex $v \in U$, the expected number of outneighbours of~$v$ outside $X$ is 
at least 
\begin{align*}
\left(y + \frac{x}{2} - \frac{\alpha n}{4}\right) \left(1 - 
\frac{(1+\alpha)x}{|G'|}\right) 
\geq &  y - \frac{\alpha n}{4} + \frac{x}{2} - 
\frac{(1+\alpha)xy}{2(1-\gamma)n} - \frac{(1+\alpha)x^2}{4(1-\gamma)n} 
\\ \geq & y - \frac{\alpha n}{4} + \frac{2xn - 2xy - x^2 - 2\gamma xn - 2 
\alpha xy - \alpha x^2}{4(1-\gamma)n} 
\\ \geq & y + \frac{x^2}{4n} - 2\alpha n
 \geq  y + \frac{\beta^2 n}{4} - 2\alpha n
 \geq  y + 2\alpha n,
\end{align*}
where in the first inequality of the third line we used the fact that 
$2n-2y-x \geq x$. 
A similar calculation shows that for each $v \in U$, the expected number of
inneighbours of $v$ outside $X$ is at least $z + 2 \alpha n$. So by
Proposition~\ref{chernoff} we find that with probability at least $2/3$, every 
vertex 
$v \in U$ has at least $y + \alpha n$ outneighbours outside $X$ and at least
$z + \alpha n$ inneighbours outside $X$. Fix a choice of $X$ such that both
these events of probability at least $2/3$ occur. 

Since every vertex of $U$ has either at least%
\COMMENT{$|X| = x(1+\alpha)$ so $(|G| - |X|)/2 \geq (2-\gamma)n/2 - (1+\alpha)x/2 \geq 
n - x/2 - \alpha n \geq x/2+y+z-\alpha n \geq y+z+\alpha n$.} 
$(|G| - |X|)/2 \geq y + z + \alpha n$
inneighbours outside $X$ or at least $y+z+\alpha n$
outneighbours outside $X$, we may choose a set $U' \subseteq U$ of size $|U'| 
\geq |U|/2 \geq \alpha \beta n/10$ such that either 
\begin{itemize}
\item[($\alpha_1$)] every $v \in U'$ has at least $y + \alpha n$ 
outneighbours outside 
$X$ and at least $y+z+\alpha n$ inneighbours outside $X$, or
\item[($\alpha_2$)] every $v \in U'$ has at least $y+z+\alpha n$ 
outneighbours 
outside $X$ and at least $z + \alpha n$ inneighbours outside $X$.
\end{itemize}

So $G'[X]$ is a $(2\eps/\beta)$-regular $(d/2)$-dense cycle of cluster 
tournaments on clusters $X_1, \dots, X_k$ of size $(1+\alpha)x/k$, and $U' 
\subseteq X_1 \cup \dots \cup X_k$ has size $|U'| \geq \alpha \beta x/10$. 
Also $T_{ext}$ is a directed tree on $x$ vertices rooted at $t_1$ and with 
$\Delta(T_{ext}) \leq \Delta^*$, and $H \subseteq V(T_{ext})$ has $|H| \leq 
n/\Delta^{k^{1/\beta\delta}} \leq \delta x/7k$ and $|\{t \in T_{ext}: 1 \leq 
d(t, \mathcal{P}_k(H)) \leq k^3\}| \leq \delta \beta n \leq \delta x$. So by 
Lemma~\ref{embed_with_restrictions} (with $\alpha \beta/10$, $\Delta^*$ and 
$d/2$ in place of $\lambda$, $\Delta$ and $d$ respectively), $G'[X]$ contains 
a copy of $T_{ext}$ in which every vertex of $H$ is embedded to a vertex of 
$U'$. 

So every vertex $t \in H$ has been embedded to some vertex $v(t) \in U'$. 
Suppose that $(\alpha_1)$ holds. Then for every $t \in H$, $v(t)$ has at 
least $y + 2n/\Delta$ outneighbours outside $X$ (and so unoccupied by 
vertices of $T_{ext}$). Since the only vertices of $T_{ext}$ which may have 
neighbours in~$T_1 - T_{ext}$ are the vertices of $H$, we may use 
Theorem~\ref{bestsofar} to extend the embedding of $T_{ext}$ in $G[X]$ to an 
embedding of $T_1$ in $G$ in the same way as in the proof of 
Lemma~\ref{onebyone} (we cannot just apply Lemma~\ref{onebyone} as vertices 
of $G$ to which we embedded $T_{ext} - H$ may not have sufficiently many 
outneighbours, but since vertices of $T_{ext} - H$ do not have any 
outneighbours outside $T_{ext}$ this does not cause any problems). Then for 
every $t \in H$, $v(t)$ has at least $z+2n/\Delta$ inneighbours outside $X$ 
which are not occupied by this embedding of $T_1$. So in the same way we may 
extend the embedding of $T_{ext}$ in $G[X]$ to an embedding of $T_2$ in the 
vertices of $G$ not occupied by $T_1 - T_{ext}$. So the embeddings of $T_1$ 
and $T_2$ only overlap in $T_{ext}$, and so together form an embedding of $T$ 
in $G$. If instead $(\alpha_2)$ holds we may embed $T$ in $G$ similarly by 
first embedding $T_2$ then $T_1$.  
\endproof

We can now deduce that if $G$ is a large almost-regular tournament and 
if $|T_\Delta|>1$, then Sumner's conjecture holds with a little room to spare
(we shall need this extra room in the proof of Lemmas~\ref{smallcorepart1} and~\ref{smallcorepart1.5}).
Indeed, we shall see that a 
large almost-regular tournament $G$ is also a robust outexpander, and so if 
$T_\Delta$ is large, then we can embed $T$ in $G$ by 
Lemma~\ref{ar_tourns_large_core_tree}. On the other hand, if $T_\Delta$ is 
small but has more than one vertex,
then we may embed $T$ in $G$ by Lemma~\ref{ar_tourns_small_core_tree}. 

In particular, together with Lemma~\ref{core_tree_size_1} (which deals with the case $|T_\Delta|=1$), 
this means that at this stage, we have proved that Sumner's 
conjecture holds for all large almost-regular tournaments.
\begin{lemma} \label{sumner_for_ar_tourns}
Suppose that $1/n \ll \gamma \ll 1/\Delta \ll 1$. Let $T$ be a directed tree on $n$ 
vertices with $|T_\Delta| > 1$.  Then every $\gamma$-almost-regular tournament 
$G$ on at least $(2-\gamma)n$ vertices contains a copy of $T$. 
\end{lemma}

\proof
Introduce constants $\mu, \nu, \eta, \Delta', \beta, \gamma'$ such that 
$$1/n \ll 1/\Delta' \ll \mu \ll \nu \ll \eta \ll \gamma \ll \beta \ll \gamma' 
\ll 1/\Delta \ll 1.$$ 
Let $G$ be a $\gamma$-almost-regular tournament on at least 
$(2-\gamma)n$ vertices. Then we shall show that $G$ is a robust $(\mu, 
\nu)$-outexpander. Indeed, let $S \subseteq V(G)$ satisfy $\nu |G| \leq |S| 
\leq 2|G|/3$. Then at least $(1-\gamma)|S|(|G|-1)/2$ edges originate in $S$. 
At most $\binom{|S|}{2}$ of these have both endvertices in $S$, so at least 
$(1-\gamma)|S|(|G|-1)/2 - \binom{|S|}{2} \geq |S|((1-\gamma)(|G|-1) - |S|)/2 
\geq \nu |G|^2/10$ edges leave $S$. So at least $\nu |G|/20 \geq 3 \mu |G|$ 
vertices outside $S$ have at least $\nu |G|/20 \geq 3 \mu |G|$ inneighbours 
in $S$. At most $2 \mu |G|$ vertices of $S$ have fewer than $\mu |G|$ 
inneighbours in $S$, and so $|RN^+_\mu(S)| \geq |S| + \mu |G|$, as desired. 
On the other hand, if $S \subseteq V(G)$ satisfies $2|G|/3 < |S| \leq 
(1-\nu)|G|$, every vertex of $G$ has at least $|G|/7 \geq \mu |G|$ 
inneighbours in $S$. So $|RN^+_\mu(S)| = |G| \geq |S| + \mu |G|$, as desired.  

So $G$ is indeed a robust $(\mu, \nu)$-outexpander. Clearly $\delta^0(G) \geq 
\eta |G|$. So if $|T_{\Delta'}| \geq \beta n$, then by 
Lemma~\ref{ar_tourns_large_core_tree}, $G$ contains a copy of $T$. So we may 
assume that $|T_{\Delta'}| \leq \beta n$. But $G$ is also a 
$\gamma'$-almost-regular tournament on at least $(2-\gamma')n$ vertices, and 
so by Lemma~\ref{ar_tourns_small_core_tree}, $G$ contains a copy of $T$. 
\endproof

\section{Embedding trees whose core tree is small} \label{sec:smallcore}

We now turn our attention to the general case of the problem. As when 
considering almost-regular tournaments, we consider the problem of embedding 
directed trees whose core trees are small separately from the case when the 
core trees are large. In this section we shall consider directed trees with 
small core trees, proving the following lemma. 

\begin{lemma} \label{main_small_core}
Suppose $1/n \ll \beta, 1/\Delta' \ll 1$. Let $T$ be a directed tree on $n$ 
vertices with $|T_{\Delta'}| \leq \beta n$, and let $G$ be a tournament on 
$2n-2$ vertices. Then $G$ contains a copy of $T$. 
\end{lemma}

We begin by showing that we may assume that the tournament $G$ consists of 
two large disjoint almost-regular tournaments, with almost all of the edges 
between them directed the same way. 

\begin{lemma} \label{smallcorepart1}
Suppose that $1/n \ll \beta, 1/\Delta \ll \gamma \ll \eta \ll 1$. Let $T$ be 
a directed tree on $n$ vertices with $|T_\Delta| \leq \beta n$, and let $G$ 
be a tournament on $2n-2$ vertices. Let $y$ be the outweight of $T_\Delta$, 
and let $z$ be the inweight of $T_\Delta$. Then the following properties 
hold. 
\begin{enumerate}
\item[(i)] If $z < \eta n$ or $y < \eta n$ then $G$ contains a copy of $T$.
\item[(ii)] Either $G$ contains a copy of $T$, or we can find disjoint sets 
$Y, Z \subseteq V(G)$ such that $|Y| \geq (2-\gamma)y$ and $|Z| \geq 
(2-\gamma)z$, $G[Y]$ and $G[Z]$ are $\gamma$-almost-regular, any vertex of 
$Y$ has at most $3 \gamma n$ outneighbours in $Z$ and any vertex of $Z$ has 
at most $3 \gamma n$ inneighbours in $Y$.   
\end{enumerate}
\end{lemma}

\proof
Introduce new constants $M, M', \eps, \eps', \alpha, \gamma^*$ and $\Delta^*$ 
such that 
$$1/n \ll \beta, 1/\Delta \ll 1/M \ll 1/M' \ll \eps \ll \eps' \ll \gamma \ll 
\alpha \ll \eta \ll \gamma^* \ll 1/\Delta^* \ll 1.$$ 
Partition the vertex set of $G$ into sets $A, B, C, D, E$ such that:
\begin{align*}
A &\subseteq  \{v \in G: d^+(v) \leq y+\eps n \},
\\B &\subseteq \{v \in G: y+\eps n < d^+(v) < n-\eps n \},
\\C &\subseteq \{v \in G: d^+(v), d^-(v) \geq n-\eps n\}, 
\\D &\subseteq \{v \in G: z+ \eps n < d^-(v) < n-\eps n\}, 
\\E &\subseteq \{v \in G: d^-(v) \leq z+\eps n \}.
\end{align*}
These subset relations may not all be equality, for example in the case where 
$z$ is very small, when we have $y + \eps n \geq n - \eps n$. However, it is 
clear that each vertex $v \in G$ lies in at least one of these five sets, so 
we may choose such a partition of $V(G)$. Let $x := |T_\Delta|$, so $x+y+z = 
n$ and $x \leq \beta n$. 

Suppose that $|B| \geq 3x$. Then by Theorem~\ref{bestsofar} we may 
embed~$T_\Delta$ in $G[B]$. 
Let $S_\Delta \subseteq B$ be the set of vertices occupied by this embedding 
of $T_\Delta$. Then every vertex of $S_\Delta$ has at least $y+\eps n-x \geq 
y + 2n/\Delta$ outneighbours outside $S_\Delta$ and at least $|G|-x-(n-\eps 
n) \geq y+z+2n/\Delta$ inneighbours outside $S_\Delta$. Let $T_1$ be the 
subtree of $T$ formed by $T_\Delta$ and all outcomponents of $T_\Delta$, and 
let $T_2$ be the subtree of $T$ formed by $T_\Delta$ and all incomponents of 
$T_\Delta$. Then $|T_1| = x+y$ and $|T_2| = x+z$. By 
Proposition~\ref{coretreeprops}(iv), all incomponents and outcomponents 
of~$T_\Delta$ contain at most $n/\Delta$ vertices, so by 
Lemma~\ref{onebyone}(c) we may extend our embedding of~$T_\Delta$ 
in~$S_\Delta$ to an embedding of $T_1$ in $G$. Then each vertex of $S_\Delta$ 
still has at least $z+2n/\Delta$ inneighbours outside $S_\Delta$ which are 
not occupied by this embedding of $T_1$, so by Lemma~\ref{onebyone}(c) we may 
also extend our embedding of $T_\Delta$ in $S_\Delta$ to an embedding of 
$T_2$ in $G$ which avoids vertices occupied by the embedding of $T_1 - 
T_\Delta$. Then these embeddings of $T_1$ and $T_2$ do not overlap outside 
$T_\Delta$, and so together form an embedding of $T$ in $G$. We may therefore 
assume that $|B| < 3x \leq 3\beta n$. By the same argument (embedding first 
$T_2$ and then $T_1$ in $G$) we may assume that $|D| < 3x \leq 3 \beta n$. 

If $|T_{\Delta^*}| = 1$, then $G$ contains a copy of $T$ by 
Lemma~\ref{core_tree_size_1}. So we may assume that $|T_{\Delta^*}| \geq 2$. 
Now, if $z < \eta n$, then every $v \in E$ satisfies $d^-(v) < (\eta+\eps) n 
< 2 \eta n$, so $|E| \leq 4 \eta n+1$, and so $|B \cup D \cup E| \leq 4 \eta 
n +1+ 6\beta n \leq 5 \eta n$. Let $G' := G[A \cup C]$. Then $|G'| \geq 
2n-2-5\eta n$, and every vertex $v \in G'$ has $d^+(v) \leq n+\eps n$. So by 
Proposition~\ref{findartourn}, $G'$ contains a $\gamma^*$-almost-regular 
subtournament $G'$ on at least $(2-\gamma^*)n$ vertices. Since 
$|T_{\Delta^*}| \geq 2$, by Lemma~\ref{sumner_for_ar_tourns} $G'$ contains a 
copy of $T$, so $G$ contains a copy of $T$ also. If instead we have $y < \eta 
n$, then we may similarly embed $T$ in $G[C \cup E]$. So if $z < \eta n$ or 
$y < \eta n$ then $G$ contains a copy of $T$, completing the proof of (i). So 
for (ii), we may assume that $y, z \geq \eta n$.   

Suppose now that $|C| \geq 5 \eps' n$. Let disjoint subsets $V_1, \dots, V_k$ 
and a subgraph $G^* \subseteq G$ satisfy the conditions of 
Corollary~\ref{clustertourns}. So $M' \leq k \leq M$, and $G^*$ is an 
$\eps$-regular cluster tournament on clusters $V_1, \dots, V_k$ of equal size 
$m$, where 
$$\frac{(1-\eps)|G|}{k} \leq m \leq \frac{|G|}{k}.$$
We shall show that $G^*$ has the property that for some $i \in [k]$ we have 
\begin{equation} \label{eq:goodchoiceofi}
\sum_{j \in ([k] \sm \{i\})} d_{ij} \geq \frac{(1 - 3\eps') k}{2} 
\textrm{\hspace{5mm} and} \sum_{j \in ([k] \sm \{i\})} d_{ji} \geq \frac{(1 - 
3\eps') k}{2}. 
\end{equation}
Indeed, if for some $i \in [k]$ we have $\sum_{j \in ([k] \sm \{i\})} d_{ij} 
< (1 - 3\eps') k/2$, then by Lemma~\ref{regularnbrhoods} all but at most 
$\eps' m$ vertices of $V_i$ have at most  
$$\sum_{j \in ([k] \sm \{i\})} d_{ij}m + \eps'km < \frac{(1 - \eps') km}{2} < 
n-8\eps n$$  
outneighbours in $\bigcup_{j \in ([k] \sm \{i\})} V_j$ (in the graph $G^*$), 
and hence at most $n - 8\eps n + (|G| - |G^*|) + |V_i| + 2\eps |G| < n-\eps 
n$ outneighbours in $G$. So at most $\eps' m$ vertices of $V_i$ lie in $C$. 
Similarly if for some $i \in [k]$ we have  $\sum_{j \in ([k] \sm \{i\})} 
d_{ji} < (1 - 3\eps') k/2$ then again at most $\eps'm$ vertices of $V_i$ lie 
in $C$. Since $|C| \geq 5\eps' n > 2 \eps' mk + (|G| - |G^*|)$, there must be 
some $i \in [k]$ which satisfies (\ref{eq:goodchoiceofi}). Fix such an~$i$. 
Then if at least $\alpha k$ values of $j \in [k] \sm \{i\}$ have $d_{ij} \geq 
\alpha$ and $d_{ji} \geq \alpha$ then $G^*$ contains a copy of~$T$ by 
Lemma~\ref{ar_tourns_part_one} (applied with $\eps'$ in the place of 
$\gamma$). Alternatively, if at most $\alpha k$ values of $j \in [k] \sm 
\{i\}$ have $d_{ij} \geq \alpha$ and $d_{ji} \geq \alpha$ then since $y, z 
\geq \eta n$,~$G^*$ contains a copy of $T$ by 
Lemma~\ref{ar_tourns_part_two}(iii) (again applied with $\eps'$ in the place 
of $\gamma$). So in either case $G$ contains a copy of $T$, and so we may 
assume that $|C| < 5 \eps' n$. 

So to prove (ii), observe that we must therefore have $|B \cup C \cup D| \leq 
5\eps' n + 6\beta n \leq 6 \eps' n$. Trivially $|A| \leq 2y+2\eps n+1$ and 
$|E| \leq 2z+2\eps n+1$, and so we must have 
\begin{align*}
|A| &\geq 2n-2-6\eps' n - 2z - 2\eps n-1 \geq 2y-7\eps' n, \textrm{ and}
\\|E| &\geq 2n-2-6\eps' n - 2y - 2\eps n-1 \geq 2z-7\eps' n.
\end{align*}
So by Proposition~\ref{findartourn}, $G[A]$ contains a 
$\gamma$-almost-regular subtournament on at least $(2-\gamma)y$ vertices, and 
$G[E]$ contains a $\gamma$-almost-regular subtournament on at least 
$(2-\gamma)z$ vertices. Let $Y$ and $Z$ be the vertex sets of these 
subtournaments respectively. 
Then any vertex of $Y$ has at least $(1-2\gamma) y$ outneighbours in $Y$, and 
so has at most $y + \eps n - (1-2\gamma) y \leq 3 \gamma n$ outneighbours in 
$Z$. Similarly any vertex of $Z$ has at least $(1-2\gamma) z$ inneighbours in 
$Z$, and so has at most $3 \gamma n$ inneighbours in $Y$. So $Y$ and $Z$ are 
as required for (ii). 
\endproof

The next lemma builds on the previous lemma and will in turn be used in the 
proof of Lemma~\ref{smallcorepart2}. 

\begin{lemma} \label{smallcorepart1.5}
Suppose that $1/n \ll \beta, 1/\Delta' \ll \alpha \ll 1/\Delta \ll 1$. Let 
$T$ be a directed tree on $n$ vertices with $|T_{\Delta'}| \leq \beta n$. Let $y$ and 
$z$ be the outweight and inweight of $T_{\Delta'}$ respectively. Suppose that 
forests~$F^-$ and $F^+$ are induced subgraphs of $T$ which partition the 
vertices of $T$, such that $|F^+| \leq y + 2\alpha n $, $|F^-| \leq z - 
\alpha n$, and every edge of $T$ between $F^-$ and $F^+$ is directed from 
$F^-$ to $F^+$. Suppose also that either 
\begin{itemize}
\item[(i)] no component of $F^+$ has order greater than $y-\alpha n$, or
\item[(ii)] the largest component $T_1$ of $F^+$ has $|(T_1)_{\Delta}| \geq 2$. 
\end{itemize} 
Then any tournament $G$ on $2n-2$ vertices contains a copy of $T$.
\end{lemma}

\proof
Let $G$ be a tournament on $2n-2$ vertices, and let $T_1$ and $T_2$ be the 
largest and second largest components of $F^+$ respectively. Introduce new 
constants $\gamma$ and $\eta$ with  
$$1/n \ll \beta, 1/\Delta' \ll \gamma \ll \alpha \ll 1/\Delta \ll \eta \ll 
1.$$  
Then by Lemma~\ref{smallcorepart1} we may assume that $y, z \geq \eta n$. 
Also by Lemma~\ref{smallcorepart1} we may find subsets $Y, Z \subseteq V(G)$ 
such that $|Y| \geq (2-\gamma)y$, $|Z| \geq (2-\gamma)z$, $G[Y]$ is 
$\gamma$-almost-regular, each vertex of $Y$ has at most $3\gamma n$ 
outneighbours in $Z$, and each vertex of $Z$ has at most $3\gamma n$ 
inneighbours in $Y$. Then $|Y| \geq 3|F^+|/2 + \alpha n \geq |F^+| + |T_2| + 
\alpha n$, and $|Z| \geq 2|F^-| + \alpha n$, and so by 
Lemma~\ref{component_by_component} any embedding of $T_1$ in $G[Y]$ may be 
extended to an embedding of $T$ in $G$. 

It therefore suffices to embed $T_1$ in $G[Y]$. If $|T_1| < y/2$ then we may 
do this by Theorem~\ref{bestsofar}. If instead $|T_1| \geq y/2 \geq \eta n/2$ 
and we also have (i), then $|T_1| \leq y - \alpha n$. Since $|Y| \geq 
(2-\gamma)y \geq 2|T_1| + \alpha n$ we may embed $T_1$ in $G[Y]$ by 
Theorem~\ref{approxversion}(i). Finally, if $|T_1| \geq \eta n/2$ and we also 
have (ii), then $|T_1| \leq |F^+| \leq y+2\alpha n$ and $|(T_1)_\Delta| \geq 
2$. Since $\gamma \leq 9\alpha/\eta$, $G[Y]$ is a 
$9\alpha/\eta$-almost-regular tournament on at least $(2 - \gamma)y \geq 
(2-9\alpha/\eta)|T_1|$ vertices, and so we may embed $T_1$ in $G[Y]$ by 
Lemma~\ref{sumner_for_ar_tourns}. So in any case we may embed $T_1$ in 
$G[Y]$, completing the proof. 
\endproof

Observe that as with Lemma~\ref{component_by_component} a `dual' form of 
Lemma~\ref{smallcorepart1.5} can be proved similarly. For this we instead 
require that $|F^+| \leq y - \alpha n$ and $|F^-| \leq z + 2\alpha n$, and 
also either that no component of $F^-$ has order greater than $z - \alpha n$ 
or that the largest component $T_1$ of $F^-$ has $|(T_1)_\Delta| \geq 2$. If 
these conditions are met then we may conclude that $G$ contains a copy 
of~$T$. As with Lemma~\ref{component_by_component}, we shall sometimes 
implicitly refer to this `dual' when referring to 
Lemma~\ref{smallcorepart1.5}.

In the next lemma we show that Lemma~\ref{main_small_core} holds for any 
directed tree $T$ whose core tree $T_\Delta$ is not a directed path in which 
most of the outweight and inweight of $T_\Delta$ lies at the endvertices of 
$T_\Delta$. We say that a vertex $t$ of a directed tree $T$ is an 
\emph{outleaf} if $t$ has one inneighbour and no outneighbours, or an 
\emph{inleaf} if $t$ has one outneighbour and no inneighbours.

\begin{lemma} \label{smallcorepart2}
Suppose that $1/n \ll \beta, 1/\Delta' \ll 1/\Delta \ll \sigma \ll 1$. Let 
$T$ be a directed tree on $n$ vertices with $|T_{\Delta'}| \leq \beta n$, and let $y$ 
and $z$ be the outweight and inweight of $T_{\Delta'}$ respectively. Let $G$ 
be a tournament on $2n-2$ vertices. Then either $G$ contains a copy of $T$, 
or $T_{\Delta}$ is a directed path whose outleaf has outweight at least $y - 
\sigma n$ and whose inleaf has inweight at least $z-\sigma n$. 
\end{lemma}

\proof 
Introduce new constants $\alpha$ and $\eta$ with 
$$1/n \ll \beta, 1/\Delta' \ll \alpha \ll 1/\Delta \ll \sigma \ll \eta \ll 1.$$ 
Then by Lemma~\ref{smallcorepart1} we may assume that $y, z \geq \eta n$. 
Also, if $|T_\Delta| = 1$ then $G$ contains a copy of $T$ by 
Lemma~\ref{core_tree_size_1}, so we may assume that $|T_{\Delta}| \geq 2$. 

Suppose that some vertex $t \in T$ has the property that $w^-(t) \leq z - 
\alpha n-1$, and also that every outcomponent of $t$ contains at most $w^+(t) 
- 3\alpha n = |V^+| - 3 \alpha n$ vertices. Then let the set $V^-$ consist of 
$t$ and every vertex in an incomponent of $t$, and let $V^+ := V(T) \sm V^-$. 
Then $|V^-| \leq w^-(t) + 1 \leq z - \alpha n$, and every edge of $T$ between 
$V^-$ and $V^+$ is directed from $V^-$ to $V^+$. Also, each component of 
$T[V^+]$ contains at most $w^+(t) - 3\alpha n$ vertices. Now, select a source 
vertex from the largest component of $T[V^+]$, delete this vertex from $V^+$, 
and add it to $V^-$. Repeat this step until we have $|V^+| \leq y+2\alpha n$ 
and $|V^-| \leq z - \alpha n$. For these final $V^+$ and $V^-$, let $F^+ := 
T[V^+]$ and let $F^- := T[V^-]$. Then $F^-$ and $F^+$ are forests which 
partition the vertices of $T$, with $|F^+| \leq y+2\alpha n$ and $|F^-| \leq 
z-\alpha n$. Also, every edge of $T$ between $F^-$ and $F^+$ is directed 
from~$F^-$ to $F^+$. Finally, since we always deleted a vertex from the 
largest component of $T[V^+]$, no component of $F^+$ contains more than 
$|F^+| - 3\alpha n \leq y-\alpha n$ vertices. So by 
Lemma~\ref{smallcorepart1.5}(i) $G$ contains a copy of $T$.  
So we may assume that

\textno there is no vertex $t \in T$ such that $w^-(t) \leq z - \alpha n-1$ 
and every outcomponent of $t$ contains at most $w^+(t) - 3\alpha n$ vertices. 
In particular, this implies that for every inleaf $t$ of $T_\Delta$, at least 
$n/2\Delta$ vertices of $T$ lie in incomponents of $t$.& (\dagger)  

Indeed, if $T_\Delta$ contains some inleaf $t$ such that fewer than 
$n/2\Delta \leq z - \alpha n-1$ vertices of~$T$ lie in incomponents of $t$, 
then by the definition of $T_\Delta$ at least $n/2\Delta-1$ vertices of $T$ 
lie in outcomponents of $t$ other than the outcomponent containing the 
remaining vertices of~$T_\Delta$. Moreover, the definition of $T_\Delta$ also 
implies that at least $n/\Delta$ vertices of $T$ lie in the one component of 
$T - t$ containing $T_\Delta-t$. Altogether this shows that every 
outcomponent of $t$ contains at most $w^+(t) - n/2\Delta+1 \leq w^+(t) - 
3\alpha n$ vertices, a contradiction. By the same argument with the roles of 
incomponents and outcomponents switched, we may assume that 

\textno there is no vertex $t \in T$ such that $w^+(t) \leq y - \alpha n-1$ 
and every incomponent of $t$ contains at most $w^-(t) - 3\alpha n$ vertices. 
It follows from this that for every outleaf $t$ of $T_\Delta$, at least 
$n/2\Delta$ vertices of $T$ lie in outcomponents of $t$. & (\dagger \dagger) 

\medskip \noindent {\bf Claim.}
\emph{If $T_\Delta$ has at least two inleaves or at least two outleaves, then 
$G$ contains a copy of~$T$.} 

\medskip \noindent
To prove the claim, suppose that $T_\Delta$ has two outleaves $t$ and $t'$ 
(the proof for inleaves is similar). Then we shall form a set $V^+$ of size 
between $n-z+\alpha n$ and $y+2\alpha n$ such that any edge of $T$ between 
$V^+$ and $V^- := V(G) \sm V^+$ is directed from $V^-$ to $V^+$. We may do 
this by repeatedly selecting a sink vertex of $T$, adding it to $V^+$ and 
removing it from $T$. Now, by ($\dagger \dagger$) at least $n/2\Delta$ 
vertices lie in outcomponents of $t$, and at least $n/2\Delta$ vertices lie 
in outcomponents of $t'$. Furthermore, if $T'$ is an outcomponent of $t$, 
then any sink vertex in $T'$ is a sink vertex in $T$, and the same is true if 
$T'$ is instead an outcomponent of $t'$. So we may form $V^+$ and $V^-$ as 
described above so that additionally $V^+$ contains at least $n/2\Delta$ 
vertices from outcomponents of $t$ and at least $n/2\Delta$ vertices from 
outcomponents of $t'$. Fix such a choice of $V^+$ and $V^-$, and let $F^+ := 
T[V^+]$ and $F^- := T[V^-]$ be the induced forests. Then $|F^+| \leq y + 
2\alpha n$ and $|F^-| = n - |F^+| \leq z - \alpha n$, and every edge of $T$ 
between $F^-$ and $F^+$ is directed from $F^-$ to $F^+$. So if every 
component of $F^+$ contains at most $y - \alpha n$ vertices, then $G$ 
contains a copy of $T$ by Lemma~\ref{smallcorepart1.5}(i). We may therefore 
assume that the largest component $T^+$ of $F^+$ contains more than $y- 
\alpha n \geq |F^+| - n/4\Delta$ vertices. Since $F^+$ includes at least 
$n/2\Delta$ vertices from outcomponents of $t$ and at least $n/2\Delta$ 
vertices from outcomponents of $t'$, it follows that $T^+$ contains at least 
$n/4\Delta$ vertices from outcomponents of $t$ and at least $n/4\Delta$ 
vertices from outcomponents of $t'$. As a consequence $T^+$ must contain $t$ 
and $t'$. Furthermore, we must have $t, t' \in (T^+)_{4\Delta}$, and so 
$|(T^+)_{4\Delta}| \geq 2$. So $G$ contains a copy of $T$ by 
Lemma~\ref{smallcorepart1.5}(ii), which proves the claim. \medskip 

We may therefore assume that $T_{\Delta}$ has at most one outleaf and at most 
one inleaf. So $T_{\Delta}$ is a path with one inleaf and one outleaf. Let 
$t_1, \dots, t_x$ be the vertices of this path, labelled so that $t_1$ is the 
inleaf of $T_{\Delta}$ (so $t_1 \rightarrow t_2$), $t_x$ is the outleaf of 
$T_{\Delta}$ (so $t_{x-1} \rightarrow t_x$), and for each $i \in [x-1]$ there 
is an edge of $T_{\Delta}$ between $t_i$ and $t_{i+1}$.

Now suppose that the inweight of $T_\Delta$ is less than $z-2\alpha n$. Let 
the set $V^-$ consist of all vertices of $T$ which lie in $T_\Delta$ or in 
incomponents of $T_\Delta$. Then $|V^-| \leq z -2 \alpha n + |T_\Delta| \leq 
z - \alpha n$ (since $|T_\Delta| \leq |T_{\Delta'}| \leq \beta n$). Also, 
every edge of $T$ between $V^-$ and $V^+ := V(T) \sm V^-$ is directed from 
$V^-$ to $V^+$. Choose a source vertex of $T[V^+]$, delete it from $V^+$, and 
add it to $V^-$, and repeat this step until we have $|V^-| \leq z-\alpha n$ 
and $|V^+| \leq y+2\alpha n$. For these final $V^-$ and $V^+$, let $F^+ := 
T[V^+]$ and $F^- := T[V^-]$ be the induced forests. Then $|F^-| \leq z - 
\alpha n$, $|F^+| \leq y + 2\alpha n$, and every edge of $T$ between $F^-$ 
and $F^+$ is directed from $F^-$ to $F^+$. Also, every component of $F^+$ is 
contained within a component of $T - T_\Delta$, and so has order at most 
$n/\Delta \leq y - \alpha n$ by Proposition~\ref{coretreeprops}. So $G$ 
contains a copy of $T$ by Lemma~\ref{smallcorepart1.5}(i). We may therefore 
assume that the inweight of $T_\Delta$ is at least $z - 2\alpha n$, and by a 
similar argument we may also assume that the outweight of $T_\Delta$ is at 
least $y - 2\alpha n$. It follows that the outweight of $T_\Delta$ is at most 
$n - (z - 2\alpha n) \leq y+3\alpha n$ and that the inweight of $T_\Delta$ is 
at most $n - (y - 2\alpha n) \leq z + 3\alpha n$. 

We now suppose that fewer than $y-\sigma n$ vertices of $T$ lie in 
outcomponents of $t_x$. Let $T_1$ be the subtree of $T$ formed by $T_\Delta$ 
and all of its outcomponents. Initially let the set $V^+ := V(T_1)$, so 
$|V^+| \leq y + 4\alpha n$, and every edge of $T$ between $V^+$ and $V^- := 
V(G) \sm V^+$ is directed from $V^-$ to $V^+$. Choose a sink vertex of 
$T[V^-]$, delete it from $V^-$ and add it to $V^+$, and repeat this step 
until we have $|V^+| \leq y + 4\alpha n$ and $|V^-| \leq z - 2\alpha n$. Fix 
these final $V^+$ and $V^-$ and let $F^- := T[V^-]$ and $F^+ := T[V^+]$ be 
the induced forests. So $|F^+| \leq y + 4\alpha n$, $|F^-| \leq z - 2\alpha 
n$, and every edge of $T$ between $F^-$ and $F^+$ is directed from $F^-$ to 
$F^+$. Also $T_1 \subseteq F^+$, so $T_1$ is contained within a single 
component $T^+$ of $F^+$. Since at least $y-2\alpha n$ vertices of $T$ lie in 
outcomponents of $T_\Delta$, at least~$\sigma n/2$ vertices of $T$ lie in 
outcomponents of $T_{\Delta}$ other than the outcomponents of $t_x$. 
Moreover, since $t_x$ is an outleaf of $T_\Delta$, by $(\dagger \dagger)$ at 
least $n/2\Delta$ vertices lie in outcomponents of $t_x$. So $t_{x-1} \in 
(T^+)_{2\Delta}$ and $t_x \in (T^+)_{2\Delta}$, and so $|(T^+)_{2\Delta}| \geq 
2$. But since the outweight of $T_\Delta$ is at least $y - 2\alpha n$ we have 
$|T^+| \geq |T_1| \geq y-2\alpha n$, and so $T^+$ must be the largest 
component of $F^+$. So $G$ contains a copy of $T$ by 
Lemma~\ref{smallcorepart1.5}(ii). 

So we may assume that at least $y-\sigma n$ vertices of $T$ lie in 
outcomponents of $t_x$, as desired. If fewer than $z-\sigma n$ vertices of 
$T$ lie in incomponents of $t_1$, then we may similarly embed $T$ in $G$, so 
we may also assume that at least $z - \sigma n$ vertices of $T$ lie in 
incomponents of $t_1$. So at most $3 \sigma n$ vertices of $T$ do not lie in 
incomponents of $t_1$ or outcomponents of $t_x$. It remains only to show that 
$T_{\Delta}$ is a directed path. So suppose for a contradiction that 
$T_{\Delta}$ is not a directed path. Then there is some $i \in [x-1]$ such 
that $t_i \leftarrow t_{i+1}$. Choose the minimal such $i$ (note $i > 1$ as 
$t_1$ is an inleaf of $T_\Delta$). Then $t_i$ has two inneighbours and no 
outneighbours in $T_\Delta$. So at least two incomponents of $t_i$ contain at 
least $n/\Delta$ vertices, and so no incomponent of $t_i$ contains more than 
$w^-(t_i) - n/\Delta \leq w^-(t_i) - 3 \alpha n$ vertices. Also, at most $3 
\sigma n \leq y - \alpha n -1$ vertices of $T$ lie in outcomponents of $t_i$, 
contradicting $(\dagger \dagger)$. 
\endproof

We can now prove that Sumner's universal tournament conjecture holds for any 
large directed tree $T$ whose core tree $T_\Delta$ contains precisely two 
vertices.

\begin{lemma} \label{core_tree_size_2}
Suppose that $1/n \ll 1/\Delta' \ll 1$. Let $T$ be a directed tree on $n$ vertices 
with $|T_{\Delta'}| = 2$, and let $G$ be a tournament on $2n-2$ vertices. 
Then $G$ contains a copy of $T$. 
\end{lemma}

\proof
Introduce new constants $\Delta, \eps, \gamma$ and $\eta$ with $$1/n \ll 
\beta, 1/\Delta' \ll 1/\Delta \ll \eps \ll \gamma \ll \eta \ll 1.$$ 
Then $|T_{\Delta'}| = 2 \leq \beta n$. Also, since $\Delta \leq \Delta'$ we 
have $T_{\Delta} \subseteq T_{\Delta'}$. If $|T_\Delta| = 1$, then by 
Lemma~\ref{core_tree_size_1} $G$ contains a copy of $T$. So we may assume 
that $T_\Delta = T_{\Delta'}$. Let $t_2$ and $t_1$ be the vertices of 
$T_\Delta$, labelled so that $t_2 \rightarrow t_1$. Let $y$ be the outweight 
of $T_{\Delta}$, and let $z$ be the inweight of $T_{\Delta}$, so $y + z = n - 
2$. Then by Lemma~\ref{smallcorepart2} (with $\eps$ in the place of 
$\sigma$), we may assume that $t_2$ has inweight at least $z - \eps n$, and 
also that $t_1$ has outweight at least $y - \eps n$. Let $T_1$ be the subtree 
of $T$ consisting of all vertices which lie in $T_\Delta$ or in outcomponents 
of $T_\Delta$, and let $T_2$ be the subtree of $T$ consisting of all vertices 
which lie in $T_\Delta$ or in incomponents of $T_\Delta$. So $|T_1| = y+2$ 
and $|T_2| = z+2$. By Lemma~\ref{smallcorepart1}(i) we may assume that $y, z 
\geq \eta n$. 

As in the proof of Lemma~\ref{smallcorepart1}, we partition the vertices of 
$G$ into sets $A, B, C, D$ and $E$, where: 
\begin{align*}
A &:= \{v \in G: d^+(v) \leq y+\eps n \},
\\B &:= \{v \in G: y+\eps n < d^+(v) < n-\eps n \},
\\C &:= \{v \in G: d^+(v), d^-(v) \geq n-\eps n\},
\\D &:= \{v \in G: z+ \eps n < d^-(v) < n-\eps n\},
\\E &:= \{v \in G: d^-(v) \leq z+\eps n \}.
\end{align*}
Since $y, z \geq \eta n$ and $\eps \ll \eta$ this is indeed a partition.
Suppose first that $|B| \geq 2$. Then we may embed $T_\Delta$ in $G[B]$. Let 
$S_\Delta \subseteq B$ be the set of vertices occupied by $T_\Delta$. Then 
every vertex of $S_\Delta$ has at least $y+\eps n -1 \geq y + 2n/\Delta$ 
outneighbours outside $S_\Delta$ and at least%
\COMMENT{I think the $-2$ here 
is correct: $d^-(v) \geq |G|-1-d^+(v) > |G| - 1 - (n-\eps n)$, then take off 
one vertex which may lie in $S_\Delta$.} $|G| - 2 - (n-\eps n) \geq 
y+z+2n/\Delta$ inneighbours outside $S_\Delta$. So by Lemma~\ref{onebyone}(c) 
we may extend the embedding of $T_\Delta$ in $S_\Delta$ to an embedding of 
$T_1$ in $G$. This embedding of $T_1$ occupies at most $y$ vertices of $G$ 
outside $S_\Delta$, and so we may apply Lemma~\ref{onebyone}(c) again to 
extend the embedding of $T_\Delta$ in $S_\Delta$ to an embedding of $T_2$ in 
$G$ so that the embeddings of $T_1$ and $T_2$ do not overlap outside 
$T_\Delta$. Then together the embeddings of $T_1$ and $T_2$ form an embedding 
of $T$ in $G$. So we may assume that $|B| \leq 1$. If $|D| \geq 2$ we may 
embed $T$ in $G$ in the same way by embedding $T_\Delta$ in $D$ and then 
extending this embedding to embeddings of first $T_2$ and then $T_1$ in $G$ 
which do not overlap outside $T_\Delta$. So we may also assume that $|D| \leq 1$. 

Now suppose that $|C| \geq 3$. Then we may choose vertices $v_2, v_1 \in C$ 
with $v_2 \rightarrow v_1$ and $|N^+(v_1) \cap N^+(v_2)| \geq \eta n \geq 
\eta n/2 + 2n/\Delta$. Embed $t_1$ to $v_1$ and $t_2$ to $v_2$. Then since 
$|N^+(v_1)|, |N^+(v_2)| \geq n-\eps n \geq y+2n/\Delta$, by 
Lemma~\ref{onebyone}(b) and (c) we may extend the embedding of $T_\Delta$ in 
$\{v_1, v_2\}$ to an embedding of $T_1$ in $G$ so that at least $\eta n/2$ 
vertices of $T_1$ are embedded in $N^+(v_1) \cap N^+(v_2)$. Then at 
most $y + 2 - \eta n/2$ vertices of 
$N^-(v_1) \cup N^-(v_2)$ are occupied by this embedding, and so in each of 
$N^-(v_1)$ and $N^-(v_2)$ at least $n-\eps n - (y + 2 - \eta n/2) \geq z + 
2n/\Delta$ vertices remain unoccupied. So by Lemma~\ref{onebyone}(a) and (c) 
we may extend the embedding of $T_\Delta$ in $\{v_1, v_2\}$ to an embedding 
of $T_2$ in $G$ which does not overlap with the embedding of $T_1$ outside 
$T_\Delta$. Then together these embeddings form an embedding of $T$ in $G$. 
So we may assume that $|C| \leq 2$, and hence that $|A \cup E| \geq 2n-6$. 

\medskip \noindent {\bf Claim.} \emph{Either some vertex of $A$ has at least 
$y$ outneighbours in $A \cup B \cup D$ or some vertex of $E$ has at least $z$ 
inneighbours in $B \cup D \cup E$.}  

\medskip \noindent Indeed, suppose for a contradiction that both of these 
statements are false. Then certainly every vertex of $A$ has fewer than $y$ 
outneighbours in $A$ and every vertex of $E$ has fewer than $z$ inneighbours 
in~$E$. So $|A| \leq 2y-1$ and $|E| \leq 2z-1$. Since $y+z = n-2$ and $|A 
\cup E| \geq  2n-6$, we must have $|A| = 2y-1$ and $|E| = 2z-1$, and also 
$|B| = 1, |D| = 1$ and $|C|=2$. Then every vertex of $A$ must have $y-1$ 
outneighbours in $A$, and so no vertex of $A$ can have an outneighbour in $B$ 
or in $D$. Likewise, every vertex of $E$ must have $z-1$ inneighbours in $E$, 
and so no vertex of $E$ can have an inneighbour in $B$ or in $D$. But then if 
we let $b$ be the vertex in $B$ and $d$ be the vertex in $D$ we have $d^+(b) 
= d^+(d) \pm 3$, contradicting the definition of $B$ and $D$. So either some 
vertex of $A$ has at least $y$ outneigbours in $A \cup B \cup D$ or some 
vertex of $E$ has at least $z$ inneighbours in $B \cup D \cup E$. This 
completes the proof of the claim. \medskip 

If some $v \in A$ has at least $y$ outneighbours in $A \cup B \cup D$, then 
we shall embed $T_1$ in $G[A]$ so that we may then embed the incomponents of 
$t_2$ and $t_1$ in the unoccupied vertices of $E$ and $A$ respectively. For 
this, note that $|E| \leq 2(z+\eps n)+1$, so $|A| \geq 2n-2z-2\eps n - 7 \geq 
2y-3\eps n$ (and similarly we have $|E| \geq 2z-3\eps n$). Since every $a \in 
A$ has at most $y+\eps n$ outneighbours in $A$, by 
Proposition~\ref{findartourn} $G[A]$ contains a $\gamma$-almost-regular 
subtournament on at least $(2-\gamma)y$ vertices. Let $Y$ be the vertex set 
of this subtournament. Now,  
$$|(A \cup B \cup D) \sm Y| \leq 2 + (2y + 2\eps n +1) - (2-\gamma)y \leq 
2\gamma y,$$  
so $v$ must have at least $(1-2\gamma)y$ outneighbours in $Y$. Also, since $v 
\in A$ we have $$(1-2\gamma) y \leq |N^+(v) \cap Y| \leq y+\eps n \leq 
(1+2\gamma)y.$$  
So at most $10\gamma y$ vertices of $N^+(v) \cap Y$ have more than 
$(1-3\gamma)y$ outneighbours in $N^+(v) \cap Y$, and at most $10\gamma y$ 
vertices of $N^+(v) \cap Y$ have more than $(1-3\gamma)y$ inneighbours in 
$N^+(v) \cap Y$. Since every vertex of $Y$ has at least $(1-2\gamma)y$ 
inneighbours in $Y$ and at least $(1-2\gamma)y$ outneighbours in $Y$, this 
means that at least $|N^+(v) \cap Y| - 20 \gamma y \geq 3n/{\Delta}$ vertices 
of $N^+(v) \cap Y$ have at least $\gamma y \geq 6n/{\Delta}$ outneighbours in 
$Y \sm N^+(v)$ and at least $6n/{\Delta}$ inneighbours in $Y \sm N^+(v)$. Let 
$T^+$ be the tree formed by $t_1$ and its outcomponents, so $|T^+| \leq y+1$. 
Then every component of $T^+ - t_1$ is a component of $T - T_{\Delta}$ and so 
has order at most $n/\Delta$ by Proposition~\ref{coretreeprops}. So by 
Lemma~\ref{roundtheback} (applied with $N := N^+(v) \cap (A \cup B \cup D)$ 
and $X := Y \sm N^+(v)$), we may embed $T^+$ in $G[A \cup B \cup D]$ so that 
$t_1$ is embedded to $v$ and at most $4 n/\Delta$ vertices are embedded 
outside $N^+(v)$.  

Since $v \in A$ we have $d^+(v) \leq y + \eps n$, and so $v$ has at least 
\begin{equation} \label{eq:innbrsofvinY} 
|Y|-1-(y+\eps n)-4n/\Delta \geq 7 \eps n
\end{equation} inneighbours in $Y$ which are not occupied by the embedding of 
$T^+$. Let $T^*$ be the tree formed by all vertices of $T$ which do not lie 
in outcomponents of $t_1$ or incomponents of $t_2$. Then every edge incident 
to $t_1$ in $T^*$ is directed towards $t_1$. Also, $|T^*| \leq n - (y - \eps 
n) - (z - \eps n) = 2\eps n+2$, so certainly every component of $T^* - t_1$ 
has order at most $2\eps n+1$. Together with (\ref{eq:innbrsofvinY}) and 
Theorem~\ref{bestsofar} this shows that we may extend the embedding of $t_1$ 
in $\{v\}$ to an embedding of $T^*$ in $\{v\} \cup (N^-(v) \cap Y)$ so that 
the embeddings of $T^+$ and $T^*$ only overlap in the vertex $t_1$. Then in 
particular $t_2$ is embedded to some vertex $v_2 \in Y$.  

To complete the embedding, observe that every vertex of~$Y$ has at least 
$(1-2\gamma)y$ outneighbours in $Y$, and therefore at most $3\gamma y$ 
outneighbours outside $Y$. So $v_2$ has at least $|E|-3\gamma y \geq 
z+2n/\Delta$ inneighbours in $E$, none of which have been occupied by the 
embeddings of $T^+$ and $T^*$. Let $T^-$ be the subtree of $T$ consisting of 
$t_2$ and all of its incomponents. Then $|T^-| \leq z+1$, and each component 
of $T^- - t_2$ is a component of $T - T_{\Delta}$ and so has order at most 
$n/{\Delta}$ by Proposition~\ref{coretreeprops}. So by 
Lemma~\ref{onebyone}(c) we may extend the embedding of $t_2$ in $\{v_2\}$ to 
an embedding of $T^-$ in $\{v_2\} \cup E$. These embeddings together form an 
embedding of $T$ in~$G$.  

If instead some $v \in E$ has at least $z$ inneighbours in $B \cup D \cup E$ 
then we may similarly embed~$T$ in $G$ by choosing $Z$ to be the vertex set 
of a $\gamma$-almost-regular subtournament of~$G[E]$ on at least 
$(2-\gamma)z$ vertices and embedding $T^-$ in~$G[B \cup D \cup E]$, then 
embedding $T^* - t_2$ in the unoccupied vertices of $Z$, before finally 
embedding $T^+ - t_1$ in $G[A]$. 
\endproof

We can now give the proof of Lemma~\ref{main_small_core}. It was necessary to 
prove Lemma~\ref{core_tree_size_2} separately from this as the method of 
proof does not hold for $|T_\Delta| = 2$ (we cannot obtain the partition of 
$V(G)$ into $Y^*$ and $Z^*$ in this case). 

\medskip \noindent {\bf Proof of Lemma~\ref{main_small_core}.}
Introduce new constants $\gamma, \alpha, \Delta$ and $\eta$ with 
$$1/n \ll \beta, 1/\Delta' \ll 1/\Delta \ll \gamma \ll \alpha \ll \eta \ll 
1.$$  
Let $y'$ be the outweight of $T_{\Delta'}$ and let $z'$ be the inweight of 
$T_{\Delta'}$. Then by Lemma~\ref{smallcorepart1} we may assume that $y', z' 
\geq \eta n$. Similarly let $y$ and $z$ be the outweight and inweight of 
$T_{\Delta}$ respectively. If $|T_{\Delta}| = 1$, then~$G$ contains a copy of 
$T$ by Lemma~\ref{core_tree_size_1}. If instead $|T_{\Delta}| = 2$ then~$G$ 
contains a copy of~$T$ by Lemma~\ref{core_tree_size_2}. So we may assume that 
$\ell := |T_\Delta| \geq 3$, and by Lemma~\ref{smallcorepart2} we may assume 
that $T_\Delta$ is a directed path. Let $t_1, \dots, t_\ell$ be the vertices 
of $T_\Delta$, labelled so that $t_i \rightarrow t_{i+1}$ for each $i \in 
[\ell-1]$. Then by Lemma~\ref{smallcorepart2} we may also assume that the 
inweight of $t_1$ is at least $z' - \gamma n$ and that the outweight of 
$t_\ell$ is at least $y' - \gamma n$. This implies that $z \geq z' - \gamma 
n$ and $y \geq y' - \gamma n$. Since $y'+z'+|T_{\Delta'}| = y+z+|T_\Delta| = 
n$ it follows that we must have  
\begin{equation} \label{eq:yneary'}
y = y' \pm 2 \gamma n \textrm{ and } z = z' \pm 2 \gamma n.
\end{equation}
Finally, by Lemma~\ref{smallcorepart1} we may assume that there are disjoint 
sets $Y, Z \subseteq V(G)$ such that: 
\begin{enumerate}
\item[(a)]$|Y| \geq (2-\gamma)y'$ and $|Z| \geq (2-\gamma)z'$,
\item[(b)] $G[Y]$ and $G[Z]$ are $\gamma$-almost-regular, and
\item[(c)] any vertex of $Y$ has at most $3 \gamma n$ outneighbours in $Z$ 
and any vertex of $Z$ has at most $3 \gamma n$ inneighbours in $Y$.   
\end{enumerate}
Let $X := V(G) \sm (Y \cup Z)$, so $|X| \leq 2 \gamma n$. Let $T^*$ be the 
subtree of $T$ formed by deleting from $T$ all vertices in outcomponents of 
$t_\ell$ or incomponents of $t_1$. So $|T^*| \leq n - (z' - \gamma n) - (y' - 
\gamma n) \leq 3\gamma n$. Let $T^+$ be the subtree of $T$ formed by $t_\ell$ 
and its outcomponents, and let $T^-$ be the subtree of $T$ formed by $t_1$ 
and its incomponents. So $|T^+| \leq y+1$ and $|T^-| \leq z+1$. Also, each 
component of $T^+ - t_\ell$ and each component of $T^- - t_1$ is a component 
of $T - T_\Delta$ and so has order at most $n/\Delta$ by 
Proposition~\ref{coretreeprops}.  

Suppose that some vertex $v \in X$ has at least $\alpha n$ inneighbours in 
$Y$ and at least $\alpha n$ outneighbours in $Z$. Since $\ell \geq 3$, we may 
choose $i$ with $1 < i < \ell$. Embed $t_i$ to $v$. Let $T_a$ be the subtree 
of $T^*$ consisting of $t_i$ and all of its outcomponents, and let $T_b$ be 
the subtree of $T^*$ consisting of $t_i$ and all of its incomponents. Then 
$|T_a|, |T_b| \leq |T^*| \leq 3 \gamma n$. So by Lemma~\ref{onebyone} we may 
extend the embedding of $t_i$ in $\{v\}$ to an embedding of $T_a$ in $Z \cup 
\{v\}$, and similarly we may extend the embedding of $t_i$ in $\{v\}$ to an 
embedding of $T_b$ in $Y \cup \{v\}$. Then in particular~$t_1$ is embedded to 
some $v_1 \in Y$ and $t_\ell$ is embedded to some $v_\ell \in Z$. So $v_1$ 
has at least $|Z| - 3\gamma n \geq z+3\gamma n + 2n/\Delta$ inneighbours in 
$Z$, at most $3 \gamma n$ of which are occupied by the embedding of $T_a$. 
Similarly $v_\ell$ has at least $|Y| - 3\gamma n \geq y+3\gamma n + 
2n/\Delta$ outneighbours in $Y$, at most $3 \gamma n$ of which are occupied 
by the embedding of $T_b$. So by Lemma~\ref{onebyone} we may extend the 
embedding of $t_1$ in $\{v_1\}$ to an embedding of $T^-$ in $\{v_1\} \cup Z$ 
and also extend the embedding of $t_\ell$ in $\{v_\ell\}$ to an embedding of 
$T^+$ in $\{v_\ell\} \cup Y$ so that these embeddings together form a copy of 
$T$ in $G$. 

So we may assume that no vertex of $X$ has at least $\alpha n$ inneighbours 
in $Y$ and at least $\alpha n$ outneighbours in $Z$. Let $X^+ \subseteq X$ 
consist of all vertices of $X$ with fewer than $\alpha n$ inneighbours in 
$Y$, and let $X^- \subseteq X \sm X^+$ consist of all vertices of $X \sm X^+$ 
with fewer than $\alpha n$ outneighbours in $Z$. Let $Y^* := Y \cup X^-$ and 
let $Z^* := Z \cup X^+$, so $Y^*$ and $Z^*$ partition the vertices of $G$. 
Then any vertex of $Y^*$ has at most $\alpha n$ outneighbours in $Z$, and 
thus at least $z + \alpha n$ inneighbours in $Z^*$ (by (a), 
(\ref{eq:yneary'}) and the fact that $z' \geq \eta n$). Similarly any vertex 
of $Z^*$ has at most $\alpha n$ inneighbours in $Y$, and therefore at least 
$y + \alpha n$ outneighbours in $Y^*$. Let $W \subseteq V(G)$ consist of all 
vertices in $Y^*$ with at least $y+\alpha n$ outneighbours in $Y^*$ and all 
vertices in $Z^*$ with at least $z+\alpha n$ inneighbours in $Z^*$.  

Now suppose that $|W| \geq |T_{\Delta}|$. Since $T_{\Delta}$ is a directed 
path, by Theorem~\ref{redeithm} we may embed~$T_{\Delta}$ in~$G[W]$. Let 
$S_\Delta \subseteq W$ be the set of vertices occupied by this embedding. 
Then $|S_\Delta| = |T_{\Delta}| \leq |T_{\Delta'}| \leq \beta n$. So every 
vertex of $S_\Delta$ has at least $y+\alpha n/2 \geq y+2n/\Delta$ 
outneighbours in $Y^* \sm S_\Delta$ and at least $z + \alpha n/2 \geq 
z+2n/\Delta$ inneighbours in $Z^* \sm S_\Delta$. Let $T_1$ be the subtree of 
$T$ consisting of~$T_\Delta$ and all of its outcomponents, and let $T_2$ be 
the subtree of $T$ consisting of~$T_\Delta$ and all of its incomponents. So 
$|T_1| = \ell+y$ and $|T_2| = \ell+z$. Also, each component of $T_1 - 
T_\Delta$ and each component of $T_2 - T_\Delta$ is a component of $T - 
T_\Delta$, and so has order at most $n/\Delta$ by 
Proposition~\ref{coretreeprops}. So by Lemma~\ref{onebyone} we may extend the 
embedding of $T_\Delta$ in $S_\Delta$ to an embedding of $T_1$ in $Y^* \cup 
S_\Delta$. Similarly by Lemma~\ref{onebyone} we may extend the embedding 
of~$T_\Delta$ in $S_\Delta$ to an embedding of $T_2$ in $Z^* \cup S_\Delta$. 
These embeddings of $T_1$ and $T_2$ do not overlap outside $T_\Delta$, and so 
together form an embedding of~$T$ in $G$.  

We may therefore assume that $|W| < |T_{\Delta}|$, and hence that $|G - W| 
\geq 2n-1-\ell$. Since $y+z = n-\ell$, we must have either $|Y^* \sm W| \geq 
2y$ or $|Z^* \sm W| \geq 2z$. Suppose that $|Y^* \sm W| \geq 2y$. Then $Y^* 
\sm W$ contains a vertex $v_\ell$ with at least $y$ outneighbours in $Y^*$. 
So we may choose a set $N \subseteq N^+(v_\ell) \cap Y^*$ with $|N| = y$. 
Then $|N \cap Y| \geq y - (|Y^*| - |Y|) \geq y - 2 \gamma n$. Now, by (a), 
(b) and (\ref{eq:yneary'}) every vertex of $Y$ has at least 
$(1-2\sqrt{\gamma})y$ inneighbours in $Y$ and at least $(1-2\sqrt{\gamma})y$ 
outneighbours in $Y$. Since $|N| = y$, at most $6 \sqrt{\gamma} y$ vertices 
of $N \cap Y$ have more than $(1-3\sqrt{\gamma})y$ inneighbours in $N \cap 
Y$, and at most $6 \sqrt{\gamma} y$ vertices of $N \cap Y$ have more than 
$(1-3\sqrt{\gamma})y$ outneighbours in $N \cap Y$. So at least $|N \cap Y| - 
12 \sqrt{\gamma} n \geq 3n/\Delta$ vertices of $N$ have at least $6n/\Delta$ 
inneighbours in $Y^* \sm (N \cup \{v_\ell\})$ and at least $6n/\Delta$ 
outneighbours in $Y^* \sm (N \cup \{v_\ell\})$. This means that by 
Lemma~\ref{roundtheback} (applied with $Y^* \sm (N \cup \{v_\ell\})$ playing 
the role of $X$) we may embed $T^+$ in $Y^*$ with $t_\ell$ embedded to 
$v_\ell$, and at most $4n/\Delta$ vertices of $T^+$ embedded outside $N$. 
Since $v_\ell \notin W$, $v_\ell$ has at most $y+\alpha n$ outneighbours in 
$Y^*$, and so $v_\ell$ has at least $|Y|-1 - (y+\alpha n) - 4n/\Delta \geq 
9\gamma n$ inneighbours in $Y$ which are not occupied by the embedding of 
$T^+$. Since $|T^*| \leq 3 \gamma n$, by Lemma~\ref{onebyone} we may extend 
the embedding of $t_\ell$ in $v_\ell$ to an embedding of $T^*$ in $Y$ which 
only overlaps the embedding of $T^+$ in $t_\ell$. The vertex $t_1$ of $T$ 
will therefore be embedded to some vertex $v_1 \in Y$. By (3), $v_1$ then has 
at least $|Z| - 3\gamma n \geq z + 2n/\Delta$ inneighbours in $Z$, none of 
which will have been occupied by the embeddings of $T^*$ and $T^+$ so far. So 
by Lemma~\ref{onebyone} we may extend the embedding of $t_1$ in $\{v_1\}$ to 
an embedding of $T^-$ in $Z \cup \{v_1\}$. Then the embeddings of $T^+$, 
$T^-$ and $T^*$ combine to form an embedding of $T$ in~$G$. 
If instead we have $|Z^* \sm W| \geq 2z$, then we may embed $T$ in $G$ 
similarly, first embedding~$T^-$ in $Z^*$, then embedding $T^*$ in the 
unoccupied vertices of $Z$, and finally embedding $T^+$ in~$Y$. So in either 
case $G$ contains a copy of $T$, completing the proof. 
\endproof

\section{Proof of Theorem~\ref{main}} \label{sec:mainproof}

Having proved that Sumner's conjecture holds for directed trees of small core, we now 
show that the same is true for directed trees of large core, which will complete the 
proof of Theorem~\ref{main}. We begin with an embedding result similar to 
Lemma~\ref{smallcorepart1.5}. 

\begin{lemma} \label{largecorepart1}
Suppose that $1/n \ll 1/\Delta \ll \mu \ll \nu \ll \eta \ll \gamma \ll \alpha 
\ll \beta \ll 1$. Let $T$ be a directed tree on $n$ vertices, and let forests~$F^-$ 
and $F^+$ be induced subgraphs of $T$ which partition the vertices of~$T$ 
such that $|F^+|  \geq 6\alpha n$. Suppose also that every edge of~$T$ 
between $F^-$ and $F^+$ is directed from $F^-$ to $F^+$. Let $Y$ and $Z$ be 
disjoint sets with $|Y| \geq 2|F^+| - 2\alpha n \textrm{ and } |Z| \geq 
2|F^-| + \alpha n,$ and let $G$ be a tournament on vertex set $Y \cup Z$ such 
that every vertex of $Y$ has at most $\gamma |G|$ outneighbours in $Z$ and 
every vertex of $Z$ has at most $\gamma |G|$ inneighbours in $Y$. Finally, 
let $T_1^+$ be the largest component of $F^+$, and suppose that either 
\begin{itemize}
\item[(i)] $|T_1^+| \leq |F^+| - 3\alpha n$,
\item[(ii)] $G[Y]$ is a robust $(\mu, \nu)$-outexpander with $\delta^0(G[Y]) 
\geq \eta |Y|$ and $|(T_1^+)_{\Delta}| \geq \beta n$, or 
\item[(iii)] $\Delta(T_1^+) \leq \Delta$.
\end{itemize} 
Then $G$ contains a copy of $T$.
\end{lemma}

\proof
First observe that if $|G| \geq 3n$, then $G$ contains a copy of $T$ by 
Theorem~\ref{bestsofar}. So we may assume that $|G| < 3n$, and hence that 
every vertex of $Y$ has at most $3 \gamma n$ outneighbours in $Z$ and every 
vertex of $Z$ has at most $3 \gamma n$ inneighbours in $Y$. Let $T^+_2$ be 
the second largest component of $F^+$. Then $|F^+| - |T^+_2| \geq |F^+|/2 
\geq 3\alpha n$, so $|Y| \geq |F^+| + |T_2^+| + \alpha n$. Since $|Z| \geq 
2|F^-| + \alpha n$, by Lemma~\ref{component_by_component} any embedding of 
$T^+_1$ in $G[Y]$ may be extended to an embedding of $T$ in $G$. So it is 
sufficient to embed $T^+_1$ in $G[Y]$.  

Note that $|Y| \geq 10\alpha n$, so if $|T_1^+| < \alpha n$, then $G[Y]$ 
contains a copy of $T_1^+$ by Theorem~\ref{bestsofar}. Alternatively, suppose 
that $|T_1^+| \geq \alpha n$. If (i) holds, then $|T_1^+| \leq |Y|/2 - 
2\alpha n$, and so $|Y| \geq (2+\alpha)|T_1^+|$. So $G[Y]$ contains a copy of 
$T_1^+$ by Theorem~\ref{approxversion}(i). If instead (ii) holds then $G$ 
contains a copy of $T_1^+$ by Lemma~\ref{ar_tourns_large_core_tree}. Finally, 
if (iii) holds then $G$ contains a copy of~$T_1^+$ by 
Theorem~\ref{approxversion}(ii), completing the proof. \endproof 

Observe that as with Lemma~\ref{component_by_component} and 
Lemma~\ref{smallcorepart1.5}, a `dual' form of Lemma~\ref{largecorepart1} can 
be proved similarly. For this we instead require that that $|F^-|  \geq 
6\alpha n$, $|Y| \geq 2|F^+| + \alpha n$ and $|Z| \geq 2|F^-| - 2\alpha n$, 
and also either that the largest component $(T_1^-)_{\Delta}$ of $F^-$ 
contains at most $|F^-| - 3\alpha n$ vertices, or that $G[Z]$ is a robust 
$(\mu, \nu)$-outexpander with $\delta^0(G[Z]) \geq \eta |Z|$ and 
$|(T_1^-)_{\Delta}| \geq \beta n$, or that $\Delta(T_1^-) \leq \Delta$. If 
these conditions are met we may conclude that $G$ contains a copy of $T$. As 
with Lemma~\ref{component_by_component}, we shall sometimes implicitly refer 
to this `dual' when referring to Lemma~\ref{largecorepart1}. 

The next lemma is our final result we need to proof Theorem~\ref{main}. It states that if we can find disjoint subsets $Y, Z \subseteq V(G)$ containing almost all of the vertices of $G$, so that $G[Y]$ and $G[Z]$ are robust outexpanders of large minimum semidegree with almost all edges between $Y$ and $Z$ directed the same way, then $G$ contains a copy of $T$.

\begin{lemma} \label{largecorepart2}
Suppose that $1/n \ll 1/\Delta \ll \mu \ll \nu \ll \eta \ll \gamma \ll \alpha 
\ll \beta \ll 1$. Let $T$ be a directed tree on $n$ vertices with $|T_\Delta| \geq 
\beta n$. Let $Y$ and $Z$ be disjoint sets with $|Y \cup Z| \geq 
(2-\alpha)n$, and let $G$ be a tournament on vertex set $Y \cup Z$ such that  
\begin{enumerate}
\item[(i)] $G[Y]$ is a robust $(\mu, \nu)$-outexpander with $\delta^0(G[Y]) 
\geq \eta|Y|$, 
\item[(ii)] $G[Z]$ is a robust $(\mu, \nu)$-outexpander with $\delta^0(G[Z]) 
\geq \eta|Z|$, and 
\item[(iii)] every vertex of $Y$ has at most $\gamma |G|$ outneighbours in 
$Z$, and every vertex of $Z$ has at most $\gamma |G|$ inneighbours in $Y$. 
\end{enumerate}
Then $G$ contains a copy of $T$.
\end{lemma}

\proof
If $|Y \cup Z| \geq (2+\alpha)n$, then $G$ contains a copy of $T$ by 
Theorem~\ref{approxversion}(i). So we may assume that $|Y \cup Z| = 
(2\pm\alpha) n$. Suppose first that $|Z| < 64 \alpha n$. Then $|Y| \geq 
(2-65\alpha) n$, and hence $G[Y]$ contains a copy of $T$ by (i) and 
Lemma~\ref{ar_tourns_large_core_tree}. Similarly if $|Y| < 64 \alpha n$, then 
by (ii) and Lemma~\ref{ar_tourns_large_core_tree} $G[Z]$ contains a copy of 
$T$. So we may assume that $|Y| \geq 64 \alpha n$ and $|Z| \geq 64 \alpha n$. 

So we may form a forest $F^+_1$ of order between $|Y|/2 + 4\alpha n$ and 
$|Y|/2 + 5 \alpha n$ by repeatedly choosing a sink vertex of $T$, deleting it 
from $T$ and adding it to $F^+_1$. Let $F^-_1 := T - F^+_1$, so that 
\begin{equation} \label{eq:sizeofF-}
\frac{|Z|}{2} - 6 \alpha n \leq n-\frac{|Y|}{2} - 5 \alpha n \leq |F_1^-| 
\leq n - \frac{|Y|}{2} - 4 \alpha n \leq \frac{|Z|}{2}- 3 \alpha n. 
\end{equation}
We therefore have $|Y| \geq 2|F_1^+| - 10\alpha n$ and $|Z| \geq 2|F_1^-| + 
6\alpha n$. Note also that $|F_1^+| \geq 36\alpha n$. Let $T'$ be the largest 
component of $F^+_1$. If $|T'| \leq |F_1^+| - 18\alpha n$ or $|T'_\Delta| 
\geq \beta n/3$ then $G$ contains a copy of $T$ by (i), (iii) and 
Lemma~\ref{largecorepart1}. So we may assume that $|T'| > |F^+_1| - 18\alpha 
n$, and that $|T'_\Delta| < \beta n/3$. 

Next we form a forest $F^-_2$ which is a subgraph of $T$ and which contains 
$F^-_1$. To do this, take $F^-_2$ initially to be $F^-_1$. Then select a 
source vertex of $F^+_1$, delete it from $F^+_1$ and add it to $F^-_2$, and 
repeat this step until $|Z|/2 + 4\alpha n \leq |F^-_2| \leq |Z|/2 + 5\alpha 
n$, and let $F^+_2 := T - F^-_2$. Then by (\ref{eq:sizeofF-}) we have $|F^+_1 
\cap F^-_2| = |F^-_2| - |F^-_1| \leq 11 \alpha n$. Also $|F^+_2| \leq |Y|/2 - 
3 \alpha n$, and so we have both $|Z| \geq 2|F_2^-| - 10\alpha n$ and $|Y| 
\geq 2|F_2^+| + 6\alpha n$. Observe also that $|F^-_2| \geq 36 \alpha n$. Let 
$T''$ be the largest component of $F^-_2$. Then if $|T''| \leq |F_2^-| - 18 
\alpha n$ then $G$ contains a copy of $T$ by (ii), (iii) and 
Lemma~\ref{largecorepart1}. So we may assume that $|T''| > |F_2^-| - 18 
\alpha n$. Clearly $|T' \cap T''| \leq |F_1^+ \cap F_2^-| \leq 11 \alpha n$, 
and so $|T' \cup T''| \geq |T'| + |T''| - |T' \cap T''| > (1-47\alpha)n$. 
This implies that $|T''_\Delta| \geq \beta n/3$, as otherwise by 
Lemma~\ref{twocoretrees} we would have $|T_\Delta| < \beta n$, a 
contradiction. Thus $G$ contains a copy of $T$ by (ii), (iii) and 
Lemma~\ref{largecorepart1}, as desired. 
\endproof

\medskip
\noindent {\bf Proof of Theorem~\ref{main}.}
Introduce new constants with 
$$1/n \ll 1/\Delta \ll \mu \ll \nu \ll \eta \ll \gamma \ll \alpha \ll \alpha' 
\ll \beta \ll 1.$$  
If $|T_\Delta| < \beta n$ then $G$ contains a copy of $T$ by 
Lemma~\ref{main_small_core}. So we may assume that $|T_\Delta| \geq \beta n$. 
Let $x := |T_\Delta|$, let $y$ be the outweight of $T_\Delta$, and let $z$ be 
the inweight of $T_\Delta$, so $x+y+z = n$. 
Also let $T_1$ be the subtree of $T$ formed by $T_\Delta$ and all 
outcomponents of $T_\Delta$, and let $T_2$ be the subtree of $T$ formed by 
$T_\Delta$ and all incomponents of~$T_\Delta$, so $|T_1| = x+y$, and $|T_2| = 
x+z$.  

By Lemma~\ref{tournament_split} we may choose disjoint subsets $S_1, \dots, 
S_r$ of $V(G)$ such that 
\begin{itemize}
\item[(i)] $|\bigcup_{i \in [r]} S_i| \geq (1-\gamma) |G|$,
\item[(ii)] for each $i \in [r]$, any vertex $v \in S_i$ has at most $\gamma 
|G|$ inneighbours in $\bigcup_{j > i} S_j$ and at most $\gamma |G|$ 
outneighbours in $\bigcup_{j < i} S_j $, and 
\item[(iii)] for each $i \in [r]$, either $G[S_i]$ is a robust $(\mu, 
\nu)$-outexpander with $\delta^0(G[S_i]) \geq \eta |G|$ or $|S_i| < \gamma 
|G|$. 
\end{itemize}
Let $i$ be maximal such that $|S_1 \cup \dots \cup S_{i-1}| < \max 
\{2(z-\alpha n), 4 \alpha n\}$, and let $j$ be minimal such that $|S_{j+1} 
\cup \dots \cup S_r| < \max \{ 2(y-\alpha n), 4\alpha n\}.$ Since $y+z \leq 
n-\beta n$, by (i) we have $i \leq j$ (though equality is possible here). Let 
$Z := S_1 \cup \dots \cup S_i$, let $Y := S_j \cup \dots \cup S_r$ and let $X 
:= S_{i+1} \cup \dots \cup S_{j-1}$. Then we have 
\begin{equation} \label{eq:7def_of_i_and_j}
|Z \sm S_i| < \max \{2(z-\alpha n), 4 \alpha n\} \textrm{ and } |Y \sm S_j| < 
\max \{ 2(y-\alpha n), 4\alpha n\}. 
\end{equation}
Also, by the maximality of $i$ and the minimality of $j$ we have
\begin{equation} \label{eq:7coro_of_i_and_j}
|Z| \geq  z+\alpha n \textrm{ and } |Y| \geq y + \alpha n.
\end{equation}

\medskip \noindent {\bf Claim.} \emph{If $|Z \sm S_i| \geq 11\alpha n$ or $|Y 
\sm S_j| \geq 11 \alpha n$ then $G$ contains a copy of $T$.} 

\medskip \noindent To prove the claim, suppose first that $|Y \sm S_j| \geq 
11\alpha n$. Let $X^-: = Z \cup X \cup S_j$ and $X^+ := Y \sm S_j$. By 
(\ref{eq:7def_of_i_and_j}) we have $|X^+| < 2y - 2\alpha n$. Also, by (ii) 
every vertex in $X^-$ has at most $\gamma |G|$ inneighbours in $X^+$ and 
every vertex in $X^+$ has at most $\gamma |G|$ outneighbours in $X^-$. Now, 
$T_1 - T_\Delta$ is a forest on $y > |X^+|/2 + \alpha n$ vertices in which 
each component has order at most $n/\Delta$ by 
Proposition~\ref{coretreeprops}(iv). So by repeatedly deleting a source 
vertex of $T_1 - T_\Delta$, we may obtain a subforest~$F^+$ on between 
$|X^+|/2 + 2\alpha n/3$ and $|X^+|/2+\alpha n$ vertices. So $|F^+| \geq 
6\alpha n$, and each component of $F^+$ has order at most $n/\Delta \leq 
|F^+| - 3 \alpha n$. Let $F^+ := T - F^-$, so every edge of $T$ between $F^-$ 
and $F^+$ is directed from $F^-$ to $F^+$. Since $|X^+| + |X^-| \geq 
(1-\gamma)|G|$ by~(i), we have  
$$|F^-| = n - |F^+| \leq n - \frac{|X^+|}{2} - \frac{2\alpha n}{3} \leq 
\frac{|X^-|}{2} - \frac{\alpha n}{2}.$$ 
So $|X^-| \geq 2|F^-| + \alpha n$, and $|X^+| \geq 2|F^+| - 2\alpha n$, and 
so $G$ contains a copy of~$T$ by Lemma~\ref{largecorepart1}(i). If instead 
$|Z \sm S_i| \geq 11\alpha n$ then $G$ contains a copy of $T$ similarly. This 
proves the claim. \medskip 

We may therefore assume that $|Z \sm S_i| < 11 \alpha n$ and $|Y \sm S_j| < 
11 \alpha n$. Suppose first that $i = j$. Then $|S_i| \geq 
(1-\gamma)|G|-22\alpha n \geq (2 - \alpha')n$, so by (iii) $G[S_i]$ is a 
robust $(\mu, \nu)$-outexpander with $\delta^0(G[S_i]) \geq \eta |G| \geq 
\eta |S_i|$. Thus $G$ contains a copy of $T$ by 
Lemma~\ref{ar_tourns_large_core_tree}. Now suppose instead that $i \neq j$, 
and also that $|X| < 12 \alpha' n$. Then $|S_i \cup S_j| \geq (1-\gamma)|G| - 
|X| - 22 \alpha n \geq (2-13\alpha')n$. Now if $|S_i| < \gamma |G|$, then we 
must have $|S_j| \geq (2-14\alpha')n$. Then by (iii) $G[S_j]$ must be a 
robust $(\mu, \nu)$-outexpander with $\delta^0(G[S_j]) \geq \eta |G| \geq 
\eta |S_j|$, so $G[S_j]$ contains a copy of $T$ by 
Lemma~\ref{ar_tourns_large_core_tree}. Alternatively, if $|S_j| < \gamma |G|$ 
then $G[S_i]$ contains a copy of $T$ similarly. Finally, if $|S_i|, |S_j| 
\geq \gamma |G|$, then by (iii) $G[S_i]$ and $G[S_j]$ must both be robust 
$(\mu, \nu)$-outexpanders with $\delta^0(G[S_i]) \geq \eta |G| \geq \eta 
|S_i|$ and $\delta^0(G[S_j]) \geq \eta |S_j|$. Also, by (ii) every vertex of 
$S_i$ has at most $\gamma |G|$ inneighbours in $S_j$, and every vertex of 
$S_j$ has at most $\gamma |G|$ outneighbours in $S_i$. So $G[S_i \cup S_j]$ 
contains a copy of $T$ by Lemma~\ref{largecorepart2}.  

So we may assume that $i \neq j$, and also that $|X| \geq 12 \alpha' n$. We next consider two cases for the size of $X$, in each 
case showing that $T$ may be embedded in $G$. 

{\medskip \noindent \bf Case 1: $|X| \geq (1+\alpha)x$.}

Since by Proposition~\ref{coretreeprops}(iii) we have $\Delta(T_\Delta) \leq 
\Delta$, by Theorem~\ref{approxversion}(ii) we may embed $T_\Delta$ 
in~$G[X]$. Let $X' \subseteq X$ consist of the vertices occupied by this 
embedding. Now, by (ii) every vertex of $X'$ has at most $\gamma |G|$ 
inneighbours in $Y$, and hence by~(\ref{eq:7coro_of_i_and_j}) at least 
$y+\alpha n/2$ outneighbours in $Y$. Since by 
Proposition~\ref{coretreeprops}(iv) every component of $T_1 - T_\Delta$ has 
order at most $n/\Delta$, by Lemma~\ref{onebyone} we may extend the embedding 
of $T_\Delta$ in $G[X']$ to an embedding of~$T_1$ in $G[X' \cup Y]$. 
Similarly by (ii) every vertex of $X'$ has at most $\gamma |G|$ outneighbours 
in~$Z$, and hence by~(\ref{eq:7coro_of_i_and_j}) at least $z+\alpha n/2$ 
inneighbours in $Z$. Since by Proposition~\ref{coretreeprops}(iv) every 
component of $T_2 - T_\Delta$ has order at most $n/\Delta$, by 
Lemma~\ref{onebyone} we may extend the embedding of $T_\Delta$ in $G[X']$ to 
an embedding of $T_2$ in $G[X' \cup Z]$. Since these embeddings of $T_1$ and 
$T_2$ only overlap in $T_\Delta$, they together form an embedding of $T$ in $G$. 

{\medskip \noindent \bf Case 2: $|X| < (1+\alpha)x$.}

Observe that if $|Z| \leq 2z+\alpha n$ and $|Y| \leq 2y+\alpha n$, then by (i)
and the fact that $x = |T_\Delta| \geq \beta n$ we have
$$|X| \geq (1-\gamma)|G| - |Z| - |Y| \geq 2n-2z-2y-3\alpha n\geq 2x - 3\alpha 
n \geq (1+\alpha) x,$$
contradicting our assumption on $X$. So at least one of $|Z| > 2z+\alpha n$ and $|Y| > 2y+\alpha n$
must hold. This gives us three further cases, which we consider separately.

{\medskip \noindent \bf Case 2(a): $|Z| > 2z+\alpha n$, $|Y| \leq 2y+\alpha n$.}

In this case it is sufficient to embed $T_2$ in $G[X \cup Z]$. Indeed, by 
(ii) every vertex of $X \cup Z$ has at most $\gamma |G|$ inneighbours in $Y$, 
and therefore by (\ref{eq:7coro_of_i_and_j}) at least $y+\alpha n/2$ 
outneighbours in $Y$. Since by Proposition~\ref{coretreeprops}(iv) every 
component of $T - T_2$ has order at most $n/\Delta$, any embedding of $T_2$ 
in $G[X \cup Z]$ can be extended to an embedding of $T$ in $G$ by 
Lemma~\ref{onebyone}. 

Now, if $|X \cup Z| \geq 2|T_2| + 2\alpha n$, then we may embed $T_2$ in $G[X 
\cup Z]$ by Theorem~\ref{approxversion}(i). So we may assume that $|X \cup Z| 
< 2|T_2| + 2\alpha n$. Also, by (i) we have 
$$|X \cup Z| \geq (1-\gamma)|G| - |Y| \geq 2n -2y - 2\alpha n = 2x+2z - 
2\alpha n = 2|T_2| - 2 \alpha n.$$ 
So $|X \cup Z| = 2|T_2| \pm 2 \alpha n$. In particular, since $|T_2| \geq 
|T_\Delta| \geq \beta n$, we have $|X \cup Z| \geq \beta n$.  
By repeatedly deleting a source vertex of $T_\Delta$, we may form 
a forest $F$ which is an induced subgraph of $T_\Delta$ (consisting of the undeleted
vertices of $T_\Delta$) so that every edge 
between $T_\Delta - F$ and $F$ is directed from $T_\Delta$ to $F$, and also 
so that 
$$\frac{|X|}{2} + \frac{2\alpha' |T_2|}{3} \leq |F| \leq \frac{|X|}{2} + 
\alpha' |T_2|.$$ 
Let $F^- := T_2 - F$. Then
$$|F^-| = |T_2| - |F| \leq |T_2| - \frac{|X|}{2} - \frac{2\alpha' |T_2|}{3} 
\leq \frac{|Z|}{2} - \frac{\alpha' |T_2|}{2}.$$ 
So $|X| \geq 2|F| - 2 \alpha'  |T_2|$ and $|Z| \geq 2|F^-| + \alpha'  |T_2|$. 
Also, $|F| \geq |X|/2 \geq 6 \alpha' |T_2|$, and since~$F$ is a subtree of 
$T_\Delta$, by Proposition~\ref{coretreeprops}(iii) each component $C$ of $F$ 
has $\Delta(C) \leq \Delta$. Since by~(ii) every vertex of $X$ has at most 
$\gamma |G| \leq 2\gamma |X \cup Z|/\beta$ outneighbours in $Z$ and every 
vertex of~$Z$ has at most $\gamma |G|  \leq 2\gamma |X \cup Z|/ \beta$ 
inneighbours in $X$, $G[X \cup Z]$ contains a copy of $T_2$ by 
Lemma~\ref{largecorepart1}, as required.  

{\medskip \noindent \bf Case 2(b): $|Z| \leq 2z+\alpha n$, $|Y| > 2y+\alpha n$.} 

In this case $T$ may be embedded in $G$ by the same method as in the previous 
case, with the roles of inneighbours and outneighbours switched. So we begin 
by embedding $T_1$ in $G[X \cup Y]$, and then use Lemma~\ref{onebyone} to 
extend this embedding to an embedding of $T$ in $G$. 

{\medskip \noindent \bf Case 2(c): $|Z| > 2z+\alpha n$, $|Y| > 2y+\alpha n$.}

In this case, we shall partition $T$ into three forests as follows. Initially 
take $F^-$ to be the forest formed by all incomponents of $T_\Delta$, and 
$F^+$ to be the forest formed by all outcomponents of~$T_\Delta$. Then select 
a source vertex of $T_\Delta$, delete it from $T_\Delta$ and add it to $F^-$. 
Repeat this step until $2|F^-| +\alpha n \leq |Z| \leq 2|F^-| + 2\alpha n$. 
Next, select a sink vertex of $T_\Delta$, delete it from~$T_\Delta$ and add 
it to $F^+$. Repeat this step until $2|F^+| +\alpha n \leq |Y| \leq 2|F^+| + 
2\alpha n$. Then let $F$ consist of all vertices remaining in~$T_\Delta$. So 
$F$ is a subgraph of $T_\Delta$. Also, by (i)  
$$|F| = n - |F^-| - |F^+| \leq n - |Y|/2 - |Z|/2 + 2\alpha n \leq |X|/2 + 
3\alpha n,$$ 
so (since $|X| \geq \alpha' n$) $|X| \geq |F| + \alpha n$.
We shall embed the components of $F^-$, $F$ and $F^+$ in turn amongst the 
vertices of $Z$, $X$ and $Y$ respectively. Indeed, the proof is similar to 
the proof of Lemma~\ref{component_by_component}, but with three forests 
instead of two. 

Let $C_1, \dots, C_s$ be the components of $F^-$, $F$ and $F^+$, ordered so 
that $C_1$ is a component of $F$, and for each $i \in [s-1]$, $C_{i+1}$ has 
precisely one neighbour in $C_1 \cup \dots \cup C_i$. We shall embed the 
$C_i$ in turn, so that each component of $F^-$ is embedded in $G[Z]$, each 
component of $F$ is embedded in $G[X]$, and each component of $F^+$ is 
embedded in $G[Y]$. We also require that after each $C_i$ is embedded, the 
embeddings of $C_1, \dots, C_i$ together form an embedding in $G$ of the 
subtree of $T$ induced by the vertices of $C_1, \dots, C_i$. So suppose that 
we have successfully embedded $C_1, \dots, C_{i-1}$ in this manner, and we 
now wish to extend this embedding to include $C_i$. Then if $i \geq 2$, there 
is precisely one edge of $T$ between $C_i$ and $C_1 \cup \dots \cup C_{i-1}$. 
Let $t$ be the endvertex of this edge in $C_1 \cup \dots \cup C_{i-1}$, and 
let $v$ be the vertex to which $t$ was embedded. If $C_i$ is a component of 
$F^-$, then $i \geq 2$, the edge between $t$ and $C_i$ is directed towards $t$
and $v \in X \cup Y$. So we may let $S$ consist of the inneighbours of $v$ 
in~$Z$. Then by (ii) we have $|S| \geq |Z| - \gamma |G|$. Let $S' \subseteq 
S$ consist of the unoccupied vertices of $S$. Since at most $|F^-| - |C_i|$ 
vertices of $S$ are occupied by the embeddings of $C_1, \dots, C_{i-1}$,  
$$|S'| \geq |Z| - \gamma |G| - |F^-| + |C_i| \geq 2|C_i| + \alpha n/2.$$
So if $|C_i| < \alpha n/2$ then $G[S']$ contains a copy of $T$ by 
Theorem~\ref{bestsofar}, and if $|C_i| \geq \alpha n/2$ then $G[S']$ contains a 
copy of $T$ by Theorem~\ref{approxversion}(i). Alternatively, if $C_i$ is a 
component of $F^+$, then $i \geq 2$, the edge between $t$ and $C_i$ is directed
towards $C_i$ and $v \in X \cup Z$. So we may let $S$ consist of the 
outneighbours of $v$ in $Y$, and let $S' \subseteq S$ consist of the 
unoccupied vertices of $S$. Then we may embed $C_i$ in $S'$ by the same 
argument as used when $C_i$ is a component of $F^-$. Finally, suppose that 
$C_i$ is a component of $F$. Then if $i \geq 2$ and $t \in F^+$, let $S$ 
consist of the inneighbours of $v$ in $X$. If instead $i \geq 2$ and $t \in 
F^-$, let $S$ consist of the outneighbours of $v$ in~$X$. If $i=1$ then let 
$S = X$. Then by (ii) we have $|S| \geq |X| - \gamma |G|$. Again let $S' 
\subseteq S$ consist of the unoccupied vertices of $S$. Then it suffices to 
embed $C_i$ in $G[S']$. Since at most $|F| - |C_i|$ vertices have been 
embedded in $S$, we have $|S'| \geq  |X| - \gamma |G| - |F| + |C_i| \geq 
|C_i| + \alpha n/2$. Now,~$C_i$ is a subtree of $T_\Delta$, so $\Delta(C_i) 
\leq \Delta$ by Proposition~\ref{coretreeprops}(iii). So if $|C_i| \geq \alpha n/4$,
then $G[S']$ contains a copy of $C_i$ by Theorem~\ref{approxversion}(ii).
On the other hand, if $|C_i| < \alpha n/4$, then $G[S']$ contains a copy of $C_i$
by Theorem~\ref{bestsofar}. So in any case we may embed~$C_i$ as desired, completing the proof. 
\endproof

\medskip

\noindent
{\footnotesize
\noindent
Daniela K\"uhn, Deryk Osthus, \\School of Mathematics,
\\University of Birmingham, \\Birmingham, \\B15 2TT, \\United Kingdom, \\
\{{\tt kuehn,osthus}\}{\tt @maths.bham.ac.uk } }

\medskip
\noindent
{\footnotesize
\noindent
Richard Mycroft, \\ School of Mathematical Sciences,
\\ Queen Mary, University of London, \\London, \\E1 4NS, \\United Kingdom, \\
{\tt r.mycroft@qmul.ac.uk } }
\end{document}